\documentclass[lettersize, 11pt]{article}
\usepackage{graphicx, bbm, latexsym, amsfonts, epsfig, verbatim, amsmath, color, epsfig, verbatim, enumerate, multirow, caption, graphicx, subcaption, algorithm, algpseudocode, amssymb, hyperref}
\usepackage{url, multicol, color, latexsym, amsfonts, epsfig, verbatim, tikz, cases, bbm}
\usepackage[title]{appendix}
\usepackage[top=1.0in, bottom=1.0in, left=1.0in, right=1.0in]{geometry}
\usepackage{algcompatible}
\usepackage{diagbox,multicol,graphicx}
\usepackage{booktabs, caption, makecell}

\usepackage{threeparttable}
\usepackage{array}
\newcolumntype{L}[1]{>{\raggedright\let\newline\\\arraybackslash\hspace{0pt}}m{#1}}
\newcolumntype{C}[1]{>{\centering\let\newline\\\arraybackslash\hspace{0pt}}m{#1}}
\newcolumntype{R}[1]{>{\raggedleft\let\newline\\\arraybackslash\hspace{0pt}}m{#1}}

\newtheorem{theorem}{Theorem}
\newtheorem{lemma}{Lemma}

\newtheorem{assumption}{Assumption}

\newtheorem{proposition}{Proposition}

\newtheorem{example}{Example}

\newtheorem{structure}{Structure}
\newenvironment{customstructure}[1]
  {\renewcommand\thestructure{#1}\structure}
  {\endstructure}

\def\tinyl{\mbox{\tiny L}}
\def\tinyu{\mbox{\tiny U}}
\newcommand{\proof}{\noindent{\em Proof: }}
\newcommand{\qed}{\hfill$\Box$ \newline\noindent}

\allowdisplaybreaks
\begin{document}
\title{ {\bf Nurse Staffing under Absenteeism: \mbox{A Distributionally Robust Optimization Approach}}
\author{\normalsize \bf{Minseok Ryu}, \bf{Ruiwei Jiang} \\
{\small Department of Industrial and Operations Engineering}\\
{\small University of Michigan, Ann Arbor, MI 48109}\\[1mm]
{\small Email: \{msryu,ruiwei\}@umich.edu}\\
}
}
\date{ }
\maketitle

\begin{abstract}
\setlength{\baselineskip}{18pt}
We study the nurse staffing problem under random nurse demand and absenteeism. While the demand uncertainty is exogenous (stemming from the random patient census), the absenteeism uncertainty is \emph{endogenous}, i.e., the number of nurses who show up for work partially depends on the nurse staffing level. 
{\color{black}
For quality of care, many hospitals have developed float pools, i.e., groups of hospital units, and trained nurses to be able to work in multiple units (termed cross-training) in response to potential nurse shortage.}
In this paper, we propose a distributionally robust nurse staffing (DRNS) model that considers both exogenous and endogenous uncertainties. We derive a separation algorithm to solve this model under a general structure of float pools. In addition, we identify several pool structures that often arise in practice and recast the corresponding DRNS model as {\color{black} a mixed-integer linear program}, which facilitates off-the-shelf commercial solvers. Furthermore, we optimize the float pool design to {\color{black} reduce cross-training while achieving specified target staffing costs}. The numerical case studies, based on the data of a collaborating hospital, suggest that the units with high absenteeism probability should be pooled together.

\mbox{}

\end{abstract}
%

\setlength{\baselineskip}{18pt}

\section{Introduction} \label{sec:intro}
Nurse staffing plays a key role in hospital management. The cost of staffing nurses accounts for over 30\% of the overall hospital annual expenditures (see, e.g.,~\cite{welton2011hospital}). Besides, {\color{black} the nurse staffing level can have significant impact} on patient safety, quality of care, and the job satisfaction of nurses (see, e.g.,~\cite{tevington2011mandatory}). In view of that, a number of governing agencies (e.g., the California Department of Health~\cite{california2004california} and the Victoria Department of Health~\cite{victoria2015victoria}) have set up minimum nurse-to-patient ratios (NPRs) for various types of hospital units to regulate the staffing decision.

In general, {\color{black} nurse planning consists of} the following four phases: (A) nurse demand forecasting and staffing, (B) nurse shift scheduling, (C) pre-shift staffing and re-scheduling, and (D) nurse-patient assignment (see~\cite{bard2010nurse, hulshof2012taxonomic, kortbeek2015flexible, bam2017planning}). In particular, phase (A) takes place weeks or months ahead of a shift and determines the nurse staffing levels based on, e.g., the forecasted patient census and the NPRs; and phase (C) takes place hours before the shift and recruits additional workforce (e.g., temporary or off-duty nurses) if any units are short of nurses. In this paper, we focus on these two phases and {\color{black} refer to the corresponding decision making process as ``nurse staffing,'' i.e., decide the staffing level of nurses for multiple units during a specific shift.} 
The outputs of our study (e.g., the nurse staffing levels) can be used in phases (B) and (D) to generate shift schedules and assignments of individual nurses.

Nurse staffing is a challenging task, largely because of the uncertainties of nurse demand and absenteeism. The demand uncertainty stems from the random patient census and has been well documented (see, e.g.,~\cite{gaynor1995uncertain,boutsioli2010demand}) and studied in the nurse staffing literature (see, e.g.,~\cite{davis2014nurse,kim2015two}). In contrast, the absenteeism uncertainty has received relatively less attention in this literature (see, e.g.,~\cite{green2013nursevendor,maass2017incorporating}), albeit commonly observed in practice. For example, according to the U.S. Bureau of Labor Statistics~\cite{Bureau2016}, the average absence rate among all nurses in the Veterans Affairs Health Care System is 6.4\%~\cite{wang2014nurse}, significantly higher than that among all occupations (2.9\%) and among health-care support occupations (4.3\%). For quality of care, many hospitals have developed float pools, which consist of multiple hospital units. A nurse assigned to a float pool is trained to be able to work in multiple units, termed cross-training, so that in phase (C) the nurse can be assigned to any unit in the pool that is short of nurses (see, e.g.,~\cite{inman2005cross}).

Unlike demand, the random number of nurses who show up for a shift partially depends on the nurse staffing level, i.e., the absenteeism uncertainty is \emph{endogenous}. 
{\color{black} For example, if the nurse staffing level is $w \in \mathbb{N}_+$, then the nurse absenteeism rate depends on $w$ as discussed in \cite{green2013nursevendor}, and the random number of nurses who show up cannot exceed $w$.}
Although failing to incorporate such endogeneity may result in understaffing (see~\cite{green2013nursevendor}), unfortunately, modeling endogeneity usually makes optimization models computationally prohibitive (see, e.g.,~\cite{dupacova2006optimization}). Due to this technical difficulty, endogenous uncertainty has received much less attention in the literature of stochastic optimization than exogenous uncertainty. Existing works often resort to exogenous uncertainty for an approximate solution. Alternatively, they employ certain parametric probability distributions to model the endogenous uncertainty (see~\cite{dupacova2006optimization}), e.g., the absence of each nurse follows \emph{independent} Bernoulli distribution with the \emph{same} probability (which may depend on the staffing level; see~\cite{green2013nursevendor}). A basic challenge to adopting parametric models is that a complete and accurate knowledge of the endogenous probability distribution is usually unavailable. Under many circumstances, we only have historical data, e.g., past nurse staffing levels and absence records, which can be considered as samples taken from the true (but ambiguous) endogenous distribution. As a result, the solution obtained by assuming a parametric model can yield unpleasant out-of-sample performance if the chosen model is biased.

In this paper, we propose an alternative, nonparametric model of both exogenous and endogenous uncertainties based on distributionally robust optimization (DRO). Our approach considers a family of probability distributions, termed an ambiguity set, based only on the support and moment information of these uncertainties. In particular, the number of nurses who show up in a unit/pool is bounded by the corresponding staffing level and its mean value is a function of this level. Then, we employ this ambiguity set in a two-stage distributionally robust nurse staffing (DRNS) model that imitates the decision making process in phases (A) and (C). Building on DRNS, we further search for \emph{sparse} pool structures that result in a minimum amount of cross-training while achieving a specified target staffing cost. To the best of our knowledge, this is the first study of the endogenous uncertainty in nurse staffing by using a DRO approach.

\subsection{Literature Review} \label{sec:literature-review}
A vast majority of the nurse staffing literature focuses on deterministic models that do not take into account the randomness of the nurse demand and/or absenteeism (see~\cite{van2013personnel}). Various (deterministic) optimization models have been employed, including linear programming (see, e.g.,~\cite{kao1981aggregate,brusco1993constrained}) and mixed-integer programming (see, e.g.,~\cite{trivedi1981mixed,yang2007evaluation,maenhout2013integrated,bam2017planning}). For example,~\cite{kao1981aggregate} assessed the need for hiring permanent staffs and temporary helpers and~\cite{trivedi1981mixed} analyzed the trade-offs among hiring full-time, part-time, and overtime nurses. More recently,~\cite{yang2007evaluation} compared cross-training and flexible work days and demonstrated that cross-training is far more effective for performance improvement than flexible work days. Similarly,~\cite{bam2017planning} identified cross-training as a promising extension from their deterministic model. Despite the potential benefit of operational flexibility brought by float pools and cross-training,~\cite{maenhout2013integrated} pointed out that the pool design and staffing are often made manually in a qualitative fashion (also see~\cite{stewart1994mathematical}). In addition, when the nurse demand and/or absenteeism is random, the deterministic models may underestimate the total staffing cost (see, e.g.,~\cite{kao1985budgeting}).

Existing stochastic nurse staffing models often consider the demand uncertainty only. For example,~\cite{campbell2011two} studied a two-stage stochastic programming model that integrates the staffing and scheduling of cross-trained workers (e.g., nurses) under demand uncertainty. Through numerical tests,~\cite{campbell2011two} demonstrated that cross-training can be even more valuable than the perfect demand information (i.e., knowing the realization of demand when making staffing decisions). In addition,~\cite{lu2016mandatory} studied how the mandatory overtime laws can negatively affect the service quality of a nursing home. Using a two-stage stochastic programming model under demand uncertainty,~\cite{lu2016mandatory} pointed out that these laws result in a lower staffing level of permanent registered nurses and a higher staffing level of temporary registered nurses. Unfortunately, as~\cite{green2013nursevendor} pointed out, ignoring nurse absenteeism may result in understaffing, which reduces the service quality and increases the operational cost because additional temporary nurses need to be called in.

When the nurse absenteeism is taken into account, the stochastic optimization models become unscalable.~\cite{green2013nursevendor} considered the staffing of a single unit under both nurse demand and absenteeism uncertainty and successfully derived a closed-form optimal staffing level. In addition,~\cite{slaugh2018consistent} studied the staffing of a single on-call pool that serves multiple units whose staffing levels are fixed and known. In a setting that regular nurses can be absent while pool nurses always show up, the authors successfully derived a closed-form optimal pool staffing level. Unfortunately, the problem becomes computationally prohibitive when multiple units and/or multiple float pools are incorporated. For example,~\cite{easton2014service} studied a multi-unit and one-pool setting\footnote{More precisely, the model in~\cite{easton2014service} allows to re-assign nurses from one unit to any other unit. In the context of this paper, that is equivalent to having a single float pool that serves all the units and assigning all nurses to this pool.}. The author showed that the proposed stochastic optimization model outperforms the (deterministic) mean value approximation.
However, the evaluation of this model ``does not scale well.'' More specifically, even when staffing levels are \emph{fixed}, one needs to solve an exponential number (in terms of the staffing level) of linear programs to evaluate the expected total cost of staffing. This renders the search of an optimal staffing level so challenging that one has to resort to heuristics.~\cite{wang2014nurse} considered a multi-unit and no-pool setting and analyzed the staffing problem based on a cohort of nurses who have heterogeneous absence rates. The authors showed that the staffing cost is lower when the nurses are heterogenous within each unit but uniform across units. Unfortunately, searching for an optimal staffing strategy is ``computationally demanding'' with a large number of nurses. Similar to~\cite{easton2014service},~\cite{wang2014nurse} resort to easy-to-use heuristics.

To mitigate the computational challenges of nurse absenteeism, the existing literature often make parametric assumptions on the endogenous probability distribution. For example,~\cite{green2013nursevendor,easton2014service,wang2014nurse} assumed that the absences of all nurses are stochastically \emph{independent} and the absence rate in~\cite{green2013nursevendor,easton2014service} is assumed \emph{homogeneous}. But the nurse absences may be positively correlated during extreme weather (e.g., heavy snow) or during day shifts (e.g., due to conflicting family obligations). In addition, {\color{black} the data analysis in~\cite{wang2014nurse}} suggests that the nurses actually have heterogeneous absence rates. Furthermore, the absenteeism can be drastically different among different units/hospitals, and even within the same unit/hospital, has high temporal variations. For example, based on the data from different hospitals,~\cite{green2013nursevendor} concluded that the absence rate depends on the staffing level and ignoring such dependency results in understaffing, while~\cite{wang2014nurse} concluded that such dependency is insignificant. A fundamental challenge to adopting parametric models is that the solution thus obtained can yield suboptimal out-of-sample performance if the adopted model is biased. In this paper, we take into account both nurse demand and absenteeism uncertainty in a multi-unit and multi-pool setting. To address the challenges on computational scalability and out-of-sample performance, we propose an alternative nonparametric model based on DRO. In particular, this model allows dependence or independence between the absence rate and the staffing level. Moreover, our model can be solved to global optimality by a separation algorithm and, in several important special cases, by solving a single mixed-integer linear program (MILP).

DRO models have received increasing attention in the recent literature. In particular, DRO has been applied to two-stage stochastic optimization (see, e.g.,~\cite{bertsimas2010models,bertsimas2018adaptive,hanasusanto2018conic}), as in this paper. In general, two-stage DRO models are computationally prohibitive. For example, suppose that the second-stage formulation is linear and continuous with right-hand side uncertainty. Then, even with \emph{fixed} first-stage decision variables,~\cite{bertsimas2010models} showed that evaluating the objective function of the DRO model is NP-hard. To mitigate the computational challenge,~\cite{hanasusanto2018conic,kong2013scheduling} recast the two-stage DRO model as a copositive program, which admits semidefinite programming approximations. In addition,~\cite{ang2014two,bertsimas2018adaptive} applied linear decision rules (LDRs) to obtain conservative and tractable approximations. In contrast to these work, our second-stage formulation involves integer variables to model the pre-shift staffing. Besides undermining the convexity of our formulation, this prevents us from applying the LDRs because fractional staffing levels are not implementable. To the best of our knowledge, there are only two existing work~\cite{luo2018distributionally,noyan2018distributionally} on DRO with endogenous uncertainty. Specifically,~\cite{luo2018distributionally} derived equivalent reformulations of the endogenous DRO model under various ambiguity sets, and~\cite{noyan2018distributionally} applied an endogenous DRO model on the machine scheduling problem. In this paper, we study a two-stage endogenous DRO model for nurse staffing and derive tractable reformulations under several practical float pool structures. We summarize our main contributions as follows:
\begin{enumerate}[1.]
\item We propose the first DRO approach for nurse staffing, considering both exogenous nurse demand and endogenous nurse absenteeism. The proposed two-stage endogenous DRO model considers multiple units, multiple float pools, and both long-term and pre-shift nurse staffing. For general pool structures, we derive a min-max reformulation of the model and a separation algorithm that solves this model to global optimality.
\item For multiple pool structures that often arise in practice, including one pool, disjoint pools, and chained pools, we provide {\color{black}a MILP} reformulation of our DRO model by deriving strong valid inequalities. The binary variables of this MILP reformulation arise from the nurse staffing decisions only. That is, under these practical pool structures, the computational burden of our DRO approach is de facto the same as that of deterministic nurse staffing.
\item Building upon the DRO model, we further study how to design sparse and effective disjoint pools. To this end, we proactively optimize the nurse pool structure to minimize the total number of cross-training, while providing a guarantee on the staffing cost.
\item We conduct extensive case studies based on the data and insights from our collaborating hospital. The results demonstrate the value of modeling nurse absenteeism and the computational efficacy of our DRO approach. In addition, we provide managerial insights on how to design sparse and effective pools.
\end{enumerate}
The remainder of the paper is organized as follows. In Section \ref{sec: Model}, we describe the two-stage DRO model with endogenous nurse absenteeism. In Section \ref{sec: Solution Approach}, we derive a solution approach for this model under general pool structures. In Section \ref{sec:tractable}, we derive strong valid inequalities and tractable reformulations under special pool structures. We extend the DRO model for optimal pool design in Section \ref{sec:pool-design}, conduct case studies in Section \ref{sec:results}, and conclude in Section \ref{sec:conclusion}. To ease the exposition, we relegate all proofs to the appendices.

\noindent {\bf Notation}: We use \textsuperscript{$\sim$} to indicate random variables and \textsuperscript{$\wedge$} to indicate realizations of the random variables. For example, $\tilde{d}$ represents a random variable and $\hat{d}^1, \ldots, \hat{d}^N$ represent $N$ realizations of $\tilde{d}$. For $a, b \in \mathbb{Z}$, we define $[a] := \{1, 2, \ldots, a\}$ and $[a, b]_{\mathbb{Z}} := \{n \in \mathbb{Z}: a \leq n \leq b\}$. For $x \in \mathbb{R}$, we define $[x]_+ := \max\{x, 0\}$. For set $S$, we define its indicator function $\mathbbm{1}_S$ such that $\mathbbm{1}_S(s) := 1$ if $s \in S$ and $\mathbbm{1}_S(s) := 0$ if $s \notin S$, and denote its convex hull by $\mbox{conv}(S)$.
{\color{black} We define $\mathbb{B} := \{0,1\}$. }

\section{Distributionally Robust Nurse Staffing} \label{sec: Model}
We consider a group of $J$ hospital units, each facing a random demand of nurses denoted by $\tilde{d}_j$ for all $j \in [J]$ {\color{black} during a specific shift.} To enhance the operational flexibility, the manager forms $I$ nurse float pools. For all $i \in [I]$, pool $i$ is associated with a set $P_i$ of units and each nurse assigned to this pool is capable of working in any unit $j \in P_i$. Due to random absenteeism, if we staff unit $j$ with $w_j$ nurses (termed unit nurses) {\color{black} during the shift}, then there will be a random number $\tilde{w}_j$ of nurses showing up for work, where $\tilde{w}_j \in [0, w_j]_{\mathbb{Z}}$. 
Likewise, $\tilde{y}_i$ nurses show up if we staff pool $i$ with $y_i$ nurses {\color{black} during the shift}, where $\tilde{y}_i \in [0, y_i]_{\mathbb{Z}}$. {\color{black}As specified in the ambiguity set $\mathcal{D}(w, y)$ defined in \eqref{ambiguity-mean}, we assume that the probability distribution $\mathbb{P}_{\tilde{w}, \tilde{y}, \tilde{d}}$ of $(\tilde{w}, \tilde{y}, \tilde{d})$ depends on decision variables $(w, y)$. In particular, the expected number of nurses who show up in a unit/pool is a piecewise linear function of the staffing level.} After the uncertain parameters $\tilde{d}_j$, $\tilde{w}_j$, and $\tilde{y}_i$ are realized, the nurses showing up in pool $i$ can be re-assigned to any units in $P_i$ to make up the nurse shortage, if any. After the re-assignment, any remaining shortage will be covered by hiring temporary nurses in order to meet the NPR requirement. Mathematically, for given $\tilde{w} := [\tilde{w}_1, \ldots, \tilde{w}_J]^{\top}$, $\tilde{y} := [\tilde{y}_1, \ldots, \tilde{y}_I]^{\top}$, and $\tilde{d} := [\tilde{d}_1, \ldots, \tilde{d}_J]^{\top}$, the total operational cost can be obtained from solving the following integer program:
\begin{subequations}
\label{sns-sub}	
\begin{align}
V(\tilde{w}, \tilde{y}, \tilde{d}) \ := \ \min_{z, x, e} \ & \ {\color{black} \sum_{j=1}^J  c^{\mbox{\tiny x}}_j x_j } \label{sns-sub-obj} \\
\mbox{s.t.} \ & \ \sum_{i \in [I]:\ j \in P_i} z_{ij} + x_j - e_j = \tilde{d}_j - \tilde{w}_j, \ \ \forall j \in [J], \label{sns-sub-con-demand} \\
& \ \sum_{j \in P_i} z_{ij} \leq \tilde{y}_i, \ \ \forall i \in [I], \label{sns-sub-con-pool} \\
& \ x_j, e_j \in \mathbb{Z}_+, \ \ \forall j \in [J], \ \ z_{ij} \in \mathbb{Z}_+, \ \ \forall i \in [I], \ \forall j \in P_i, \label{sns-sub-con-integer}
\end{align}
\end{subequations}
where variables $z_{ij}$ represent the number of nurses re-assigned from pool $i$ to unit $j$, variables $x_j$ represent the number of temporary nurses hired in unit $j$, variables $e_j$ represent the excessive number of nurses in unit $j$, {\color{black} and parameters $c^{\mbox{\tiny x}}_j$ represent the unit cost of hiring temporary nurses in unit $j$.}
In the above formulation, objective function \eqref{sns-sub-obj} minimizes the cost of hiring temporary nurses. Constraints \eqref{sns-sub-con-demand} describe three ways of satisfying the nurse demand in each unit: (i) assigning unit nurses, (ii) re-assigning pool nurses, and (iii) hiring temporary nurses. Constraints \eqref{sns-sub-con-pool} ensure that the number of nurses re-assigned from each pool does not exceed the number of nurses showing up in that pool. Constraints \eqref{sns-sub-con-integer} describe integrality restrictions.
{\color{black} In formulation \eqref{sns-sub}, we assume that the resource of temporary nurses is uncapacitated  because hospitals often have sufficient time to recruit temporary nurses either from the previous shift or from collaborating agencies. In Appendix \ref{apx-capacitated_temp_nurses}, we discuss a variant of formulation \eqref{sns-sub} that relaxes this assumption.}

\begin{figure}[H]
	\centering
	\includegraphics[scale=0.35]{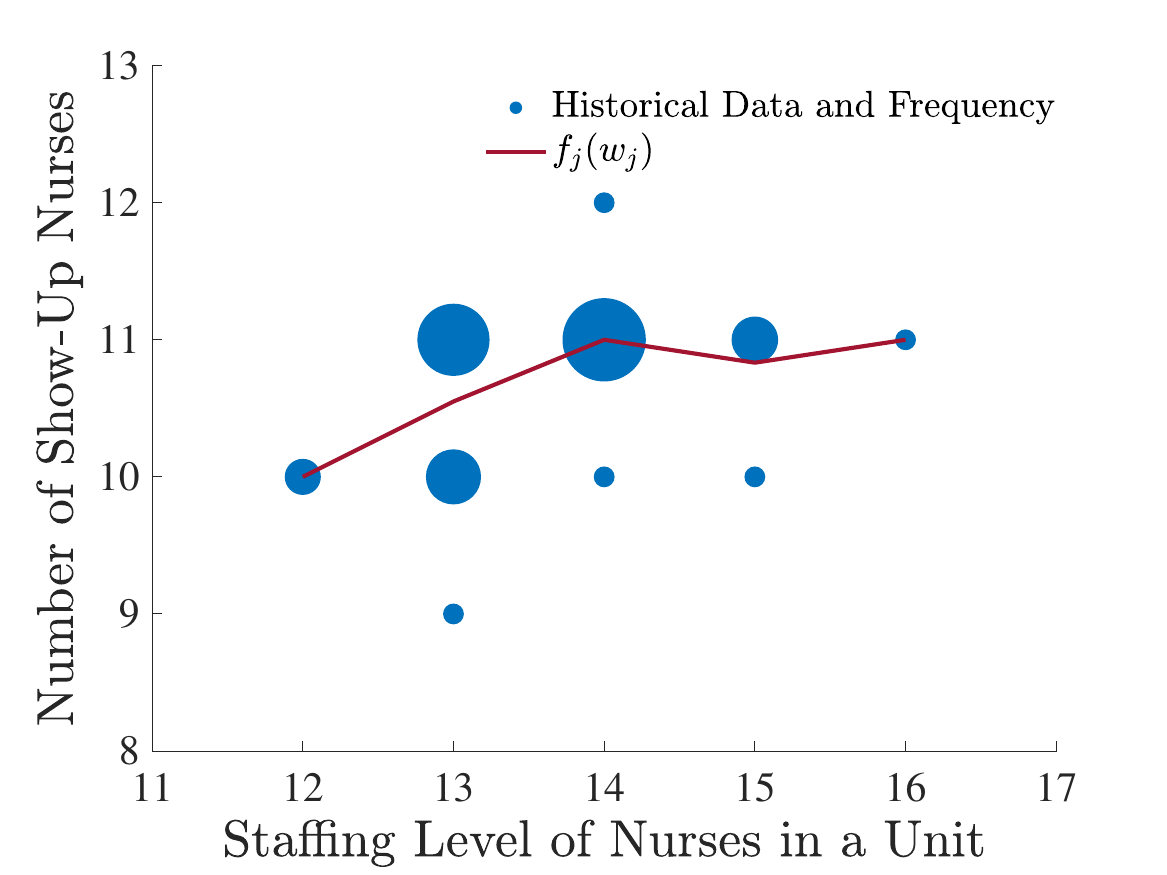}	
	\includegraphics[scale=0.35]{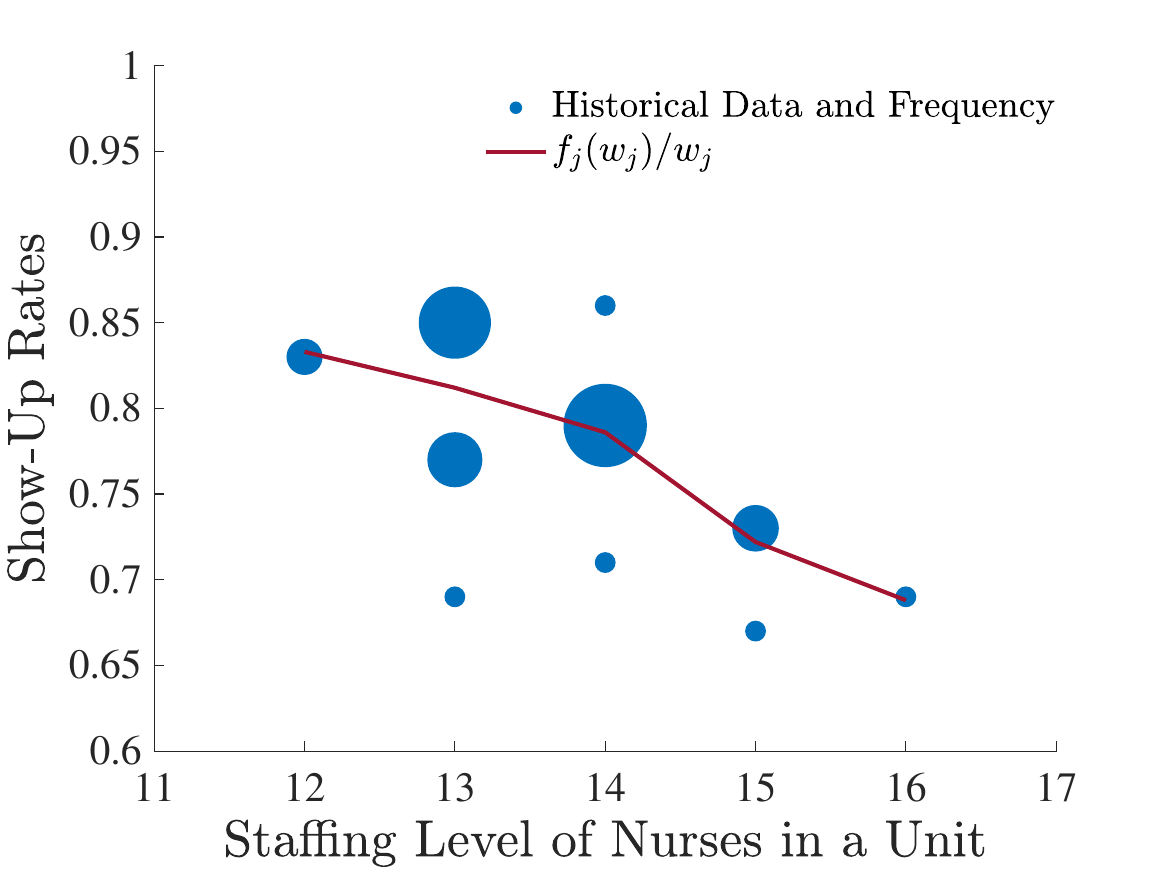}
	\caption{{\color{black}Number of show-up nurses and show-up rates, as functions of staffing levels, based on the data from a collaborating hospital. Dots represent historical data samples and the size of dots represent frequency.}} \label{Fig:segmented}
\end{figure}

In reality, it is often challenging to obtain an accurate estimate of the true probability distribution $\mathbb{P}_{\tilde{w}, \tilde{y}, \tilde{d}}$ of $(\tilde{w}, \tilde{y}, \tilde{d})$. For example, the historical data of the nurse demand (via patient census and NPRs) can typically be explained by multiple (drastically) different distributions. More importantly, because of the endogeneity of $\tilde{w}$ and $\tilde{y}$, $\mathbb{P}_{\tilde{w}, \tilde{y}, \tilde{d}}$ is in fact a \emph{conditional distribution} depending on the nurse staffing levels. This further increases the difficulty of estimation. {\color{black} We emphasize that not only the number of show-up nurses but also the show-up rates depend on the staffing level (see, e.g.,~\cite{green2013nursevendor} and Figure~\ref{Fig:segmented}).} Using a biased estimate of $\mathbb{P}_{\tilde{w}, \tilde{y}, \tilde{d}}$ can yield post-decision disappointment. For example, if one simply ignores the endogeneity of $\tilde{w}$ and $\tilde{y}$ and employs their empirical distribution based on historical data, then the nurse staffing thus obtained may lead to disappointing out-of-sample performance. In this paper, we assume that $\mathbb{P}_{\tilde{w}, \tilde{y}, \tilde{d}}$ is ambiguous and it belongs to the following {\color{black}endogenous moment} ambiguity set:
\begin{subequations}
\label{ambiguity-mean}
\begin{align}
{\color{black} \mathcal{D} (w,y) } \ = \ \Bigl\{\mathbb{P} \in \mathcal{P}(\Xi): \ & \ \mathbb{E}_{\mathbb{P}}[\tilde{d}_j^q] = \mu_{jq}, \ \ \forall j \in [J], \ \forall q \in [Q], \label{demand-moments} \\
& \ \mathbb{E}_{\mathbb{P}}[\tilde{w}_j] = f_j(w_j), \ \ \forall j \in [J], \label{nurse-linear-unit} \\
& \ \mathbb{E}_{\mathbb{P}}[\tilde{y}_i] = g_i(y_i), \ \ \forall i \in [I] \Bigr\}, \label{nurse-linear-pool}
\end{align}
\end{subequations}
where $\Xi$ represents the support of $(\tilde{w}, \tilde{y}, \tilde{d})$ and $\mathcal{P}(\Xi)$ represents the set of {\color{black} probability distributions} supported on $\Xi$. We consider a box support $\Xi := \Xi_{\tilde{w}} \times \Xi_{\tilde{y}} \times \Xi_{\tilde{d}}$, where $\Xi_{\tilde{w}} = \Pi_{j=1}^J [0, w_j]_{\mathbb{Z}}$, $\Xi_{\tilde{y}} = \Pi_{i=1}^I [0, y_i]_{\mathbb{Z}}$, $\Xi_{\tilde{d}} = \Pi_{j=1}^J [d^{\tinyl}_j, d^{\tinyu}_j]_{\mathbb{Z}}$, and $d^{\tinyl}_j$ and $d^{\tinyu}_j$ represent lower and upper bounds of the nurse demand in unit $j$. In addition, for $Q \in \mathbb{N}_+$, all $q \in [Q]$, and all $j \in [J]$, $\mu_{jq}$ represents the $q^{\mbox{\tiny th}}$ moment of $\tilde{d}_j$. Furthermore, for all $j \in [J]$ and $i \in [I]$, $f_j: \mathbb{N}_+ \rightarrow \mathbb{R}_+$ and $g_i: \mathbb{N}_+ \rightarrow \mathbb{R}_+$ represent two {\color{black}piecewise linear} functions such that $f_j(0) = g_i(0) = 0$. We note that these functions can model arbitrary dependence of $(\tilde{w}, \tilde{y})$ on the staffing levels. {\color{black} This is because the staffing levels $(w, y)$ are integer-valued. As a consequence, any $f_j$ and $g_i$ can be equivalently re-expressed as piecewise linear functions. For example, function $f_j(w_j)$ can be re-expressed as $\hat{f}_j(w_j) := f_j(k) + \left[f_j(k+1)-f_j(k)\right](w_j - k)$ if $k \leq w_j \leq k+1$ for $k \in [0, w^{\mbox{\tiny U}}_j]_{\mathbb{Z}}$.} In addition, the assumption $f_j(0) = g_i(0) = 0$ ensures that if we assign no nurses in a unit/pool then nobody will show up.

The ambiguity set {\color{black} $\mathcal{D}(w,y)$} can be conveniently calibrated. First, suppose that $\mathbb{P}_{\tilde{w}, \tilde{y}, \tilde{d}}$ is observed through nurse demand data $\{\hat{d}^1_j, \ldots, \hat{d}^N_j\}_{j=1}^J$ and attendance records $\{(w^1_j, \hat{w}^1_j), \ldots, (w^N_j, \hat{w}^N_j)\}_{j=1}^J$ and $\{(y^1_i, \hat{y}^1_j), \ldots, (y^N_i, \hat{y}^N_i)\}_{i=1}^I$ during the past $N$ days, where, in each pair $(w^n_j, \hat{w}^n_j)$, $w^n_j$ represents the staffing level of unit $j$ in day $n$ and $\hat{w}^n_j$ represents the corresponding number of nurses who actually showed up. Then, $\mu_{jq}$ can be obtained from empirical estimates (e.g., $\mu_{j1} = (1/N)\sum_{n=1}^N \hat{d}^n_j$, $\mu_{j2} = (1/N)\sum_{n=1}^N (\hat{d}^n_j)^2$, etc.), and $f_j$ and $g_i$ can be obtained by performing segmented linear regression on the attendance data, using the staffing levels $\{w^1_j, \ldots, w^N_j\}$ and $\{y^1_i, \ldots, y^N_i\}$ as breakpoints, respectively (see Figure~\ref{Fig:segmented} for an example). Second, if $\tilde{w}$ and $\tilde{y}$ are believed to follow certain parametric models, then we can follow such models to calibrate $\{f_j(w_j)\}_{j=1}^J$ and $\{g_i(y_i)\}_{i=1}^I$. For example, if $\tilde{w}_j$ is modeled as a Binomial random variable $B(w_j, 1 - a(w_j))$ as in~\cite{green2013nursevendor}, where $a(w_j)$ represents the absence rate, i.e., the probability of any scheduled nurse in unit $j$ being absent from work, then we have $f_j(w_j) = w_j(1 - a(w_j))$. {\color{black}Finally, although it appears intuitive that the nurse demand (i.e., anticipated workload) may influence the show-up rates, this linkage can be statistically insignificant (see, e.g., [41]). In this regard, $\mathcal{D}(w, y)$ is robust against the potential linkage by allowing an arbitrary correlation between $\tilde{d}$ and $(\tilde{w}, \tilde{y})$.}

We seek nurse staffing levels that minimize the expected total cost with regard to the worst-case probability distribution in {\color{black}$\mathcal{D}(w,y)$}, i.e., we consider the following two-stage DRO model:
\begin{subequations} \label{drns}
\begin{align}
(\mbox{\bf DRNS}): \ \ \ \ \min_{w, y} \ & \ {\color{black} \sum_{j=1}^J c^{\mbox{\tiny w}}_j w_j + \sum_{i=1}^I c^{\mbox{\tiny y}}_i y_i + \sup_{\mathbb{P} \in \mathcal{D}(w,y)}\mathbb{E}_{\mathbb{P}}\left[ V(\tilde{w}, \tilde{y}, \tilde{d}) \right] } \label{drns-obj} \\
\mbox{s.t.} \ & \ w^{\tinyl}_j \leq w_j \leq w^{\tinyu}_j, \ \ \forall j \in [J], \label{sns-con-w-bd} \\
& \ y^{\tinyl}_i \leq y_i \leq y^{\tinyu}_i, \ \ \forall i \in [I], \label{sns-con-y-bd} \\
& \ y, w \in R \cap \mathbb{Z}^{I+J}_+, \label{sns-con-integer}
\end{align}
\end{subequations}
where parameters {\color{black} $c^{\mbox{\tiny w}}_j$ and $c^{\mbox{\tiny y}}_i$ } represent the unit cost of hiring unit and pool nurses, respectively, constraints \eqref{sns-con-w-bd}--\eqref{sns-con-y-bd} designate lower and upper bounds on staffing levels, and set $R$ represents all remaining restrictions, which we assume can be represented via mixed-integer linear inequalities. {\color{black}We remark that (DRNS) decides the staffing levels during a specific shift and does not decide the schedule for individual nurses across multiple shifts. Hence, multi-stage models do not apply to our problem.} (DRNS) is computationally challenging because (i) {\color{black} $\mathcal{D}(w,y)$} involves exponentially many probability distributions, all of which depend on the decision variables $w_j$ and $y_i$ and (ii) it is a two-stage DRO model with integer recourse variables. In the next two sections, we shall derive equivalent reformulations of (DRNS) that facilitate a separation algorithm, and identify practical pool structures that admit more tractable solution approaches.

\section{Solution Approach: General Pool Structure} \label{sec: Solution Approach}
In this section, we consider general pool structures, recast (DRNS) as a min-max formulation, and derive a separation algorithm that solves this model to global optimality.

We start by noticing that the integrality restrictions \eqref{sns-sub-con-integer} in the second-stage formulation of (DRNS) can be relaxed without loss of generality.
\begin{lemma} \label{prop:sns-sub-tu}
For any given $(\tilde{w}, \tilde{y}, \tilde{d}) \in \Xi$, the value of $V(\tilde{w}, \tilde{y}, \tilde{d})$ remains unchanged if constraints \eqref{sns-sub-con-integer} are replaced by non-negativity restrictions, i.e., $x_j, e_j \geq 0, \forall j \in [J]$ and $z_{ij} \geq 0, \forall i \in [I], \ \forall j \in P_i$.
\end{lemma}
\proof See Appendix \ref{apx-prop:sns-sub-tu}.  \qed
Thanks to Lemma \ref{prop:sns-sub-tu}, we are able to rewrite $V(\tilde{w}, \tilde{y}, \tilde{d})$ as the following dual formulation:
\begin{subequations}
\begin{align}
V(\tilde{w}, \tilde{y}, \tilde{d}) \ = \ \max_{\alpha, \beta} \ & \ \sum_{j=1}^J (\tilde{d}_j - \tilde{w}_j) \alpha_j + \sum_{i=1}^I \tilde{y}_i \beta_i \label{sns-sub-dual-obj} \\
\mbox{s.t.} \ & \ \beta_i + \alpha_j \leq 0, \ \ \forall i \in [I], \ \forall j \in P_i, \label{sns-sub-dual-con-max} \\
& \ {\color{black} 0 \leq \alpha_j \leq c^{\mbox{\tiny x}}_j, \ \ \forall j \in [J],} \label{sns-sub-dual-con-bd}
\end{align}
\end{subequations}
where dual variables $\alpha_j$ and $\beta_i$ are associated with primal constraints \eqref{sns-sub-con-demand} and \eqref{sns-sub-con-pool}, respectively, and dual constraints \eqref{sns-sub-dual-con-max} and \eqref{sns-sub-dual-con-bd} are associated with primal variables $z_{ij}$ and $(x_j, e_j)$, respectively. We let $\Lambda$ denote the dual feasible region for variables $(\alpha, \beta)$, i.e., $\Lambda := \{(\alpha, \beta): \mbox{\eqref{sns-sub-dual-con-max}--\eqref{sns-sub-dual-con-bd}}\}$. Strong duality between formulations \eqref{sns-sub-obj}--\eqref{sns-sub-con-integer} and \eqref{sns-sub-dual-obj}--\eqref{sns-sub-dual-con-bd} holds valid because \eqref{sns-sub-obj}--\eqref{sns-sub-con-integer} has a finite optimal value.

We are now ready to recast (DRNS) as a min-max formulation. To this end, we consider $\mathbb{P}$ as a decision variable and take the dual of the worst-case expectation in \eqref{drns-obj}. For strong duality, we make the following technical assumption on the ambiguity set {\color{black}$\mathcal{D}(w,y)$}.
\begin{assumption} \label{assump:technical}
For any given $w := [w_1, \ldots, w_J]^{\top}$ and $y := [y_1, \ldots, y_I]^{\top}$ that are feasible to (DRNS), {\color{black}$\mathcal{D}(w,y)$} is non-empty.
\end{assumption}
Assumption \ref{assump:technical} is mild. For example, it holds valid whenever the moments of demands $\{\mu_{jq}: j \in [J], q \in [Q]\}$ are obtained from empirical estimates and the {\color{black} decision-dependent} moments $\{g_i(y_i), f_j(w_j): i \in [I], j \in [J]\}$ lie in the convex hull of their support, i.e., $f_j(w_j) \in [0, w_j]$ and $g_i(y_i) \in [0, y_i]$. In Appendix \ref{apx-prop:non-empty}, we present an approach to verify Assumption \ref{assump:technical} by solving $J$ linear programs. The reformulation is summarized in the following proposition.
\begin{proposition} \label{prop:ref}
Under Assumption \ref{assump:technical}, the (DRNS) model \eqref{drns} yields the same optimal value and the same set of optimal solutions as the following min-max optimization problem:
\begin{subequations}
\begin{align}
\min_{\substack{w, y\\ \gamma, \lambda, \rho}} \ & \ \max_{(\alpha, \beta) \in \Lambda} F(\alpha, \beta) + \sum_{i=1}^I ({\color{black} c^{\mbox{\tiny y}}_i } y_i + g_i(y_i) \lambda_i) + \sum_{j=1}^J \left[{\color{black} c^{\mbox{\tiny w}}_j} w_j + \sum_{q=1}^Q \mu_{jq} \rho_{jq} + f_j(w_j) \gamma_j \right] \label{ref-note-8} \\
\mbox{s.t.} \ & \ \mbox{\eqref{sns-con-w-bd}--\eqref{sns-con-integer}}, \label{ref-note-con}
\end{align}
where
\begin{equation}
F(\alpha, \beta) \ := \ \sum_{j=1}^J \Bigl[(- \alpha_j - \gamma_j) w_j\Bigr]_+ + \sum_{i=1}^I \Bigl[ (\beta_i - \lambda_i) y_i \Bigr]_+ + \sum_{j=1}^J \sup_{\tilde{d}_j \in [d^{\tinyl}_j, d^{\tinyu}_j]_{\mathbb{Z}}} \Bigl\{\alpha_j \tilde{d}_j - \sum_{q=1}^Q\rho_{jq}\tilde{d}^q_j \Bigr\}. \label{ref-f-def}
\end{equation}
\end{subequations}
\end{proposition}
\proof See Appendix \ref{apx-prop:ref}. \qed
In the min-max reformulation \eqref{ref-note-8}--\eqref{ref-note-con}, the additional variables $\gamma$, $\lambda$, $\rho$ are generated in the process of taking dual. In addition, function $F(\alpha, \beta)$ is jointly convex in $(\alpha, \beta)$ because, as presented in \eqref{ref-f-def}, $F(\alpha, \beta)$ is the pointwise maximum of functions affine in $(\alpha, \beta)$. This min-max reformulation is not directly computable because (i) for fixed $(w, y, \gamma, \lambda, \rho)$, evaluating the objective function \eqref{ref-note-8} needs to solve a \emph{convex maximization} problem $\max_{(\alpha, \beta) \in \Lambda} F(\alpha, \beta)$, which is in general NP-hard, and (ii) the formulation includes nonlinear and non-convex terms $g_i(y_i) \lambda_i$ and $f_j(w_j) \gamma_j$. We shall address these two challenges before presenting a separation algorithm for solving (DRNS).

First, we analyze the convex maximization problem $\max_{(\alpha, \beta) \in \Lambda} F(\alpha, \beta)$ and derive the following optimality conditions.
{\color{black}
\begin{lemma} \label{thm:extreme-points}
Without loss of generality, we assume that $0 \leq c^{\text{\tiny x}}_1 \leq \cdots \leq c^{\text{\tiny x}}_J$.	
For fixed $(w, y, \gamma, \lambda, \rho)$, there exists an optimal solution $(\bar{\alpha}, \bar{\beta})$ to problem $\max_{(\alpha, \beta) \in \Lambda} F(\alpha, \beta)$ such that (a) $\bar{\alpha}_j \in \{0, c^{\mbox{\tiny x}}_1, \cdots, c^{\mbox{\tiny x}}_j \}$ for all $j \in [J]$ and (b) $\bar{\beta}_i = - \max\{ \bar{\alpha}_j: j \in P_i\}$ for all $i \in [I]$.
\end{lemma}
\proof See Appendix \ref{apx-thm:extreme-points}. \qed

\noindent\textbf{Remark 1 (Interpretation on $\bar{\alpha}$ and $\bar{\beta}$)} $\bar{\alpha}$ are the shadow prices associated with constraints (1b), i.e., $\bar{\alpha}_j$ evaluates the increase in staffing cost when the nurse shortage $(\tilde{d}_j-\tilde{w}_j)$ of unit $j$ increases by $1$. When this happens, a new optimal solution to formulation~\eqref{sns-sub} re-assigns a pool nurse from a unit $k \in \{1, \ldots, j\}$ to unit $j$ due to the ascending staffing costs $c^{\text{\tiny x}}_1 \leq \cdots \leq c^{\text{\tiny x}}_j$. As a result, it is equivalent to transferring the additional shortage from unit $j$ to unit $k$ and accordingly the staffing cost increases by $\bar{\alpha}_j \in \{0, c^{\text{\tiny x}}_1, \ldots, c^{\text{\tiny x}}_j\}$. In addition, $\bar{\beta}_i$ evaluates the decrease in staffing cost when an additional nurse shows up in pool $i$. When that happens, a new optimal solution to \eqref{sns-sub} assigns this nurse to a unit in which the shortage is most costly. As a result, $-\bar{\beta}_i = \max\{\bar{\alpha}_j: j \in P_i\}$.

}
Lemma \ref{thm:extreme-points} enables us to avoid enumerating the infinite number of elements in $\Lambda$ and focus only on a finite set of $(\bar{\alpha}, \bar{\beta})$ values. In addition, we introduce binary variables to encode the special structure identified in the optimality conditions.
{\color{black} 
Specifically, for all $j \in [J]$ and $k \in [j]$, we define binary variable $t_{jk}$ such that $t_{jk} := 1$ if $\bar{\alpha}_j = c^{\mbox{\tiny x}}_k$, and $t_{jk} := 0$ otherwise.
For all $i \in [I]$, we denote  by $J(i)$ the largest index in $P_i$, i.e., $J(i) = \max\{j: j \in P_i\}$. For all $i \in [I]$ and $k \in [J(i)]$, we define binary variable $s_{ik} := 1$ if $\bar{\beta}_i = - c^{\mbox{\tiny x}}_k$ and $s_{ik} := 0$ otherwise. In other words, $s_{ik} = 1$ if $k$ is the largest index in $[J(i)]$ with $t_{jk}=1$ for some $j \in P_i$, in line with condition (b) in Lemma \ref{thm:extreme-points}.}
Variables $(t,s)$ need to satisfy the following constraints to make the encoding well-defined:
{\color{black} 
\begin{subequations}
	\begin{align}
	& \sum_{k = 1}^j t_{jk} \leq 1, \ \forall j \in [J], \label{encode-0} \\
	& \sum_{k=1}^{J(i)} s_{ik} \leq 1, \ \forall i \in [I], \label{encode-1}\\
	& s_{ik} \leq \sum_{j \in P_i: j \geq k} t_{jk}, \ \forall i \in [I], \forall k \in [J(i)], \label{encode-2} \\
	& s_{ik} + t_{j \ell} \leq 1, \ \forall i \in [I], \forall k  \in [J(i)], \forall j \in P_i, \forall \ell \in [k+1,j]_{\mathbb{Z}}, \label{encode-3} \\
	& t_{jk} \leq \sum_{\ell = 1}^{J(i)} s_{i \ell}, \ \forall i \in [I], \forall j \in P_i, \forall k \in [j], \label{encode-4}\\
	& s_{ik} \in \mathbb{B}, \ \forall i \in [I], \forall k \in [J(i)], \ \ t_{jk} \in \mathbb{B}, \ \forall j \in [J], \forall k \in [j],	\label{encode-5} 
	\end{align}	
\end{subequations}
where 
constraints \eqref{encode-0} describe that, for all $j \in [J]$, $t_{jk}=1$ holds for at most one $k \in [j]$; 
constraints \eqref{encode-1} describe that, for all $i \in [I]$, $s_{ik}=1$ holds for at most one $k \in [J(i)]$; 
constraints \eqref{encode-2} designate that if $s_{ik}=1$, then there is a $j \in P_i$ with $t_{jk}=1$; 
constraints \eqref{encode-3} describe that if $s_{ik} = 1$, then there cannot be any $j \in P_i$ and $\ell \geq k+1$ with $t_{j\ell} = 1$; 
and constraints \eqref{encode-4} ensure that all $t_{jk} = 0$ if $s_{i \ell}=0$ for all $\ell \in [J(i)]$.
}
This recasts $\max_{(\alpha, \beta) \in \Lambda} F(\alpha, \beta)$ as an integer linear program presented in the following theorem. For ease of exposition, we introduce auxiliary variables $r_j := 1 - \sum_{k=1}^j t_{jk}$ and $p_i := 1 - \sum_{k=1}^{J(i)} s_{ik}$, whose values depend entirely on variables $t$ and $s$.

{\color{black}
\begin{theorem} \label{thm:ip}
For fixed $(w, y, \gamma, \lambda, \rho)$, problem $\max_{(\alpha, \beta) \in \Lambda} F(\alpha, \beta)$ yields the same optimal value as the following integer linear program:
\begin{subequations}
	\label{IntegerProgram}
	\begin{align}
	\max_{t,s,r,p} \ & \sum_{j=1}^J \bigg( c ^{\text{\tiny r}}_j r_j + \sum_{k=1}^j c^{\text{\tiny t}}_{jk} t_{jk} \bigg) + \sum_{i=1}^I \bigg( c^{\text{\tiny p}}_i p_i + \sum_{k=1}^{J(i)} c^{\text{\tiny s}}_{ik} s_{ik} \bigg) \label{IntegerProgram-0} \\
	\mbox{s.t.} \
	& (t,s,r,p) \in \mathcal{H} : = \Big\{  \mbox{\eqref{encode-0}--\eqref{encode-5}}, \ r_j + \sum_{k=1}^j t_{jk} = 1, \ \forall j \in [J], \label{IntegerProgram-1} \\
	& \hspace{33mm}  p_i + \sum_{k=1}^{J(i)} s_{ik} = 1, \ \forall i \in [I] \label{IntegerProgram-3}
	\Big\},
	\end{align}
	where
$c^{\text{\tiny r}}_j := [ - \gamma_j w_j ]_+ + \sup_{\tilde{d}_j \in [d^{\tinyl}_j, d^{\tinyu}_j]_{\mathbb{Z}}} \bigg\{ - \sum_{q=1}^Q \rho_{jq} \tilde{d}^q_j \bigg\}$, 
$c^{\text{\tiny t}}_{jk} :=[ (-c^{\text{\tiny x}}_k - \gamma_j ) w_j ]_+ + \sup_{\tilde{d}_j \in [d^{\tinyl}_j, d^{\tinyu}_j]_{\mathbb{Z}}} \bigg\{ c^{\text{\tiny x}}_k \tilde{d}_j - \sum_{q=1}^Q \rho_{jq} \tilde{d}^q_j \bigg\}$, 
$c^{\text{\tiny p}}_i := [ - \lambda_i y_i  ]_+$, and
$c^{\text{\tiny s}}_{ik} := [ (-c^{\text{\tiny x}}_k - \lambda_i ) y_i  ]_+$. 
\end{subequations}
\end{theorem}
\proof See Appendix \ref{apx-thm:ip}. \qed 
}
\indent Second, we linearize the terms $f_j(w_j) \gamma_j$ and $g_i(y_i) \lambda_i$. For all $j \in [J]$, although $f_j(w_j)$ can be nonlinear and non-convex, thanks to the integrality of $w_j$, we can rewrite $f_j(w_j)$ as an affine function based on a binary expansion of $w_j$. 
Specifically, we introduce binary variables $\{u_{j \ell}: \ell \in [w^{\mbox{\tinyu}}_j - w^{\mbox{\tinyl}}_j]\}$ such that $w_j =  w^{\mbox{\tinyl}}_j + \sum_{\ell=1}^{w^{\mbox{\tinyu}}_j - w^{\mbox{\tinyl}}_j} u_{j\ell}$, where we interpret $u_{j\ell}$ as whether we assign at least $w^{\mbox{\tinyl}}_j + \ell$ nurses to unit $j$. That is, $u_{j\ell} = 1$ if $w_j \geq w^{\mbox{\tinyl}}_j + \ell$ and $u_{j\ell} = 0$ otherwise. Then, defining $\Delta_{j\ell} := f_j(w^{\tinyl}_j + \ell) - f_j(w^{\tinyl}_j + \ell - 1)$ for all $\ell \in [w^{\mbox{\tinyu}}_j - w^{\mbox{\tinyl}}_j]$, we have
\begin{align*}
f_j(w_j) \ = & \ f_j(w^{\tinyl}_j) + \sum_{\ell=1}^{w_j - w^{\mbox{\tinyl}}_j} \left[ f_j(w^{\tinyl}_j + \ell) - f_j(w^{\tinyl}_j + \ell - 1) \right] \\
= & \ f_j(w^{\tinyl}_j) + \sum_{\ell=1}^{w^{\mbox{\tinyu}}_j - w^{\mbox{\tinyl}}_j} \left[ f_j(w^{\tinyl}_j + \ell) - f_j(w^{\tinyl}_j + \ell - 1) \right] \mathbbm{1}_{[w^{\mbox{\tinyl}}_j + \ell, w^{\mbox{\tinyu}}_j]}(w_j) \\
= & \ f_j(w^{\tinyl}_j) + \sum_{\ell=1}^{w^{\mbox{\tinyu}}_j - w^{\mbox{\tinyl}}_j} \Delta_{j\ell} u_{j\ell}.
\end{align*}
It follows that $f_j(w_j) \gamma_j = f_j(w^{\tinyl}_j) \gamma_j + \sum_{\ell=1}^{w^{\mbox{\tinyu}}_j - w^{\mbox{\tinyl}}_j} \Delta_{j\ell} u_{j\ell} \gamma_j$. We can linearize the bilinear terms $u_{j \ell} \gamma_j$ by defining continuous variables $\varphi_{j \ell} := u_{j \ell} \gamma_j$ and incorporating the following standard McCormick inequalities (see~\cite{mccormick1976computability}):
\begin{subequations}
\begin{align}
\gamma_j - M(1 - u_{j\ell}) \leq \varphi_{j\ell} \leq M u_{j\ell}, \ \ \forall j \in [J], \ \forall \ell \in [w^{\mbox{\tinyu}}_j - w^{\mbox{\tinyl}}_j], \label{ref-note-36} \\
- M u_{j\ell} \leq \varphi_{j\ell} \leq \gamma_j + M(1 - u_{j\ell}), \ \ \forall j \in [J], \ \forall \ell \in [w^{\mbox{\tinyu}}_j - w^{\mbox{\tinyl}}_j], \label{ref-note-37}
\end{align}
where $M$ represents a sufficiently large positive constant. Likewise, for all $i \in [I]$, we rewrite $g_i(y_i) \lambda_i$ as $g_i(y^{\tinyl}_i) \lambda_i + \sum_{\ell=1}^{y^{\mbox{\tinyu}}_i - y^{\mbox{\tinyl}}_i} \delta_{i\ell} v_{i\ell} \lambda_i$ by using constants $\delta_{i\ell} := g_i(y^{\tinyl}_i + \ell) - g_i(y^{\tinyl}_i + \ell - 1)$ for all $\ell \in [y^{\mbox{\tinyu}}_i - y^{\mbox{\tinyl}}_i]$ and binary variables $\{v_{i\ell}: \ell \in [y^{\mbox{\tinyu}}_i - y^{\mbox{\tinyl}}_i]\}$, where $v_{i\ell} = 1$ if $y_i \geq y^{\mbox{\tinyl}}_i + \ell$ and $v_{i\ell} = 0$ otherwise. We linearize the bilinear terms $v_{i\ell} \lambda_i$ by continuous variables $\nu_{i\ell} := v_{i\ell} \lambda_i$ and the McCormick inequalities
\begin{align}
\lambda_i - M(1 - v_{i\ell}) \leq \nu_{i\ell} \leq M v_{i\ell}, \ \ \forall i \in [I], \ \forall \ell \in [y^{\mbox{\tinyu}}_i - y^{\mbox{\tinyl}}_i], \label{ref-note-40} \\
- M v_{i\ell} \leq \nu_{i\ell} \leq \lambda_i + M(1 - v_{i\ell}), \ \ \forall i \in [I], \ \forall \ell \in [y^{\mbox{\tinyu}}_i - y^{\mbox{\tinyl}}_i]. \label{ref-note-41}
\end{align}
In computation, a large big-M coefficient $M$ can significantly slow down the solution of (DRNS). Theoretically, for the correctness of the linearization \eqref{ref-note-36}--\eqref{ref-note-41}, $M$ needs to be larger than $|\gamma_j|$ and $|\lambda_i|$ for all $j \in [J]$ and $i \in [I]$, respectively. The following proposition derives uniform lower and upper bounds of $\gamma_j$ and $\lambda_i$, leading to a small value of $M$.
\begin{proposition} \label{prop:big-M}
For fixed $w$ and $y$, there exists an optimal solution $(\gamma^*, \lambda^*, \rho^*)$ to formulation \eqref{ref-note-8}--\eqref{ref-note-con} such that $\gamma^*_j \in [-c^{\mbox{\tiny x}}_j, 0]$ for all $j \in [J]$ and {\color{black} $\lambda^*_i \in [-c^{\mbox{\tiny x}}_{J(i)}, 0]$ } for all $i \in [I]$.
\end{proposition}
\proof See Appendix \ref{apx-prop:big-M}. \qed 
Proposition \ref{prop:big-M} indicates that (i) we can set {\color{black}$M := c^{\mbox{\tiny x}}_{J}$} in the McCormick inequalities \eqref{ref-note-36}--\eqref{ref-note-41} without loss of optimality and (ii) as all $\gamma_j$ and $\lambda_i$ are non-positive at optimality, we can replace McCormick inequalities \eqref{ref-note-36} and \eqref{ref-note-40} as $\gamma_j \leq \varphi_{j\ell} \leq 0$ and $\lambda_i \leq \nu_{i\ell} \leq 0$ respectively, both of which are now big-M-free. In addition, we incorporate the following constraints to break the symmetry among binary variables:
\begin{align}
& u_{j\ell} \geq u_{j(\ell+1)}, \ \ \forall j \in [J], \ \forall \ell \in [w^{\mbox{\tinyu}}_j - w^{\mbox{\tinyl}}_j - 1], \label{ref-note-32} \\
& v_{i\ell} \geq v_{i(\ell+1)}, \ \ \forall i \in [I], \ \forall \ell \in [y^{\mbox{\tinyu}}_i - y^{\mbox{\tinyl}}_i - 1]. \label{ref-note-33}
\end{align}
\end{subequations}

The above analysis recasts (DRNS) into a mixed-integer program, which is summarized in the following theorem without proof.
{\color{black}
\begin{theorem} \label{thm:ref-linearized}
Under Assumption \ref{assump:technical}, the (DRNS) model \eqref{drns} yields the same optimal value as the following mixed-integer program:
\begin{subequations}
\label{drns-milp-exp}	
\begin{align}
\min_{\substack{u, v, \varphi, \nu\\ \gamma, \lambda, \rho, \theta}} \ & \ \theta + \sum_{i=1}^I \Bigl(c^{\mbox{\tiny y}}_i y^{\tinyl}_i + g_i(y^{\tinyl}_i) \lambda_i + \sum_{\ell=1}^{y^{\mbox{\tinyu}}_i - y^{\mbox{\tinyl}}_i} \bigl(c^{\mbox{\tiny y}}_i v_{i\ell} + \delta_{i\ell} \nu_{i\ell}  \bigr) \Bigr) \nonumber \\
& \ + \sum_{j=1}^J \Biggl[ \sum_{q=1}^Q \mu_{jq} \rho_{jq} + c^{\mbox{\tiny w}}_j w^{\tinyl}_j + f_j(w^{\tinyl}_j) \gamma_j + \sum_{\ell=1}^{w^{\mbox{\tinyu}}_j - w^{\mbox{\tinyl}}_j} \bigl( c^{\mbox{\tiny w}}_j u_{j \ell} + \Delta_{j \ell} \varphi_{j \ell} \bigr) \Biggr] \label{ref-linear-1} \\
\mbox{s.t.} \ & \ \mbox{\eqref{ref-note-36}--\eqref{ref-note-33}}, \label{ref-linear-2} \\
& \ \theta \geq \sum_{j=1}^J \bigg( c ^{\text{\tiny r}}_j r_j + \sum_{k =1}^j c^{\text{\tiny t}}_{jk} t_{jk} \bigg) + \sum_{i=1}^I \bigg( c^{\text{\tiny p}}_i p_i + \sum_{k = 1}^{J(i)} c^{\text{\tiny s}}_{ik} s_{ik} \bigg), \ \ \forall (t, s, r, p) \in \mathcal{H}, \label{ref-linear-3} \\
& \ u_{j\ell} \in \mathbb{B}, \ \ \forall j \in [J], \ \forall \ell \in [w^{\mbox{\tinyu}}_j - w^{\mbox{\tinyl}}_j], \ \ v_{i\ell} \in \mathbb{B}, \ \ \forall i \in [I], \ \forall \ell \in [y^{\mbox{\tinyu}}_i - y^{\mbox{\tinyl}}_i], \label{ref-linear-4}
\end{align}
where set $\mathcal{H}$ is defined in \eqref{IntegerProgram-1}--\eqref{IntegerProgram-3} and coefficients $c^{\mbox{\tiny r}}_j$, $c^{\mbox{\tiny t}}_{jk}$, $c^{\mbox{\tiny p}}_i$, and $c^{\mbox{\tiny s}}_{ik}$ are represented through
\begin{align}
c^{\text{\tiny r}}_j &= \Bigl[- \gamma_j w^{\tinyl}_j - \sum_{\ell=1}^{w^{\mbox{\tinyu}}_j - w^{\mbox{\tinyl}}_j} \varphi_{j\ell} \Bigr]_+ + \sup_{\tilde{d}_j \in [d^{\tinyl}_j, d^{\tinyu}_j]_{\mathbb{Z}}} \bigg\{ - \sum_{q=1}^Q \rho_{jq} \tilde{d}^q_j \bigg\}, \label{ref-linear-5} \\ 
c^{\text{\tiny t}}_{jk} &=\Bigl[(- c^{\mbox{\tiny x}}_k - \gamma_j) w^{\tinyl}_j - \sum_{\ell=1}^{w^{\mbox{\tinyu}}_j - w^{\mbox{\tinyl}}_j} (c^{\mbox{\tiny x}}_k u_{j \ell} + \varphi_{j \ell}) \Bigr]_+ + \sup_{\tilde{d}_j \in [d^{\tinyl}_j, d^{\tinyu}_j]_{\mathbb{Z}}} \bigg\{ c^{\text{\tiny x}}_k \tilde{d}_j - \sum_{q=1}^Q \rho_{jq} \tilde{d}^q_j \bigg\},\label{ref-linear-6} \\
c^{\text{\tiny p}}_i &= \Bigl[ - \lambda_i y^{\tinyl}_i - \sum_{\ell=1}^{y^{\mbox{\tinyu}}_i - y^{\mbox{\tinyl}}_i} \nu_{i\ell} \Bigr]_+, \label{ref-linear-7}\\
c^{\text{\tiny s}}_{ik} &= \Bigl[(- c^{\mbox{\tiny x}}_k - \lambda_i) y^{\tinyl}_i - \sum_{\ell=1}^{y^{\mbox{\tinyu}}_i - y^{\mbox{\tinyl}}_i} (c^{\mbox{\tiny x}}_k v_{i\ell} + \nu_{i\ell}) \Bigr]_+. \label{ref-linear-8}
\end{align}
\end{subequations}
\end{theorem}
}
\begin{algorithm}[ht!]
\caption{A separation algorithm for solving the (DRNS) model \eqref{drns}}
\label{algo-sep}
{\color{black}
\begin{algorithmic}[1]
\STATE Initialization: Set the set of cuts $\mathcal{H}_{\mbox{\tiny sep}} \leftarrow \emptyset$.
\STATE Solve the master problem
\begin{align*}
(\mbox{\bf MP}): \ \ \min_{\substack{u, v, \varphi, \nu\\ \gamma, \lambda, \rho}} \ & \ \theta + \sum_{i=1}^I \Bigl(c^{\mbox{\tiny y}}_i y^{\tinyl}_i + g_i(y^{\tinyl}_i) \lambda_i + \sum_{\ell=1}^{y^{\mbox{\tinyu}}_i - y^{\mbox{\tinyl}}_i} \bigl(c^{\mbox{\tiny y}}_i v_{i\ell} + \delta_{i\ell} \nu_{i\ell}  \bigl) \Bigr) \nonumber \\
& \ + \sum_{j=1}^J \Biggl[ \sum_{q=1}^Q \mu_{jq} \rho_{jq} + c^{\mbox{\tiny w}}_j w^{\tinyl}_j + f_j(w^{\tinyl}_j) \gamma_j + \sum_{\ell=1}^{w^{\mbox{\tinyu}}_j - w^{\mbox{\tinyl}}_j} \bigl( c^{\mbox{\tiny w}}_j u_{j\ell} + \Delta_{j\ell} \varphi_{j\ell} \bigr) \Biggr] \\
\mbox{s.t.} \ & \ \mbox{\eqref{ref-linear-2}, \eqref{ref-linear-4}}, \\
& \ \theta \geq \sum_{j=1}^J \bigg( c ^{\text{\tiny r}}_j r_j + \sum_{k =1}^j c^{\text{\tiny t}}_{jk} t_{jk} \bigg) + \sum_{i=1}^I \bigg( c^{\text{\tiny p}}_i p_i + \sum_{k = 1}^{J(i)} c^{\text{\tiny s}}_{ik} s_{ik} \bigg), \ \ \forall (t, s, r, p) \in \mathcal{H}_{\mbox{\tiny sep}},
\end{align*}
and record an optimal solution $(u^*, v^*, \varphi^*, \nu^*, \gamma^*, \lambda^*, \rho^*, \theta^*)$.
\STATE Compute $c^{\mbox{\tiny r}*}_j$, $c^{\mbox{\tiny t}*}_{jk}$, $c^{\mbox{\tiny p}*}_i$, and $c^{\mbox{\tiny s}*}_{ik}$ based on \eqref{ref-linear-5}--\eqref{ref-linear-8} and the values of $(u^*, v^*, \varphi^*, \nu^*, \gamma^*, \lambda^*, \rho^*)$.
\STATE Solve the integer linear program \eqref{IntegerProgram} using objective coefficients $c^{\mbox{\tiny r}*}_j$, $c^{\mbox{\tiny t}*}_{jk}$, $c^{\mbox{\tiny p}*}_i$, and $c^{\mbox{\tiny s}*}_{ik}$. Record an optimal solution $(t^*, s^*, r^*, p^*)$.
\IF{$\theta^* \geq \sum_{j=1}^J \big( c^{\text{\tiny r}*}_j r^*_j + \sum_{k =1}^j c^{\text{\tiny t}*}_{jk} t^*_{jk} \big) + \sum_{i=1}^I \big( c^{\text{\tiny p}*}_i p^*_i + \sum_{k = 1}^{J(i)} c^{\text{\tiny s}*}_{ik} s^*_{ik} \big)$}
\STATE Stop and return $(u^*, v^*)$ as an optimal solution to (DRNS).
\ELSE
\STATE Add a cut in the form of \eqref{ref-linear-3} into ({\bf MP}) by setting $\mathcal{H}_{\mbox{\tiny sep}} \leftarrow \mathcal{H}_{\mbox{\tiny sep}} \cup \{(t^*, s^*, r^*, p^*)\}$. Go To Step 2.
\ENDIF
\end{algorithmic}
}
\end{algorithm}

The reformulation \eqref{ref-linear-1}--\eqref{ref-linear-4} facilitates the separation algorithm (see, e.g.,~\cite{nemhauser1999integer}), also known as delayed constraint generation. We notice that \eqref{ref-linear-3} involve {\color{black} $(J+1)!$ } many constraints, making it computationally prohibitive to solve \eqref{ref-linear-1}--\eqref{ref-linear-4} in one shot. Instead, the separation algorithm incorporates constraints \eqref{ref-linear-3} on-the-fly. Specifically, this algorithm first solves a relaxation of the reformulation by overlooking constraints \eqref{ref-linear-3}. Then, we check if the optimal solution thus obtained violates any of \eqref{ref-linear-3}. If yes, then we add one violated constraint back into the relaxation and {\color{black}resolve}. We call this added constraint a ``cut'' and note that each cut describes a convex feasible region. This procedure {\color{black} is repeated} until an optimal solution is found to satisfy all of constraints \eqref{ref-linear-3}. We present the pseudo code in Algorithm \ref{algo-sep}.

We close this section by confirming the correctness of Algorithm \ref{algo-sep}.
\begin{theorem} \label{thm:finite}
Algorithm \ref{algo-sep} finds a globally optimal solution to the (DRNS) model \eqref{drns} in a finite number of iterations.
\end{theorem}
\proof See Appendix \ref{apx-thm:finite}. \qed
{\color{black}
To strengthen the integer program \eqref{IntegerProgram}, we derive valid inequalities for its feasible region.
\begin{lemma} \label{lem:valid}
	The following inequalities are valid for all $(t, s, r, p) \in \mathcal{H}$:
	\begin{align}
	\sum_{\ell = k }^j t_{j \ell} \leq \sum_{\ell = k}^{J(i)} s_{i \ell}, \ \forall i \in [I], \forall j \in P_i, \forall k \in [j],
	\label{VI-general}
	\end{align}
	where $J(i)$ is the largest index in $P_i$. In particular, inequalities \eqref{VI-general} imply constraints \eqref{encode-4}.
\end{lemma}
\proof See Appendix \ref{apx-lem:valid}. \qed 
In the following section, we show that, under special but practical pool structures, inequalities \eqref{VI-general} are strong enough to describe the convex hull of $\mathcal{H}$. In Section \ref{sec:results_performance}, we demonstrate that \eqref{VI-general} significantly accelerate the solution of (DRNS) under both general and special pool structures.}

\begin{figure}[H]
	\centering
	\begin{subfigure}[b]{0.25\textwidth}
		\includegraphics[width=\textwidth]{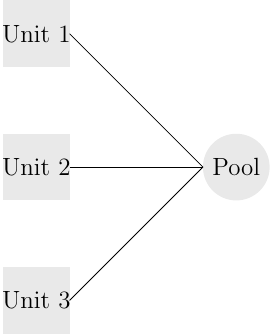}
		\caption{Structure \ref{assump:one-pool}}
		\label{Fig:one-pool}
	\end{subfigure} \hspace{0.5cm}
	\begin{subfigure}[b]{0.25\textwidth}
		\includegraphics[width=\textwidth]{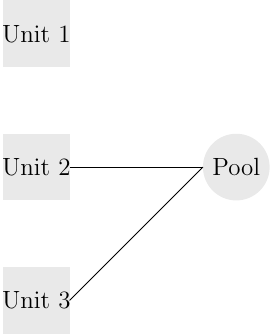}
		\caption{Structure \ref{assump:disjoint-pools}}
		\label{Fig:disjoint-pool}
	\end{subfigure} \hspace{0.5cm}
	\begin{subfigure}[b]{0.25\textwidth}
		\includegraphics[width=\textwidth]{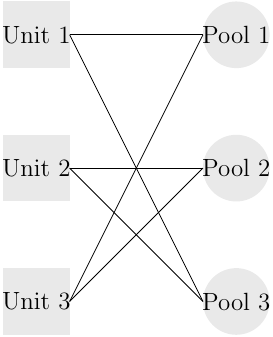}
		\caption{Structure \ref{assump:chained-pools}}
		\label{Fig:chained-pool}
	\end{subfigure}
	\caption{Examples of practical pool structures.} \label{Fig:chained-pools}
\end{figure}

\section{Tractable Cases: Practical Pool Structures} \label{sec:tractable}
In this section, we consider the following three special nurse pool structures.
\begin{customstructure}{1} \label{assump:one-pool}
\emph{(One Pool)} $I = 1$, i.e., there is one single nurse pool shared among all units \emph{(see Figure~\ref{Fig:one-pool} for an example)}.
\end{customstructure}

\begin{customstructure}{D} \label{assump:disjoint-pools}
\emph{(Disjoint Pools)} All nurse pools are disjoint, i.e., for all $i_1, i_2 \in [I]$ and $i_1 \neq i_2$, it holds that $P_{i_1} \cap P_{i_2} = \emptyset$ \emph{(see Figure~\ref{Fig:disjoint-pool} for an example)}.
\end{customstructure}

\begin{customstructure}{C} \label{assump:chained-pools}
\emph{(Chained Pools)} The nurse pools form a long chain, i.e., there are $I = J$ pools and the bipartite graph consisting of the units, pools, and pool-unit assignment forms one undirected cycle. Figure~\ref{Fig:chained-pool} displays an example with cycle $P_1$--Unit 1--$P_3$--Unit 2--$P_2$--Unit 3--$P_1$.
\end{customstructure}

{\color{black} Structure \ref{assump:one-pool}} can be utilized when all units have similar functionalities and so they can all share one nurse pool. Accordingly, every nurse assigned to this pool should be cross-trained for all units so that he/she is able to undertake the tasks in them. 
{\color{black} Structure \ref{assump:disjoint-pools}} is less demanding than one pool, as each pool covers only a subset of units which, e.g., have distinct functionalities. Accordingly, the amount of cross-training under this structure significantly decreases from that under one pool. 
{\color{black} Structure \ref{assump:chained-pools}} has been applied in the production systems to increase the operational flexibility (see, e.g.,~\cite{jordan1995principles,wang2015process,chen2015optimal,chen2019optimal}). Under this structure, every unit is covered by two nurse pools. Accordingly, every pool nurse needs to be cross-trained for only two units. All three structures have been considered and compared in a nurse staffing context (see, e.g.,~\cite{inman2005cross}). Under these practical pool structures, we derive tractable reformulations of the (DRNS) model \eqref{drns}. Our derivation leads {\color{black} to MILP} reformulations that facilitate off-the-shelf software like GUROBI.

\subsection{One Pool} \label{subsec:one-pool}
{\color{black}
Under Structure 1, i.e., $I=1$, $P_1=[J]$, and $J(1)=J$, we rewrite the feasible region $\mathcal{H}$ of integer program \eqref{IntegerProgram} as follows:
\begin{subequations}
	\begin{align}
	\mathcal{H}_{\text{\tiny 1}} = \bigg\{ (t,s,r,p) : \
	& r_j + \sum_{k =1}^j t_{jk} = 1, \ \forall j \in [J], \ \ p + \sum_{k =1}^J s_{k} = 1, \label{H-onepool-0} \\	
	& \sum_{k = 1}^j t_{jk} \leq 1, \ \forall j \in [J], \label{H-onepool-1} \\ 
	& \sum_{k=1}^{J} s_{k} \leq 1, \label{H-onepool-2} \\ 	
	& s_{k} \leq \sum_{j =k}^J t_{jk}, \ \forall k \in [J], \label{H-onepool-3}\\
	& s_{k} + t_{j \ell} \leq 1, \ \forall k \in [J], \forall j \in [J],  \forall \ell \in [k+1, j]_{\mathbb{Z}}, \label{H-onepool-4}\\
	& t_{jk} \leq \sum_{\ell = 1}^{J} s_{\ell}, \ \forall j \in [J], \forall k \in [j], \label{H-onepool-5} \\
	& s_{k} \in \mathbb{B}, \ \forall k \in [J], \ \ t_{jk} \in \mathbb{B}, \ \forall j \in [J], \forall k \in [j]  \label{H-onepool-6}
	\bigg\}.
	\end{align}
\end{subequations}
We claim that valid inequalities \eqref{VI-general}, in conjunction with constraints \eqref{H-onepool-0}--\eqref{H-onepool-3},
are sufficient to describe the convex hull of $\mathcal{H}_{\text{\tiny 1}}$.

\begin{proposition} \label{prop:onepool-convex}
Consider a polyhedron $\overline{\mathcal{H}}_{\text{\tiny 1}}$ given by
\begin{subequations}
	\begin{align}
	\overline{\mathcal{H}}_{\text{\tiny 1}} := \bigg\{ (t,s,r,p) : \ 
	& \mbox{\eqref{H-onepool-0}--\eqref{H-onepool-3}}, \nonumber \\
	& \sum_{\ell = k }^j t_{j \ell} \leq \sum_{\ell = k}^{J} s_{\ell}, \ \forall j \in [J], \forall k \in [j], \label{onepool-VI} \\
	& s_{k} \in \mathbb{R}_+, \ \forall k \in [J], \ \ t_{jk} \in \mathbb{R}_+, \ \forall j \in [J], \forall k \in [j]  \label{onepool-nonneg}
	\bigg\}.
	\end{align}	
\end{subequations}
Then, $\overline{\mathcal{H}}_{\text{\tiny 1}} = \text{conv} (\mathcal{H}_{\text{\tiny 1}})$.
\end{proposition}
\proof See Appendix \ref{apx-prop:onepool-convex}.  \qed
Better still, this yields a closed-form solution to the convex maximization problem $\displaystyle \max_{(\alpha,\beta)\in \Lambda} F(\alpha,\beta)$.

\begin{theorem} \label{thm:integrality-one}
Under Structure \ref{assump:one-pool},
for fixed $(u, v, \gamma, \lambda, \rho)$, the optimal value of $\displaystyle \max_{(\alpha, \beta) \in \Lambda} F(\alpha, \beta)$ equals
\begin{align*}
\max \Bigg\{
c^{\text{\tiny p}}_1 + \sum_{j=1}^J c^{\text{\tiny r}}_j, \ \ \max_{ \substack{j^* \in [J], \\ k^* \in [j^*]} } \bigg\{ \sum_{j=1}^J c^{\text{\tiny r}}_j + c^{\text{\tiny s}}_{1k^*}  +  (c^{\text{\tiny t}}_{j^* k^*} - c^{\text{\tiny r}}_{j^*}) + \sum_{j \in [J] \setminus \{j^*\} } \Big( \max_{k \in [\min \{j, k^* \}] } \ \big\{ c^{\text{\tiny t}}_{jk} - c^{\text{\tiny r}}_j \big\} \Big) \bigg\}
\Bigg\},
\end{align*}
where $c^{\mbox{\tiny r}}_j$, $c^{\mbox{\tiny t}}_{jk}$, $c^{\mbox{\tiny p}}_1$, and $c^{\mbox{\tiny s}}_{1k}$ are defined through \eqref{ref-linear-5}--\eqref{ref-linear-8}.
\end{theorem}
\proof See Appendix \ref{apx-thm:integrality-one}. \qed
}
{\color{black}Theorem \ref{thm:integrality-one}} enables us to reduce the {\color{black}$(J+1)!$} many constraints \eqref{ref-linear-3} in the reformulation of (DRNS) to {\color{black}$\frac{1}{2}J(J+1) + 1$} many, thanks to the closed-form solution of $\max_{(\alpha, \beta) \in \Lambda} F(\alpha, \beta)$. This leads to the {\color{black}following MILP} reformulation of (DRNS).
\begin{proposition} \label{prop:one-pool-ref}
Under Assumption \ref{assump:technical} and {\color{black}Structure \ref{assump:one-pool}}, the (DRNS) model \eqref{drns} yields the same optimal objective value as the following MILP:
\begin{subequations}
{\color{black}
\begin{align}
Z^{\star}_{1} := \min \ & \ \theta + c^{\mbox{\tiny y}}_1 y^{\tinyl}_1 + g_1(y^{\tinyl}_1) \lambda_1 + \sum_{\ell=1}^{y^{\mbox{\tinyu}}_1 - y^{\mbox{\tinyl}}_1} \bigl(c^{\mbox{\tiny y}}_1 v_{1\ell}  + \delta_{1\ell} \nu_{1\ell} \bigl) \nonumber \\
& \ + \sum_{j=1}^J \Biggl[ \sum_{q=1}^Q \mu_{jq} \rho_{jq} + c^{\mbox{\tiny w}}_j w^{\tinyl}_j + f_j(w^{\tinyl}_j) \gamma_j + \sum_{\ell=1}^{w^{\mbox{\tinyu}}_j - w^{\mbox{\tinyl}}_j} \bigl( c^{\mbox{\tiny w}}_j u_{j\ell} + \Delta_{j\ell} \varphi_{j\ell} \bigr) \Biggr] \nonumber \\
\mbox{s.t.} \ & \ \mbox{\eqref{ref-note-36}--\eqref{ref-note-33}}, \ \mbox{\eqref{ref-linear-4}}, \nonumber \\
& \theta \geq \phi^{\mbox{\tiny e}}_1 + \sum_{j=1}^J (\zeta^{\mbox{\tiny e}}_j + \eta^{\mbox{\tiny e}}_j), \nonumber \\
& \theta \geq \sum_{j=1}^J (\zeta^{\mbox{\tiny e}}_j + \eta^{\mbox{\tiny e}}_j) + \phi^{\mbox{\tiny x}}_{1k^*} + (\zeta^{\mbox{\tiny x}}_{j^* k^*} + \eta^{\mbox{\tiny x}}_{j^* k^*}) - (\zeta^{\mbox{\tiny e}}_{j^*} + \eta^{\mbox{\tiny e}}_{j^*}) + \sum_{j \in  [J] \setminus \{j^*\}} \chi_{jk^*}, \ \forall j^* \in [J], \forall k^* \in [j^*], \nonumber \\	
& \chi_{jk^*} \geq  (\zeta^{\mbox{\tiny x}}_{jk} + \eta^{\mbox{\tiny x}}_{jk}) - (\zeta^{\mbox{\tiny e}}_j + \eta^{\mbox{\tiny e}}_j), \ \forall j \in [J], \forall k^* \in [J], \forall k \in [\min\{j,k^*\}], \label{onepool-ref-1} \\
& \hspace{-0.15cm} \left.\begin{array}{l}
\zeta^{\mbox{\tiny e}}_j \geq  - \gamma_j w^{\tinyl}_j - \sum_{\ell=1}^{w^{\mbox{\tinyu}}_j - w^{\mbox{\tinyl}}_j}  \varphi_{j\ell} , \ \ \ \zeta^{\mbox{\tiny e}}_j \geq 0, \\[0.1cm]
\eta^{\mbox{\tiny e}}_j \geq  - \sum_{q=1}^Q\tilde{d}^q_j\rho_{jq}, \ \forall \tilde{d}_j \in [d^{\tinyl}_j, d^{\tinyu}_j]_{\mathbb{Z}},
\end{array}\right\} \ \ \forall j \in [J], \label{onepool-ref-3} \\
& \hspace{-0.15cm} \left.\begin{array}{l}
\zeta^{\mbox{\tiny x}}_{jk} \geq (- c^{\mbox{\tiny x}}_k - \gamma_j) w^{\tinyl}_j - \sum_{\ell=1}^{w^{\mbox{\tinyu}}_j - w^{\mbox{\tinyl}}_j} ( c^{\mbox{\tiny x}}_k u_{j \ell} + \varphi_{j \ell}), \ \ \ \zeta^{\mbox{\tiny x}}_{jk} \geq 0, \\[0.1cm]
\eta^{\mbox{\tiny x}}_{jk} \geq c^{\mbox{\tiny x}}_k \tilde{d}_j - \sum_{q=1}^Q\tilde{d}^q_j\rho_{jq}, \ \forall \tilde{d}_j \in [d^{\tinyl}_j, d^{\tinyu}_j]_{\mathbb{Z}}, \\[0.1cm]
\end{array}\right\} \ \ \forall j \in [J], \forall k \in [j], \label{onepool-ref-4} \\
& \ \phi^{\mbox{\tiny e}}_1 \geq  - \lambda_1 y^{\tinyl}_1 - \sum_{\ell=1}^{y^{\mbox{\tinyu}}_1 - y^{\mbox{\tinyl}}_1} \nu_{1\ell}, \ \ \ \phi^{\mbox{\tiny e}}_1 \geq 0, \nonumber \\
& \ \phi^{\mbox{\tiny x}}_{1k} \geq (- c^{\mbox{\tiny x}}_k - \lambda_1) y^{\tinyl}_1 - \sum_{\ell=1}^{y^{\mbox{\tinyu}}_1 - y^{\mbox{\tinyl}}_1} (c^{\mbox{\tiny x}}_k v_{1\ell} + \nu_{1\ell} ), \ \ \phi^{\mbox{\tiny x}}_{1k} \geq 0, \ \forall k \in [J]. \nonumber
\end{align}}
\end{subequations}
\end{proposition}
{\color{black} \proof See Appendix \ref{apx-prop:one-pool-ref}. \qed }

A special case of {\color{black}Structure \ref{assump:one-pool}} is when there are no nurse float pools. Mathematically, this is equivalent to assigning all units to a pool with no nurses. We hence call it {\color{black}Structure 0} as there is zero pool nurse. Under this structure, $y^{\tinyl}_1 = y^{\tinyu}_1 = 0$ and accordingly $g_1(y^{\tinyl}_1) = 0$. A MILP reformulation of (DRNS) under {\color{black}Structure 0} follows from Proposition \ref{prop:one-pool-ref}:
{\color{black}
\begin{align*}
Z^{\star}_{0} \ := \ \min \ & \ \theta + \sum_{j=1}^J \Biggl[ \sum_{q=1}^Q \mu_{jq} \rho_{jq} + c^{\mbox{\tiny w}}_j w^{\tinyl}_j + f_j(w^{\tinyl}_j) \gamma_j + \sum_{\ell=1}^{w^{\mbox{\tinyu}}_j - w^{\mbox{\tinyl}}_j} \bigl( c^{\mbox{\tiny w}}_j u_{j\ell} + \Delta_{j\ell} \varphi_{j\ell} \bigr) \Biggr] \nonumber \\
\mbox{s.t.} \ & \ \mbox{\eqref{ref-note-36}}, \ \mbox{\eqref{ref-note-37}}, \ \mbox{\eqref{ref-note-32}}, \ \mbox{\eqref{ref-linear-4}}, \ \mbox{\eqref{onepool-ref-1}--\eqref{onepool-ref-4}}, \nonumber \\
& \theta \geq \sum_{j=1}^J (\zeta^{\mbox{\tiny e}}_j + \eta^{\mbox{\tiny e}}_j), \nonumber \\
& \theta \geq \sum_{j=1}^J (\zeta^{\mbox{\tiny e}}_j + \eta^{\mbox{\tiny e}}_j) + (\zeta^{\mbox{\tiny x}}_{j^* k^*} + \eta^{\mbox{\tiny x}}_{j^* k^*}) - (\zeta^{\mbox{\tiny e}}_{j^*} + \eta^{\mbox{\tiny e}}_{j^*}) + \sum_{j \in  [J] \setminus \{j^*\}} \chi_{jk^*}, \ \forall j^* \in [J], \forall k^* \in [j^*]. \nonumber
\end{align*}}
\noindent
We notice that, whenever $y^{\tinyl}_1 = 0$, any feasible nurse staffing levels under {\color{black}Structure 0} are also feasible to (DRNS) under {\color{black}Structure \ref{assump:one-pool}}. It then follows that {\color{black}$Z^{\star}_{1} \leq Z^{\star}_{0}$}. In addition, since {\color{black}Structure \ref{assump:one-pool}} provides the most operational flexibility and {\color{black}Structure 0} has zero flexibility, we  interpret the difference {\color{black}$Z^{\star}_{0} - Z^{\star}_{1}$} as the (maximum) value of operational flexibility.

\subsection{Disjoint Pools} \label{subsec:disjoint-pools}
{\color{black}Under Structure \ref{assump:disjoint-pools}, we denote by $\mathcal{H}_{\text{\tiny D}}$ the feasible region $\mathcal{H}$ of the integer program \eqref{IntegerProgram}.
We can once again obtain the convex hull of $\mathcal{H}_{\text{\tiny D}}$
by incorporating inequalities \eqref{VI-general}. }
Intuitively, as the nurse pools are disjoint, {\color{black}$\mathcal{H}_{\text{\tiny D}}$} becomes separable in index $i$, i.e., separable among the nurse pools and the units under each pool. Hence, {\color{black}$\mbox{conv}(\mathcal{H}_{\text{\tiny D}})$} can be obtained by convexifying the projection of {\color{black}$\mathcal{H}_{\text{\tiny D}}$} in each pool and then taking their Cartesian product. 
It follows that, once again, the convex maximization problem $\max_{(\alpha, \beta)}F(\alpha, \beta)$ admits a closed-form solution and (DRNS) can be recast as {\color{black}a MILP}. In particular, we reduce the exponentially many constraints \eqref{ref-linear-3} in the reformulation of (DRNS) to {\color{black}$\frac{1}{2}J(J+1) + I$} many. We summarize these results in the following proposition.
{\color{black}
\begin{proposition} \label{prop:exact-disjoint}
Under Structure \ref{assump:disjoint-pools}, it holds that $\mbox{conv}(\mathcal{H}_{\text{\tiny D}}) \ = \ \{(t,s,r,p): \mbox{\eqref{IntegerProgram-1}--\eqref{IntegerProgram-3}, \eqref{VI-general}}\}$. In addition, for fixed $(u, v, \gamma, \lambda, \rho)$, the optimal value of $\max_{(\alpha,\beta) \in \Lambda} F(\alpha, \beta)$ equals
\begin{align*}
& \sum_{i=1}^I \Bigg[
\max \Bigg\{
c^{\text{\tiny p}}_i + \sum_{j \in P_i} c^{\text{\tiny r}}_j, \ \ \max_{ \substack{j^* \in P_i,\\ k^* \in [j^*]} } \bigg\{ \sum_{j \in P_i} c^{\text{\tiny r}}_j + c^{\text{\tiny s}}_{ik^*}  +  (c^{\text{\tiny t}}_{j^* k^*} - c^{\text{\tiny r}}_{j^*}) + \sum_{j \in P_i \setminus \{j^*\} } \bigg( \max_{k \in [\min \{j, k^* \} ] } \ \Big\{ c^{\text{\tiny t}}_{jk} - c^{\text{\tiny r}}_j \Big\} \bigg) \bigg\}
\Bigg\}
\Bigg],
\end{align*}
where $c^{\mbox{\tiny r}}_j$, $c^{\mbox{\tiny t}}_{jk}$, $c^{\mbox{\tiny p}}_i$, and $c^{\mbox{\tiny s}}_{ik}$ are defined through \eqref{ref-linear-5}--\eqref{ref-linear-8}. Furthermore, under Assumption \ref{assump:technical}, the (DRNS) model \eqref{drns} yields the same optimal objective value as the following MILP:
\begin{subequations}
\begin{align}
\hspace{-0.0cm} Z^{\star}_{\text{\tiny D}} \ := \ \min \ & \ \sum_{i=1}^I \Bigl( \theta_i + c^{\mbox{\tiny y}}_i y^{\tinyl}_i + g_i(y^{\tinyl}_i) \lambda_i + \sum_{\ell=1}^{y^{\mbox{\tinyu}}_i - y^{\mbox{\tinyl}}_i} \bigl( c^{\mbox{\tiny y}}_i v_{i\ell} + \delta_{i\ell} \nu_{i\ell} \bigr) \Bigr) \nonumber \\
& \ + \sum_{j=1}^J \Biggl[ \sum_{q=1}^Q \mu_{jq} \rho_{jq} + c^{\mbox{\tiny w}}_j w^{\tinyl}_j + f_j(w^{\tinyl}_j) \gamma_j + \sum_{\ell=1}^{w^{\mbox{\tinyu}}_j - w^{\mbox{\tinyl}}_j} \bigl( c^{\mbox{\tiny w}}_j u_{j\ell} + \Delta_{j\ell} \varphi_{j\ell} \bigr) \Biggr] \nonumber \\
\mbox{s.t.} \ & \ \mbox{\eqref{ref-note-36}--\eqref{ref-note-33}}, \ \mbox{\eqref{ref-linear-4}}, \ \mbox{\eqref{onepool-ref-1}--\eqref{onepool-ref-4}}, \nonumber \\ 
& \theta_i \geq \phi^{\mbox{\tiny e}}_i + \sum_{j \in P_i} (\zeta^{\mbox{\tiny e}}_j + \eta^{\mbox{\tiny e}}_j), \ \forall i \in [I], \nonumber \\
& \theta_i \geq \sum_{j \in P_i} (\zeta^{\mbox{\tiny e}}_j + \eta^{\mbox{\tiny e}}_j) + \phi^{\mbox{\tiny x}}_{ik^*} + (\zeta^{\mbox{\tiny x}}_{j^* k^*} + \eta^{\mbox{\tiny x}}_{j^* k^*}) - (\zeta^{\mbox{\tiny e}}_{j^*} + \eta^{\mbox{\tiny e}}_{j^*}) + \sum_{j \in P_i \setminus \{j^*\} } \chi_{jk^*}, \nonumber \\
& \hspace{9mm} \forall i \in [I], \forall j^* \in P_i, \forall k^* \in [j^*], \nonumber \\
& \phi^{\mbox{\tiny e}}_i \geq   - \lambda_i y^{\tinyl}_i - \sum_{\ell=1}^{y^{\mbox{\tinyu}}_i - y^{\mbox{\tinyl}}_i} \nu_{i\ell}, \ \ \ \phi^{\mbox{\tiny e}}_i \geq 0, \ \forall i \in [I], \label{disjoint-ref-1} \\
& \phi^{\mbox{\tiny x}}_{ik} \geq (- c^{\mbox{\tiny x}}_k - \lambda_i) y^{\tinyl}_i - \sum_{\ell=1}^{y^{\mbox{\tinyu}}_i - y^{\mbox{\tinyl}}_i} (c^{\mbox{\tiny x}}_k v_{i\ell} + \nu_{i\ell} ), \ \ \phi^{\mbox{\tiny x}}_{ik} \geq 0, \ \forall i \in [I], \forall k \in [J(i)]. \label{disjoint-ref-2}
\end{align}
\end{subequations}
\end{proposition}
\proof See Appendix \ref{apx-prop:exact-disjoint}. \qed
}

\subsection{Chained Pools} \label{subsec:chained-pools}
Under Structure \ref{assump:chained-pools}, the valid inequalities \eqref{VI-general} can still be incorporated to strengthen and simplify the mixed-integer set $\mathcal{H}$. Unfortunately, unlike under Structures \ref{assump:one-pool} and \ref{assump:disjoint-pools}, the strengthened $\mathcal{H}$ is no longer integral. We demonstrate this fact in Appendix~\ref{apx:example-chained-pool}. Despite the loss of integrality, we adopt an alternative approach to recast the integer program \eqref{IntegerProgram}, and hence the convex maximization problem $\displaystyle\max_{(\alpha, \beta) \in \Lambda} F(\alpha, \beta)$, as a linear program.
{\color{black}
We start by noticing that, by re-numbering the units if needed, the chained pools can be made $P_i = \{i, i+1\}$ for all $i \in [I-1]$ and $P_I = \{I, 1\}$. Note that after re-numbering, the costs of hiring temporary nurses $c^{\mbox{\tiny x}}_i$ may no longer be monotone. For ease of exposition, we denote by $\sigma(i)$ the original index of unit $i$ before re-numbering. We define an order $\prec$ on $[I]$ such that $i \prec j$ if and only if $\sigma(i) < \sigma(j)$, and we say $i \preceq j$ if $i \prec j$ or $i = j$. In addition, we let $i \curlyvee j$ represent the index of higher order between units $i$ and $j$, i.e., $i \curlyvee j = j$ if $i \preceq j$ and $i \curlyvee j = i$ otherwise. In addition, we create a dummy unit $I+1$, which is a duplicate of unit $1$, and refer to these two units interchangeably.
\begin{figure}[H]
	\centering
	\begin{subfigure}[b]{0.25\textwidth}
		\includegraphics[width=\textwidth]{./Figure/Chained_Pool_Before.pdf}
		\caption{Before re-numbering}
		\label{Fig:chained-pool-before}
	\end{subfigure} \hspace{3.0cm}
	\begin{subfigure}[b]{0.25\textwidth}
		\includegraphics[width=\textwidth]{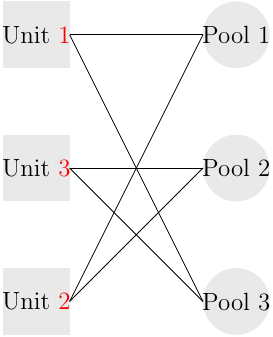}
		\caption{After re-numbering}
		\label{Fig:chained-pool-after}
	\end{subfigure}
	\caption{Re-numbering units in Structure \ref{assump:chained-pools}} \label{Fig:chained-pools-re-numbering}
\end{figure}
\begin{example} \label{ex:chained-pool}
Consider the chained pools in Figure~\ref{Fig:chained-pool-before} with monotone costs $c^{\mbox{\tiny x}}_1 \leq c^{\mbox{\tiny x}}_2 \leq c^{\mbox{\tiny x}}_3$. After re-numbering (see Figure~\ref{Fig:chained-pool-after}), we have $\sigma(1) = 1$, $\sigma(2) = 3$, and $\sigma(3) = 2$. In addition, $1 \prec 3 \prec 2$ and $2 \curlyvee 3 = 2$.
\end{example}
Next, since $P_i = \{i, i+1\}$ for all $i \in [I]$, variables $s_{ik}$ in the integer program \eqref{IntegerProgram} depend on $t_{[i]}$ and $t_{[i+1]}$ only, where $t_{[i]} := \{t_{ik}: k \preceq i\}$. Indeed, the definition of $s_{ik}$ implies
\begin{align*}
s_{ik} = h_k(t_{[i]}, t_{[i+1]}) := & \left\{\begin{array}{ll}
1 & \mbox{if $\max\{t_{ik}, t_{(i+1)k}\} = 1$, $t_{i\ell} = 0$ for all $\ell \succ k$,} \\
  & \mbox{and $t_{(i+1)\ell} = 0$ for all $\ell \succ k$}, \\[0.2cm]
0 & \mbox{otherwise}.
\end{array}\right.
\end{align*}}
Plugging these representations into formulation \eqref{IntegerProgram} yields a reformulation of $\displaystyle\max_{(\alpha, \beta) \in \Lambda} F(\alpha, \beta)$ based on variables $t$ only:
{\color{black} 
\begin{align}
\max_{(\alpha, \beta) \in \Lambda} F(\alpha, \beta) \ = \ \sum_{i=1}^I(c^{\mbox{\tiny r}}_i + c^{\mbox{\tiny p}}_i) + \sum_{i=1}^I\left[ \sum_{k \preceq i} (c^{\mbox{\tiny t}}_{ik} - c^{\mbox{\tiny r}}_i) t_{ik} + \sum_{k \preceq \ i \curlyvee (i+1)} (c^{\mbox{\tiny s}}_{ik} - c^{\mbox{\tiny p}}_i) h_k(t_{[i]}, t_{[i+1]}) \right]. \label{chained-pool-note-1}
\end{align}
}
The reformulation \eqref{chained-pool-note-1} decomposes objective function based on index $i \in [I]$ and enables us to solve $\displaystyle\max_{(\alpha, \beta) \in \Lambda} F(\alpha, \beta)$ by a dynamic program (DP), i.e., we sequentially optimize $t_{[1]}, t_{[2]}, \ldots, t_{[I]}$. To this end, we define the state of the DP in stage 1 as $t_{[1]} \equiv t_{11} \in \mathcal{B}_1$ and the states in stage $i$ as $(t_{[1]}, t_{[i]}) \in \mathcal{B}_1 \times \mathcal{B}_i$ for all $i \in [2, I]_{\mathbb{Z}}$, {\color{black}where $\mathcal{B}_i := \{t_{[i]} \in \mathbb{B}^{\sigma(i)}: \sum_{k \preceq i} t_{ik} \leq 1\}$ represents the feasible region of $t_{[i]}$ described by constraints \eqref{encode-0}.} In addition, we formulate the DP as {\color{black}$\displaystyle\max_{(t_{[1]}, t_{[I]}) \in \mathcal{B}_1 \times \mathcal{B}_I}\{V_I(t_{[1]}, t_{[I]}) + \sum_{k \preceq (1 \curlyvee I)} (c^{\mbox{\tiny s}}_{Ik} - c^{\mbox{\tiny p}}_I) h_k(t_{[I]}, t_{[1]})\}$} where the value functions $V_i(\cdot)$ are recursively defined through
{\color{black}
\begin{align*}
V_1(t_{[1]}) & := \sum_{i=1}^I(c^{\mbox{\tiny r}}_i + c^{\mbox{\tiny p}}_i) + (c^{\mbox{\tiny t}}_{11} - c^{\mbox{\tiny r}}_1) t_{11} \ \ \mbox{and} \\
V_i(t_{[1]}, t_{[i]}) & := \max_{t_{[i-1]} \in \mathcal{B}_{i-1}}\bigg\{V_{i-1}(t_{[1]}, t_{[i-1]}) + \sum_{k \preceq i}(c^{\mbox{\tiny t}}_{ik} - c^{\mbox{\tiny r}}_i) t_{ik} + \sum_{k \preceq (i-1) \curlyvee i}(c^{\mbox{\tiny s}}_{(i-1)k} - c^{\mbox{\tiny p}}_{i-1}) h_k(t_{[i-1]}, t_{[i]}) \bigg\}
\end{align*}
for all $i \in [2, I]_{\mathbb{Z}}$ and $(t_{[1]}, t_{[i]}) \in \mathcal{B}_1 \times \mathcal{B}_i$.} For all $i \in [I]$, value function $V_i(t_{[1]}, t_{[i]})$ represents the ``cumulative reward'' up to stage $i$, i.e., the terms in \eqref{chained-pool-note-1} that involve $t_{[1]}, \ldots, t_{[i]}$ only. We note that, as $t_{[1]}$ is involved in the final-stage reward, the DP stores the value of $t_{[1]}$ in the state throughout stages $2, \ldots, I$.

\begin{figure}[H]
	\centering 
	\includegraphics*[scale=0.7]{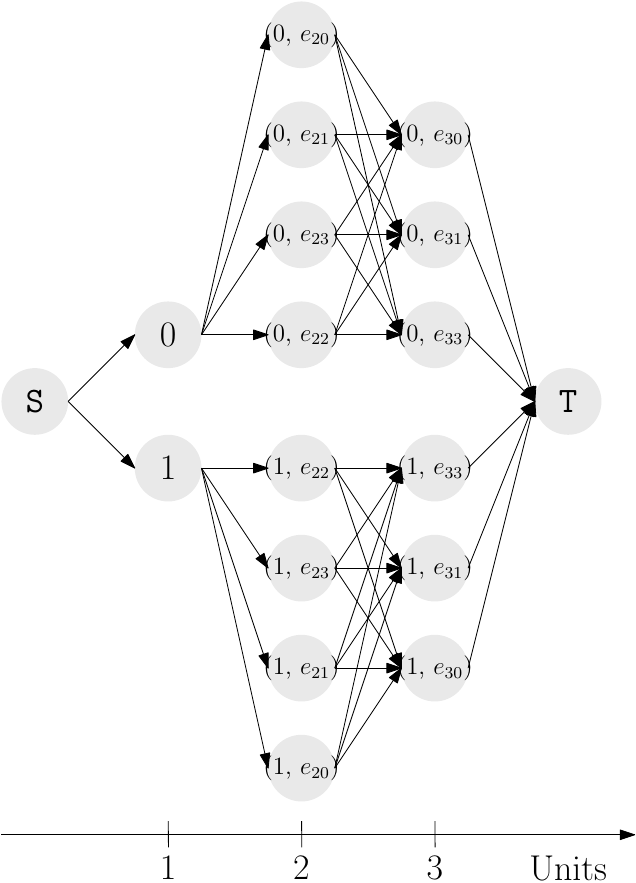}
	\caption{Longest-path problem on an acyclic directed network.}
	\label{fig:chained_pools_example}
\end{figure}

We further interpret the DP as a longest-path problem on an acyclic directed network $(\mathcal{N},\mathcal{A})$.
Specifically, the set of nodes $\mathcal{N}$ consists of $I$ layers, denoted by $\{N_i\}_{i=1}^I$.
For all $i \in [I]$, layer $i$ consists of the states of the DP in stage $i$, i.e., {\color{black}$\mathcal{N}_1 = \mathcal{B}_1 = \{0,1\}$} and {\color{black}$\mathcal{N}_i= \mathcal{B}_1 \times \mathcal{B}_i$ for all $i \in [2,I]_{\mathbb{Z}}$.}
In addition, $\mathcal{A}$ consists of arcs that connect two nodes in neighboring layers, as long as the two nodes share a common {\color{black}$t_{[1]}$ value, i.e., $\mathcal{A}= \{ [t_{[1]}, (t_{[1]}, t_{[2]}): t_{[1]} \in \mathcal{B}_1, t_{[2]} \in \mathcal{B}_2\} \cup \{ [(t_{[1]}, t_{[i-1]}), (t_{[1]}, t_{[i]})]: t_{[i-1]} \in \mathcal{B}_{i-1}, t_{[i]} \in \mathcal{B}_{i}, \forall i \in [3,I]_{\mathbb{Z}}\}$.} Finally, we incorporate into $\mathcal{N}$ a starting node $\texttt{S}$ and a terminal node $\texttt{T}$, and into $\mathcal{A}$ arcs from $\texttt{S}$ to all nodes in $\mathcal{N}_1$ and from all nodes in $\mathcal{N}_I$ to $\texttt{T}$. Then, the DP is equivalent to the longest-path problem from $\texttt{S}$ to $\texttt{T}$ on $(\mathcal{N},\mathcal{A})$. We formally state this result in the following theorem. {\color{black}For notation brevity, we say $t_{[i]} = e_{ik}$ if $t_{ik} = 1$ for a $k \preceq i$, while $t_{i\ell} = 0$ for all $\ell \preceq i$ and $\ell \neq k$; and we say $t_{[i]} = e_{i0}$ if $t_{ik} = 0$ for all $k \preceq i$.} In Figure~\ref{fig:chained_pools_example}, {\color{black}we depict $(\mathcal{N},\mathcal{A})$ for Example \ref{ex:chained-pool}.}

\begin{theorem} \label{thm:longest-path}
	Define $\{c_{[m,n]} : [m,n] \in \mathcal{A} \}$, the length of the arcs in network $(\mathcal{N},\mathcal{A})$, such that
{\color{black}
	\begin{align*}
	c_{[\texttt{S},t_{[1]}]} &= \sum_{i=1}^I(c^{\mbox{\tiny r}}_i + c^{\mbox{\tiny p}}_i) + (c^{\text{\tiny t}}_{11} - c^{\text{\tiny r}}_1) t_{11} \quad \forall t_{[1]} \in \mathcal{B}_1, \\
	c_{[t_{[1]}, (t_{[1]}, e_{2\ell})]} &=
	\begin{cases}
	( c^{\text{\tiny s}}_{11} - c^{\text{\tiny p}}_{1} ) t_{11} & \text{if } \ell = 0 \\
	(c^{\text{\tiny t}}_{2k} - c^{\text{\tiny r}}_{2} ) + (c^{\text{\tiny s}}_{1k} - c^{\text{\tiny p}}_{1}) & \text{if } \ell \neq 0 \\
	\end{cases}
	\quad \forall t_{[1]} \in \mathcal{B}_1, \forall \ell \preceq 2, \\
	c_{[(t_{[1]}, e_{(i-1)k}), (t_{[1]}, e_{i\ell})]} &=
	\begin{cases}
	0 & \text{if } k = 0, \ell = 0\\
	c^{\text{\tiny s}}_{(i-1)k} - c^{\text{\tiny p}}_{i-1} , & \text{if } k \neq 0, \ell = 0 \\
	(c^{\text{\tiny t}}_{i\ell} - c^{\text{\tiny r}}_{i}) + (c^{\text{\tiny s}}_{(i-1)(k \curlyvee \ell)} - c^{\text{\tiny p}}_{i-1}), & \text{if } \ell \neq 0 \\
	\end{cases} 
	\\
	& \ \ \ \ \ \
	\forall t_{[1]} \in \mathcal{B}_1, \forall i \in [3,I]_{\mathbb{Z}}, \forall k \preceq i-1, \forall \ell \preceq i, \\
	\text{ and } \ c_{ [(t_{[1]}, e_{Ik}), \texttt{T}] }  &=
	\begin{cases}
	(c^{\text{\tiny s}}_{I1} - c^{\text{\tiny p}}_{I}) t_{11} & \text{if } k = 0 \\
	c^{\text{\tiny s}}_{Ik} - c^{\text{\tiny p}}_{I} & \text{if } k \neq 0
	\end{cases}
	\quad \forall t_{[1]} \in \mathcal{B}_1, \forall k \preceq I.
	\end{align*}
}
\noindent Then, for fixed $(u, v, \gamma, \lambda, \rho)$, $\max_{(\alpha, \beta) \in \Lambda} F(\alpha, \beta)$ equals the length of the longest \texttt{S}-\texttt{T} path on $(\mathcal{N}, \mathcal{A})$, that is,
	\begin{align*}
	\max_{(\alpha, \beta) \in \Lambda} F(\alpha, \beta) = \max_{x \in [0,1]^{\mathcal{A}}} \ & \sum_{[m,n] \in \mathcal{A}} c_{[m,n]} x_{[m,n]} \\
	\mbox{s.t.} \
	& \sum_{n: [m,n] \in \mathcal{A}} x_{[m,n]} - \sum_{n: [n,m] \in \mathcal{A}} x_{[n,m]} = 
	\begin{cases}
	1   & \text{if } m = \texttt{S} \\
	0   & \text{if } m \neq \texttt{S},\texttt{T} \\
	-1  & \text{if } m = \texttt{T} \\
	\end{cases}
	, \ \forall m \in \mathcal{N}.
	\end{align*}
\end{theorem} 
\proof See Appendix \ref{apx-thm:longest-path}. \qed

\noindent
We note that $(\mathcal{N}, \mathcal{A})$ is acyclic and it consists of {\color{black}$O(I^2)$} nodes and {\color{black}$O(I^3)$} arcs.
Hence, the longest-path problem, as well as $\max_{(\alpha,\beta) \in \Lambda} F(\alpha, \beta)$ can be solved in time polynomial of the problem input.
Accordingly, we are able to replace the exponentially many constraints \eqref{ref-linear-3} in the formulation of (DRNS) with {\color{black}$O(I^3)$} many linear constraints. This yields the {\color{black}following MILP} reformulation.
\begin{proposition} \label{prop:two-longest}
Under Structure \ref{assump:chained-pools} and Assumption \ref{assump:technical}, the (DRNS) model \eqref{drns} yields the same optimal objective value as the following MILP:
{\color{black}
\begin{align*}
\min \ & \ \theta + \sum_{i=1}^I \Bigl(c^{\mbox{\tiny y}}_i y^{\tinyl}_i + g_i(y^{\tinyl}_i) \lambda_i + \sum_{\ell=1}^{y^{\mbox{\tinyu}}_i - y^{\mbox{\tinyl}}_i} \bigl( c^{\mbox{\tiny y}}_i v_{i\ell} + \delta_{i\ell} \nu_{i\ell}  \bigr) \Bigr) \nonumber \\
& \ + \sum_{j=1}^I \Biggl[ \sum_{q=1}^Q \mu_{jq} \rho_{jq} + c^{\mbox{\tiny w}}_j w^{\tinyl}_j + f_j(w^{\tinyl}_j) \gamma_j + \sum_{\ell=1}^{w^{\mbox{\tinyu}}_j - w^{\mbox{\tinyl}}_j} \bigl( c^{\mbox{\tiny w}}_j u_{j\ell} + \Delta_{j\ell} \varphi_{j\ell} \bigr) \Biggr] \\
\mbox{s.t.} \ & \ \mbox{\eqref{ref-note-36}--\eqref{ref-note-33}}, \ \mbox{\eqref{ref-linear-4}}, \ \mbox{\eqref{onepool-ref-3}--\eqref{onepool-ref-4}}, \ \mbox{\eqref{disjoint-ref-1}--\eqref{disjoint-ref-2}}, \\
& \theta \geq \sum_{i=1}^I (c^{\text{\tiny r}}_i  + c^{\text{\tiny p}}_i) + \pi_{\texttt{S}} - \pi_{\texttt{T}}, \\
& \pi_{\texttt{S}} - \pi_{t_{[1]}} \geq ((\zeta^{\mbox{\tiny x}}_{11} + \eta^{\mbox{\tiny x}}_{11}) - (\zeta^{\mbox{\tiny e}}_1 + \eta^{\mbox{\tiny e}}_1)) t_{11} \quad \forall t_{[1]} \in \mathcal{B}_1, \\
& \pi_{t_{[1]}} - \pi_{(t_{[1]}, e_{2\ell})} \geq \\
& \begin{cases}
( \phi^{\mbox{\tiny x}}_{11} - \phi^{\mbox{\tiny e}}_1 ) t_{11} & \text{if } \ell = 0 \\
(\zeta^{\mbox{\tiny x}}_{2\ell} + \eta^{\mbox{\tiny x}}_{2\ell})- (\zeta^{\mbox{\tiny e}}_{2} + \eta^{\mbox{\tiny e}}_{2}) + (\phi^{\mbox{\tiny x}}_{1\ell} - \phi^{\mbox{\tiny e}}_1) & \text{if } \ell \neq 0 \\
\end{cases} 
\quad \forall t_{[1]} \in \mathcal{B}_1, \forall \ell \preceq 2, \\
& \pi_{(t_{[1]}, e_{(i-1)k})} - \pi_{(t_{[1]}, e_{i\ell})} \geq \\
& \begin{cases}
0, & \text{if } k = 0, \ell = 0\\
\phi^{\mbox{\tiny x}}_{(i-1)k} - \phi^{\mbox{\tiny e}}_{i-1}, & \text{if } k \neq 0, \ell = 0 \\
(\zeta^{\mbox{\tiny x}}_{i\ell} + \eta^{\mbox{\tiny x}}_{i\ell}) - (\zeta^{\mbox{\tiny e}}_{i} + \eta^{\mbox{\tiny e}}_{i}) + (\phi^{\mbox{\tiny x}}_{(i-1) (k \curlyvee \ell)} - \phi^{\mbox{\tiny e}}_{i-1}), & \text{if } \ell \neq 0 \\
\end{cases} \\
& \hspace{2mm}
\forall t_{[1]} \in \mathcal{B}_1, \forall i \in [3,I]_{\mathbb{Z}}, \forall k \preceq i-1, \forall \ell \preceq i, \\
& \pi_{ (t_{11}, e_{Ik})}- \pi_{\texttt{T}} \geq
\begin{cases}
(\phi^{\mbox{\tiny x}}_{I1} - \phi^{\mbox{\tiny e}}_{I}) t_{11} & \text{if } k = 0 \\
\phi^{\mbox{\tiny x}}_{Ik} - \phi^{\mbox{\tiny e}}_{I} & \text{if } k \neq 0
\end{cases}
\quad \forall t_{[1]} \in \mathcal{B}_1, \forall k \preceq I.
\end{align*}
}
\end{proposition}
\proof See Appendix \ref{apx-prop:two-longest}. \qed

\section{Optimal Nurse Pool Design} \label{sec:pool-design}
Of all the three practical nurse pool structures, {\color{black}Structure \ref{assump:one-pool}} is most flexible as every pool nurse is capable of working in all units. However, this incurs a high need for cross-training. For example, to enable a nurse working in a unit to be a pool nurse, he/she needs to be cross-trained for all the remaining $J-1$ units. As a result, enabling all nurses needs as many as $J(J-1)/2$ pairs of cross-training. In contrast, {\color{black}Structure \ref{assump:chained-pools}} needs $J$ pairs of cross-training because every pool consists of exactly two units. {\color{black}Structure \ref{assump:disjoint-pools}} needs even less cross-training if we adopt a ``sparse'' design, e.g., pooling together a small subset of units. In this section, we examine how to design a sparse but effective pool structure that is disjoint. Specifically, we search for a disjoint pool structure that needs as few cross-training as possible, while achieving a pre-specified performance guarantee in terms of {\color{black} the worst-case expected total staffing cost, i.e., optimal objective value of (DRNS)}.\footnote{We notice that there exist multiple alternative quantities that can be used to quantify the effort of cross-training. In this paper, we pick the number of pairs of cross-training as a representative objective function. Alternative objectives can be similarly modeled and computed.} To this end, we define binary variables $a_{ij}$ such that $a_{ij}=1$ if unit $j$ is assigned to pool $i$ and $a_{ij}=0$ otherwise, binary variables $o_i$ such that $o_i=1$ if any units are assigned to pool $i$ (i.e., if pool $i$ is ``opened'') and $o_i = 0$ otherwise, and binary variables $p_{jk}$ such that $p_{jk}=1$ if units $j$ and $k$ are assigned to the same pool and $p_{jk}=0$ otherwise. Then, the total amount of needed cross-training equals $\sum_{j=1}^J \sum_{k=j+1}^J p_{jk}$. In addition, these binary variables satisfy the following constraints:
\begin{subequations}
\begin{align}
& \ \sum_{i = 1}^{I+1} a_{ij} = 1, \ \ \forall j \in [J], \label{onpd-note-1} \\
& \ a_{ij} \leq o_i, \ \ \forall i \in [I], \ \forall j \in [J], \label{onpd-note-2} \\
& \ p_{jk} \geq a_{ij} + a_{ik} - 1, \ \ \forall i \in [I], \ \forall j, k \in [J] \mbox{ and } j < k, \label{onpd-note-3} 
\end{align}
\end{subequations}
where constraints \eqref{onpd-note-1} designate that each unit is assigned to exactly one pool (we create a dummy pool $I+1$ that collects all units that are not covered by any existing pools), constraints \eqref{onpd-note-2} ensure that no units can be assigned to a pool if it is not opened, and constraints \eqref{onpd-note-3} designate that $p_{jk}=1$ if there is a pool $i$ such that $a_{ij} = a_{ik} = 1$. If no such a pool $i$ exists, then constraints \eqref{onpd-note-2} {\color{black}and \eqref{onpd-note-3}} reduce to $p_{jk} \geq 0$ and $p_{jk}$ equals zero at optimality due to the objective function \eqref{onpd-obj}. Based on Proposition \ref{prop:exact-disjoint}, the optimal nurse pool design (OPD) model is formulated as
{\color{black}
\begin{subequations}
\label{onpd}
\begin{align}
(\mbox{\bf OPD}): \ \ \min \ & \ \sum_{j=1}^J \sum_{k=j+1}^J p_{jk} \label{onpd-obj} \\
\mbox{s.t.} \ & \ \mbox{\eqref{ref-note-36}--\eqref{ref-note-33}}, \ \mbox{\eqref{ref-linear-4}}, 
\ \mbox{\eqref{onepool-ref-3}}, \ \mbox{\eqref{onepool-ref-4}}, \ \mbox{\eqref{disjoint-ref-1}},
\ \mbox{\eqref{onpd-note-1}--\eqref{onpd-note-3}}, \\
& \ \sum_{i=1}^I \Bigl( \theta_i + c^{\mbox{\tiny y}}_i y^{\tinyl}_i + g_i(y^{\tinyl}_i) \lambda_i + \sum_{\ell=1}^{y^{\mbox{\tinyu}}_i - y^{\mbox{\tinyl}}_i} \bigl(c^{\mbox{\tiny y}}_i v_{i\ell} + \delta_{i\ell} \nu_{i\ell}  \bigr) \Bigr) o_i \nonumber \\
& \ + \sum_{j=1}^J \Biggl[ \sum_{q=1}^Q \mu_{jq} \rho_{jq} + c^{\mbox{\tiny w}}_j w^{\tinyl}_j + f_j(w^{\tinyl}_j) \gamma_j + \sum_{\ell=1}^{w^{\mbox{\tinyu}}_j - w^{\mbox{\tinyl}}_j} \bigl( c^{\mbox{\tiny w}}_j u_{j\ell} + \Delta_{j\ell} \varphi_{j\ell} \bigr) \Biggr] \ \leq \ T, \label{onpd-con-target} \\
& \theta_i \geq \phi^{\mbox{\tiny e}}_i + \sum_{j =1}^J (\zeta^{\mbox{\tiny e}}_j + \eta^{\mbox{\tiny e}}_j) a_{ij}, \ \forall i \in [I+1], \label{onpd-con-1} \\
& \theta_i \geq \sum_{j =1}^J (\zeta^{\mbox{\tiny e}}_j + \eta^{\mbox{\tiny e}}_j)a_{ij} + \Big\{ \phi^{\mbox{\tiny x}}_{ik^*} + (\zeta^{\mbox{\tiny x}}_{j^* k^*} + \eta^{\mbox{\tiny x}}_{j^* k^*}) - (\zeta^{\mbox{\tiny e}}_{j^*} + \eta^{\mbox{\tiny e}}_{j^*}) \Big\} a_{ij^*} + \sum_{j \in [J] \setminus \{j^*\} } \chi_{jk^*}a_{ij}, \nonumber \\
& \hspace{9mm} \forall i \in [I+1], \forall j^* \in [J], \forall k^* \in [j^*], \label{onpd-con-2}  \\
& \chi_{jk^*} \geq  (\zeta^{\mbox{\tiny x}}_{jk} + \eta^{\mbox{\tiny x}}_{jk}) - (\zeta^{\mbox{\tiny e}}_j + \eta^{\mbox{\tiny e}}_j), \ \forall j \in [J], \forall k^* \in [J], \forall k \in [\min\{j,k^*\}], \\
& \hspace{-0.15cm} \left.
\begin{array}{l}
\Big\{ \phi^{\mbox{\tiny x}}_{ik} + ( c^{\mbox{\tiny x}}_k + \lambda_i) y^{\tinyl}_i + \sum_{\ell=1}^{y^{\mbox{\tinyu}}_i - y^{\mbox{\tinyl}}_i}  (c^{\mbox{\tiny x}}_k v_{i\ell} + \nu_{i\ell} ) \Big\} a_{ij} \geq 0,  \\ [0.1cm]
\phi^{\mbox{\tiny x}}_{ik} a_{ij} \geq 0, 
\end{array}\right\} \ \ \forall i \in [I+1], \forall j \in [J],\forall k\in [j], \label{onpd-con-3}
\end{align}
\end{subequations}
}
where constraint \eqref{onpd-con-target} ensures that the {\color{black} worst-case expected total staffing cost} does not exceed a given target $T$. If $y^{\tinyl}_i = 0$ for all $i \in [I]$, i.e., if there is no minimum staffing requirement for pool nurses, then we shall pick $T$ from the interval {\color{black}$[Z^{\star}_{1}, Z^{\star}_{0}]$}, where {\color{black}$Z^{\star}_{1}$} represents the {\color{black} worst-case expected total staffing cost} with maximum flexibility and {\color{black}$Z^{\star}_{0}$} represents that with minimum flexibility. By gradually decreasing this target from {\color{black}$Z^{\star}_{0}$ to $Z^{\star}_{1}$}, the amount of cross-training grows and accordingly we obtain a cost-training frontier that can clearly illustrate the trade-off between these two performance measures (see Section \ref{sec:results_poolstructure} for the numerical demonstration).

{\color{black} Formulation \eqref{onpd} poses computational challenges due to the bilinear terms in constraints \eqref{onpd-con-target}--\eqref{onpd-con-3}. To solve \eqref{onpd} more effectively, we recast it as a MILP in the following proposition.
\begin{proposition} \label{prop:opd}
Under Assumption \ref{assump:technical}, the (OPD) model \eqref{onpd} yields the same optimal objective value and the same set of optimal solutions as the following MILP:
\begin{align*}
\min \ & \ \sum_{j=1}^J \sum_{k=j+1}^J p_{jk} \\
\mbox{s.t.} \ & \ \mbox{\eqref{ref-note-36}--\eqref{ref-note-33}}, \ \mbox{\eqref{ref-linear-4}}, 
\ \mbox{\eqref{onepool-ref-3}}, \ \mbox{\eqref{onepool-ref-4}}, \ \mbox{\eqref{disjoint-ref-1}}, 
\ \mbox{\eqref{onpd-note-1}--\eqref{onpd-note-3}}, \\
& \ \sum_{i=1}^I \Bigl( \theta_i + c^{\mbox{\tiny y}}_i y^{\tinyl}_i o_i + g_i(y^{\tinyl}_i) \lambda_i + \sum_{\ell=1}^{y^{\mbox{\tinyu}}_i - y^{\mbox{\tinyl}}_i} \bigl( c^{\mbox{\tiny y}}_i v_{i\ell} + \delta_{i\ell} \nu_{i\ell} \bigr) \Bigr) \\
& \ + \sum_{j=1}^J \Biggl[ \sum_{q=1}^Q \mu_{jq} \rho_{jq} + c^{\mbox{\tiny w}}_j w^{\tinyl}_j + f_j(w^{\tinyl}_j) \gamma_j + \sum_{\ell=1}^{w^{\mbox{\tinyu}}_j - w^{\mbox{\tinyl}}_j} \bigl( c^{\mbox{\tiny w}}_j u_{j\ell} + \Delta_{j\ell} \varphi_{j\ell} \bigr) \Biggr] \ \leq \ T, \\
& \theta_i \geq  \phi^{\mbox{\tiny e}}_i + \sum_{j =1}^J (\zeta^{\mbox{\tiny e}}_{ij} + \eta^{\mbox{\tiny e}}_{ij}), \ \forall i \in [I], \\
& \theta_i \geq  \sum_{j =1}^J (\zeta^{\mbox{\tiny e}}_{ij} + \eta^{\mbox{\tiny e}}_{ij}) + \Big\{ \phi^{\mbox{\tiny x}}_{ij^*k^*} + (\zeta^{\mbox{\tiny x}}_{ij^* k^*} + \eta^{\mbox{\tiny x}}_{ij^* k^*}) - (\zeta^{\mbox{\tiny e}}_{ij^*} + \eta^{\mbox{\tiny e}}_{ij^*}) \Big\} +  \sum_{j \in [J] \setminus \{j^*\} } \chi_{ijk^*}, \\
& \hspace{9mm} \forall i \in [I], \forall j^* \in [J], \forall k^* \in [j^*], \\
& \chi_{ijk^*} \geq (\zeta^{\mbox{\tiny x}}_{ij k} + \eta^{\mbox{\tiny x}}_{ij k}) - (\zeta^{\mbox{\tiny e}}_{ij} + \eta^{\mbox{\tiny e}}_{ij}), \ \forall i \in [I+1], \forall j \in [J], \forall k^* \in [J], \forall k \in [\min \{j,k^*\}],  \\
& \hspace{-0.15cm} \left.
\begin{array}{l}
\phi^{\mbox{\tiny x}}_{ijk} + ( c^{\mbox{\tiny x}}_k + \lambda^{\text{\tiny a}}_{ij}) y^{\tinyl}_i + \sum_{\ell=1}^{y^{\mbox{\tinyu}}_i - y^{\mbox{\tinyl}}_i}  (c^{\mbox{\tiny x}}_k v^{\text{\tiny a}}_{ij\ell} + \nu^{\text{\tiny a}}_{ij\ell} )  \geq 0,  \\ [0.1cm]
\phi^{\mbox{\tiny x}}_{ijk} \geq 0, 
\end{array}\right\} \ \ \forall i \in [I+1], \forall j \in [J],\forall k\in [j], \\
& \hspace{-0.15cm} \left.
\begin{array}{l}
\zeta^{\text{\tiny e}}_j = \sum_{i=1}^{I+1} \zeta^{\text{\tiny e}}_{ij}, \ \ 
\eta^{\text{\tiny e}}_j = \sum_{i=1}^{I+1} \eta^{\text{\tiny e}}_{ij},  \\ [0.1cm]
\zeta^{\text{\tiny x}}_{jk} = \sum_{i=1}^{I+1} \zeta^{\text{\tiny x}}_{ijk}, \ \ 
\eta^{\text{\tiny x}}_{jk} = \sum_{i=1}^{I+1} \eta^{\text{\tiny x}}_{ijk}, \ \forall k \in [j], \\[0.1cm]
\end{array}\right\} \ \ \forall j \in [J], \\
& \hspace{-0.15cm} \left.
\begin{array}{l}
0 \leq \zeta^{\text{\tiny e}}_{ij} \leq K a_{ij}, \ \ -K a_{ij} \leq \eta^{\text{\tiny e}}_{ij} \leq K a_{ij},  \\ [0.1cm]
0 \leq \zeta^{\text{\tiny x}}_{ijk} \leq K a_{ij}, \ \ -K a_{ij} \leq \eta^{\text{\tiny x}}_{ijk} \leq K a_{ij}, \ \forall k \in [j], \\[0.1cm]
0 \leq \phi^{\text{\tiny x}}_{ijk} \leq K a_{ij}, \ \ \phi^{\text{\tiny x}}_{ik} - K (1-a_{ij}) \leq \phi^{\text{\tiny x}}_{ijk} \leq \phi^{\text{\tiny x}}_{ik} + K (1-a_{ij}), \ \forall k \in [j], \\[0.1cm]
- c^{\text{\tiny x}}_{J} a_{ij} \leq \lambda^{\text{\tiny a}}_{ij} \leq 0, \ \ \lambda_i \leq \lambda^{\text{\tiny a}}_{ij} \leq \lambda_i + c^{\text{\tiny x}}_{J} (1- a_{ij}), \\ [0.1cm]
0 \leq v^{\text{\tiny a}}_{ij \ell} \leq a_{ij}, \ \ v_{i \ell} - (1-a_{ij}) \leq v^{\text{\tiny a}}_{ij \ell} \leq v_{i \ell} + (1-a_{ij}), \ \forall \ell \in [y^{\text{\tiny U}}_i - y^{\text{\tiny L}}_i], \\ [0.1cm]
-c^{\text{\tiny x}}_{J} a_{ij} \leq \nu^{\text{\tiny a}}_{ij \ell} \leq c^{\text{\tiny x}}_{J} a_{ij}, \ \forall \ell \in [y^{\text{\tiny U}}_i - y^{\text{\tiny L}}_i], \\ [0.1cm]
\nu_{i \ell} - c^{\text{\tiny x}}_{J}(1-a_{ij}) \leq \nu^{\text{\tiny a}}_{ij \ell} \leq \nu_{i \ell} + c^{\text{\tiny x}}_{J}(1-a_{ij}), \ \forall \ell \in [y^{\text{\tiny U}}_i - y^{\text{\tiny L}}_i], \\ [0.1cm]
\end{array}\right\} \ \ \forall i \in [I+1], \forall j \in [J], \\
& -c^{\mbox{\tiny x}}_{J} o_i \leq \lambda_i \leq 0 , \ \ v_{i1} \leq o_i, \ \forall i \in [I],
\end{align*}
where $K$ represents a sufficiently large positive constant.
\end{proposition}
\proof See Appendix \ref{apx-prop:opd}. \qed}
In addition, the above formulation involves symmetric binary solutions. In Appendix \ref{apx-symmetry-breaking}, we derive symmetry breaking inequalities to enhance its computational efficacy.

\section{Numerical Experiments} \label{sec:results}
In this section, we report numerical experiments on (DRNS) and (OPD) models. We summarize our main findings as follows:
\begin{enumerate}[1.]
\item Under the practical nurse pool structures as introduced in Section \ref{sec:tractable}, the MILP reformulations of (DRNS) lead to significant speed-up over the separation algorithm.
\item {\color{black}Staffing decisions produced by (DRNS) provide better out-of-sample performance than those produced by the stochastic programming (SP) approach.}
\item Modeling nurse absenteeism improves the out-of-sample performance of staffing decisions. The improvement becomes more significant as the value of operational flexibility increases.
\item Even a very sparse nurse pool design can harvest most of the operational flexibility.
\item An optimal nurse pool design tends to pool together the units with higher variability, e.g., higher standard deviation of nurse demand and/or higher absence rate. In particular, the variability of nurse absenteeism plays a more important role in optimal pool design.
\end{enumerate}
In all experiments, we solve optimization models by {\color{black}GUROBI 9.0.1} via Python 2.7 on a personal laptop with an Intel(R) Core(TM) i7-4850HQ CPU@2.3GHz and 16GB RAM. 
{\color{black} The number of threads is set to be $1$.}

\begin{table}[H]
	\centering
	\caption{Wall-clock seconds taken to solve (DRNS) by separation algorithm (Sep), separation algorithm with valid inequalities (Sep$_{\text{\tiny VI}}$), and MILP reformulations in Section \ref{sec:tractable}.}
	\label{Table:Computation}	
	\begin{threeparttable}	
		\begin{tabular}{|c||r|r||r|r|r||r|r|r||r|r|r|}	
			\hline
			& \multicolumn{2}{c||}{{General\tnote{a}}}  & \multicolumn{3}{c||}{{Structure~\ref{assump:one-pool}}} & \multicolumn{3}{c||}{{Structure~\ref{assump:disjoint-pools}\tnote{b}}} & \multicolumn{3}{c|}{{Structure~\ref{assump:chained-pools}}} \\ \hline	
			$J$ & \multicolumn{1}{c|}{Sep}  	& \multicolumn{1}{c||}{Sep$_{\text{\tiny VI}}$}  & \multicolumn{1}{c|}{Sep}  	& \multicolumn{1}{c|}{Sep$_{\text{\tiny VI}}$} 	& \multicolumn{1}{c||}{MILP} & \multicolumn{1}{c|}{Sep}  	& \multicolumn{1}{c|}{Sep$_{\text{\tiny VI}}$} 	& \multicolumn{1}{c||}{MILP} & \multicolumn{1}{c|}{Sep}  	& \multicolumn{1}{c|}{Sep$_{\text{\tiny VI}}$} 	& \multicolumn{1}{c|}{MILP} \\ \hline
			$5$ 	& $0.41$ & $0.27$  	& $0.22$	& $0.20$ 	& $0.18$ 	
			& $0.23$	& $0.22$ 	& $0.16$ & $0.37$	& $0.26$ 	& $0.25$ 			 \\
			$10$ 	& $4.02$ & $0.93$ 	& $0.57$  	& $0.46$  	& $0.52$  	
			& $0.71$  	& $0.40$  	& $0.22$ & $23.15$  & $0.83$  	& $0.67$ 			 \\
			$50$ 	& $285.19$ & $50.93$ & $81.27$	& $15.72$  	& $17.20$ 	
			& $165.85$	& $44.62$  	& $8.22$ & LIMIT\tnote{c} & $3187.35$ 		& $1078.28$ \\ \hline
		\end{tabular}
		\begin{tablenotes}\footnotesize			
			\item[a] $P_1 = \{ j \in [J] : j \text{ is even} \}$, $P_2 = \{ j \in [J] : j \text{ is odd} \}$, $P_3 = \{ j \in [J] : j \leq J/2 \}$, and $P_4 = \{ j \in [J] : j > J/2 \}$.
			\item[b] $P_i=\{4i-3, 4i-2, 4i-1 ,4i\}, \ \forall i \in [\lfloor J/4 \rfloor]$, and $P_{\lfloor J/4 \rfloor +1}=\{4 \lfloor J/4 \rfloor + 1, \ldots, J\}$.
			\item[c] Time Limit = 2 hours. 			
		\end{tablenotes}
	\end{threeparttable}
\end{table}

\subsection{Computational efficacy}\label{sec:results_performance}
{\color{black}We evaluate the computational efficacy of three approaches on solving (DRNS) under general and special pool structures: the separation algorithm describe in Algorithm \ref{algo-sep} (denoted by Sep), the separation algorithm with valid inequalities \eqref{VI-general} (denoted by Sep$_{\text{\tiny VI}}$), and the MILP reformulations derived in Section \ref{sec:tractable} in case of special pool structures (denoted by MILP). We design test instances using data and insights obtained from our collaborating hospital. In particular, we obtain nurse demand and attendance data for 5 hospital units during July-2010 and June-2014. Using this data, we calibrate the ambiguity set $\mathcal{D}(w, y)$ and other problem parameters (see Appendix \ref{apx-5-unit} for details). In addition, we duplicate this 5-unit system to create 10- and 50-unit systems. For each system, we consider four different pool structures, including general pool and Structures~\ref{assump:one-pool}, \ref{assump:disjoint-pools}, and \ref{assump:chained-pools}. We report the computing time (in wall-clock seconds) in Table \ref{Table:Computation}. From this table, we observe that valid inequalities \eqref{VI-general} significantly improve the computational efficacy of the separation algorithm. For example, the computing time of Sep$_{\text{\tiny VI}}$ is a small fraction of that of Sep in most instances, and Sep$_{\text{\tiny VI}}$ solves the most challenging instance ($J = 50$ under Structure \ref{assump:chained-pools}) to optimum within the time limit, while Sep cannot. In addition, MILP prevails under special pool structures. For example, it solves all instances to global optimum in 18 minutes. For completeness of this evaluation, we also demonstrate that Sep is significantly more effective than solving the exponentially-sized reformulation \eqref{drns-milp-exp}. The results are reported in Appendix \ref{apx-seperation}.}


\begin{table}[H]
	\centering
	\caption{Out-of-sample comparison between (DRNS) and SP.}
	\label{Table:SP_DRO_OOS}	
	\begin{threeparttable}	
		\begin{tabular}{|l||C{1.6cm}|C{1.6cm}|C{1.6cm}||C{2,5cm}|c|}
			\hline 
			& \multicolumn{3}{c||}{Average staffing cost} & \multicolumn{2}{c|}{Average number of temp nurses hired} \\ \cline{2-6}
			&  SP & DRNS  & GAP(\%)\tnote{d} & SP & DRNS  \\ \hline
			Structure 0			& $34805$  & $33996$& $2.38$ & $16.94$ & $18.87$  \\
			Structure~\ref{assump:one-pool}			& $27720$  & $24489$& $13.19$ & $4.54$ & $3.32$  \\
			Structure~\ref{assump:disjoint-pools} & $31891$  & $30228$& $5.50$ & $11.56$ & $12.65$  \\
			Structure~\ref{assump:chained-pools}			& $27768$  & $24556$& $13.08$ & $4.59$ & $3.39$  \\
			General	& $27772$  & $24580$& $12.99$ & $4.59$ & $3.41$  \\  \hline
		\end{tabular}
		\begin{tablenotes}\footnotesize
			\item[d] GAP$=(Z^{\text{\tiny SP}}-Z^{\text{\tiny DRNS}})/Z^{\text{\tiny DRNS}} \times 100\%$, where $Z^{\text{\tiny SP}}$ and $Z^{\text{\tiny DRNS}}$ are average staffing costs of SP and (DRNS) in the out-of-sample simulation, respectively.
		\end{tablenotes}
	\end{threeparttable}	
\end{table}

\subsection{Value of modeling endogenous absenteeism} \label{sec:results_SP_DRO}
{\color{black}
We demonstrate the value of modeling endogenous absenteeism by comparing the out-of-sample performance of DRNS and a stochastic programming (SP) approach, which models nurse absenteeism as exogenous uncertainty. Specifically, the SP approach assumes a constant nurse absenteeism rate, independent of the staffing levels (see Appendix \ref{apx-SPformulation} for the formulation).

We report the out-of-sample average staffing costs of DRNS and SP in Table \ref{Table:SP_DRO_OOS}. From this table, we observe that DRNS leads to a lower staffing cost on average and tends to need less temporary nurses than SP. Under a general pool structure, as well as under Structures~\ref{assump:one-pool} and \ref{assump:chained-pools}, the gap in out-of-sample performance exceeds 10\%. This demonstrates the value of modeling absenteeism as endogenous uncertainty. To examine this observation further, we compare the optimal staffing levels $(w^*, y^*)$ of DRNS and SP in Table \ref{Table:SP_DRO_Staffing}. From this table, we observe that DRNS staffs more pool nurses and assigns less unit nurses. This is because the nurses become less likely to show up as the staffing level increases (see Figures~\ref{Fig:segmented} and \ref{Fig:Showuprates}). The increasing absenteeism rate discourages DRNS to staff more unit nurses. In contrast, the SP approach assumes an exogenous nurse absenteeism and over-staffs unit nurses, resulting in out-of-sample disappointment.

\begin{table}[H]
	\centering
	\caption{Optimal staffing levels produced by SP and DRNS.}
	\label{Table:SP_DRO_Staffing}	
	\begin{threeparttable}
		\begin{tabular}{|l||c|c||c|c|} 
			\hline 
			& \multicolumn{2}{c||}{SP}  & \multicolumn{2}{c|}{DRNS}  \\ \cline{2-5}
			& $w^*$ & $y^*$ & $w^*$ & $y^*$ \\ \hline 
			Structure 0			& $12, 13, 17, 13, 16$ & $0$  				& $10,11, 14, 11, 14$ & $0$  \\
			Structure 1			& $12, 13, 17, 13, 16$ & $13$  				& $9, 10, 14, 8, 12$ & $19$  \\
			Structure D& $12, 13, 17, 13, 16$ & $6$  				& $9, 11, 14, 11, 14$ & $7$  \\
			Structure C			& $12, 13, 17, 13, 16$ & $3, 3, 3, 2, 2$  	& $9, 10, 14, 8, 12$ & $2, 2, 5, 5, 5$  \\
			General	& $12, 13, 17, 13, 16$ & $7, 6$  			& $9, 10, 14, 8, 12$ & $7, 12$  \\ \hline
		\end{tabular}
	\end{threeparttable}
\end{table}
}

\subsection{Value of modeling nurse absenteeism} \label{sec:results_simulation}

As discussed in Section \ref{sec:intro}, modeling nurse absenteeism incurs endogenous uncertainty and computational challenges. It is hence worth examining what (DRNS) buys us, i.e., the value of modeling nurse absenteeism.
{\color{black}
To this end, we generate $400$ instances by tuning parameters of the 5-unit system under Structure 1. Specifically, we (i) generate random $c^{\text{\tiny x}}_j$ from the interval $[450, 1500]$, (ii) perturb $\mu_{j1}$ in Table \ref{Table:5Unit-Parameters} by adding a random noise $\tilde{\alpha}_j$ supported on the interval $[-1, 1]$, and (iii) generate random $\text{sd}_j$ from the interval $[1,2]$.}
In addition, we consider a variant of (DRNS) that overlooks the nurse absenteeism, in which we assume that all assigned nurses show up. Then, we compare the out-of-sample performance of the optimal nurse staffing decisions produced by (DRNS) and that produced by overlooking absenteeism. Fixing the nurse staffing levels as in a (DRNS) optimal solution $(w^*, y^*)$, we generate a large number of scenarios for nurse demand and absenteeism. Exposing $(w^*, y^*)$ under these scenarios produces an out-of-sample estimate of the average staffing cost with absenteeism, which we denote by $Z_{\mbox{\tiny abs}}$. Using the same set of scenarios, we examine the optimal solution produced by overlooking absenteeism and obtain an out-of-sample average cost without absenteeism, denoted by $Z_{\mbox{\tiny w/o}}$. Using the same out-of-sample procedures, we compute the average number of temporary nurses hired when considering absenteeism (denoted by $x_{\mbox{\tiny abs}}$) and when overlooking it (denoted by $x_{\mbox{\tiny w/o}}$).

\begin{figure}[H]
	\centering
	\begin{subfigure}[b]{0.4\textwidth}
		\includegraphics[width=\textwidth]{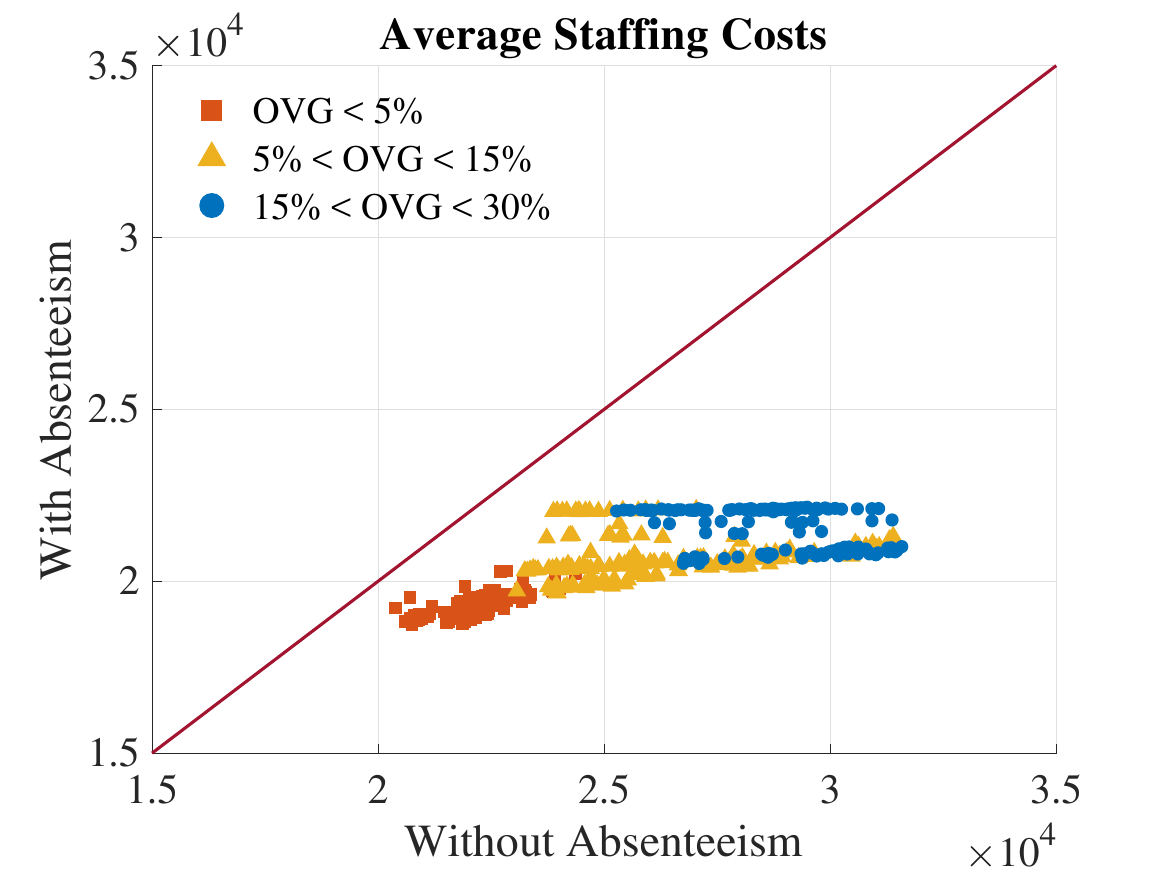}
		\caption{Average staffing costs}
		\label{Fig: absenteeism-a}
	\end{subfigure} \hspace{-0.0cm}
	\begin{subfigure}[b]{0.4\textwidth}
		\includegraphics[width=\textwidth]{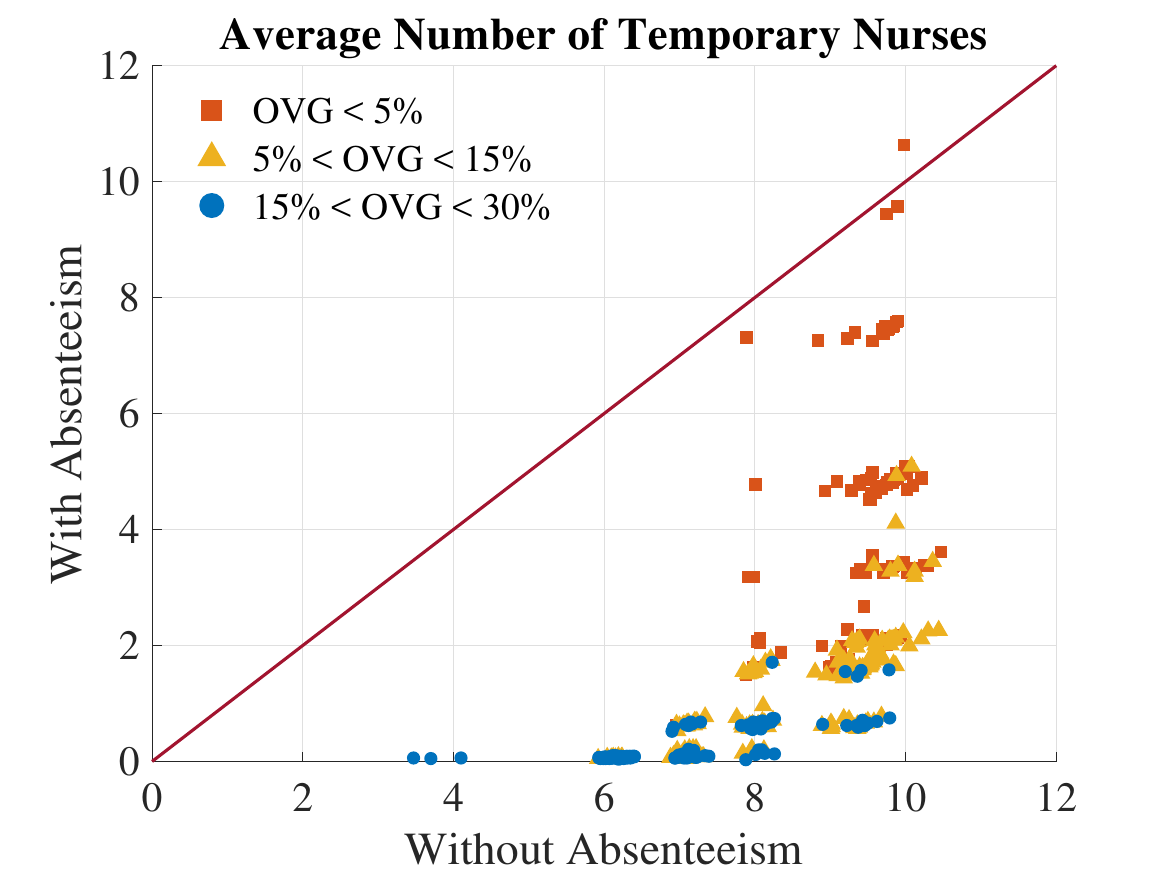}
		\caption{Average number of temporary nurses}
		\label{Fig: absenteeism-b}
	\end{subfigure}
	\caption{An out-of-sample performance comparison of considering versus overlooking nurse absenteeism. {\color{black}Note that $\mbox{OVG} := (Z^{\star}_{0} - Z^{\star}_{1})/Z^{\star}_{0} \times 100\%$, where $Z^{\star}_{0}$ (respectively, $Z^{\star}_{1}$) is the optimal value of (DRNS) with no nurse pools (respectively, under Structure 1).}}
	\label{Fig: absenteeism}
\end{figure}

We depict the values of $Z_{\mbox{\tiny w/o}}$ ($x$-coordinate) and $Z_{\mbox{\tiny abs}}$ ($y$-coordinate) obtained in {\color{black}400} replications in Figure~\ref{Fig: absenteeism-a}. From this figure, we observe that most dots are below the 45-degree line, indicating that $Z_{\mbox{\tiny w/o}} - Z_{\mbox{\tiny abs}} > 0$, i.e., modeling nurse absenteeism yields nurse staffing levels with better out-of-sample performance. In addition, we group the dots based on the relative value of operational flexibility {\color{black}$\mbox{OVG} := (Z^{\star}_{0} - Z^{\star}_{1})/Z^{\star}_{0} \times 100\%$.}
From Figure~\ref{Fig: absenteeism-a}, we observe that the difference $Z_{\mbox{\tiny w/o}} - Z_{\mbox{\tiny abs}}$ shows an increasing trend as OVG increases. That is, modeling nurse absenteeism becomes more valuable as the value of operational flexibility increases. This makes sense because when a unit is short of supply due to nurse absenteeism, making it up with pool nurses {\color{black} is less expensive} than doing so with temporary nurses. As a result, setting up nurse pools can effectively mitigate the impacts of nurse absenteeism. In Figure~\ref{Fig: absenteeism-b}, we depict the values of $x_{\mbox{\tiny w/o}}$ and $x_{\mbox{\tiny abs}}$ obtained in the {\color{black}400} replications and make similar observations.

\begin{figure}[H]
	\centering
	\begin{subfigure}[b]{0.4\textwidth}
		\includegraphics[width=\textwidth]{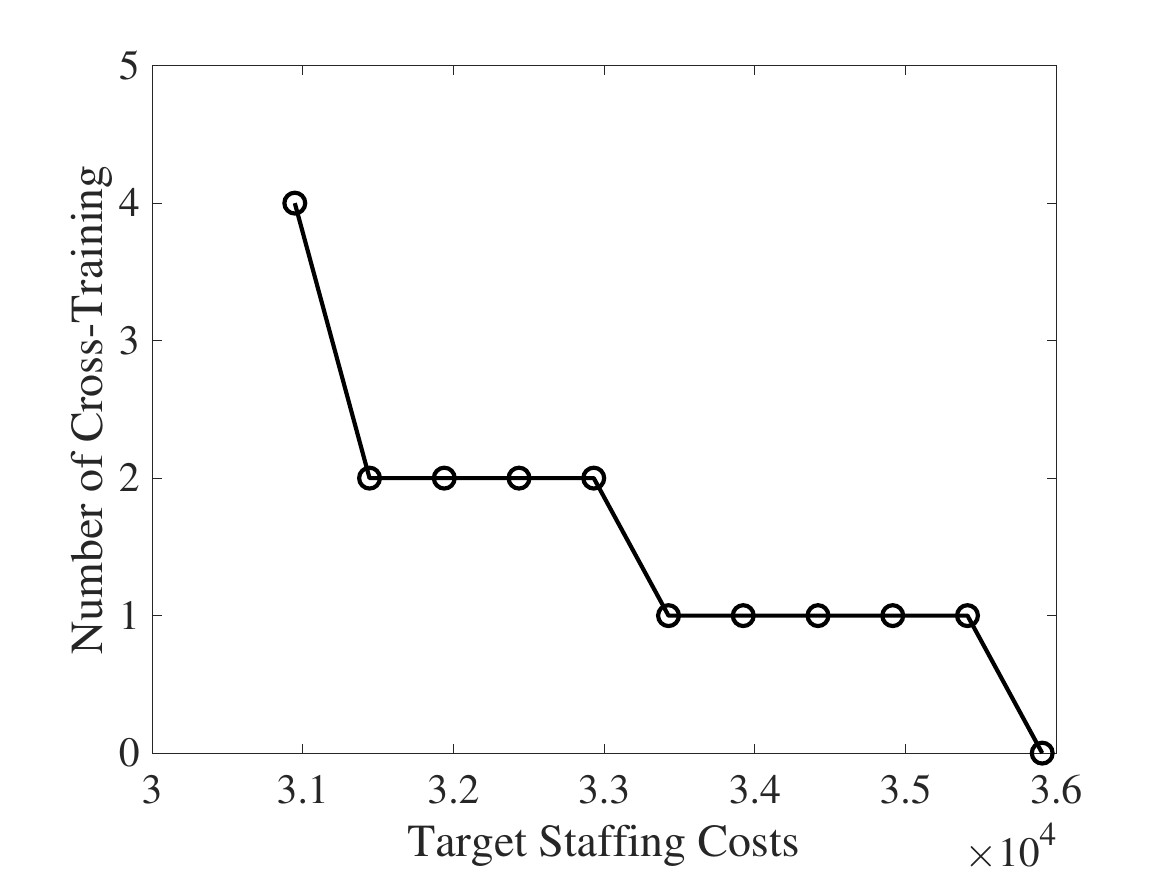}
		\caption{5-unit system}
	\end{subfigure} \hspace{-0.0cm}
	\begin{subfigure}[b]{0.4\textwidth}
		\includegraphics[width=\textwidth]{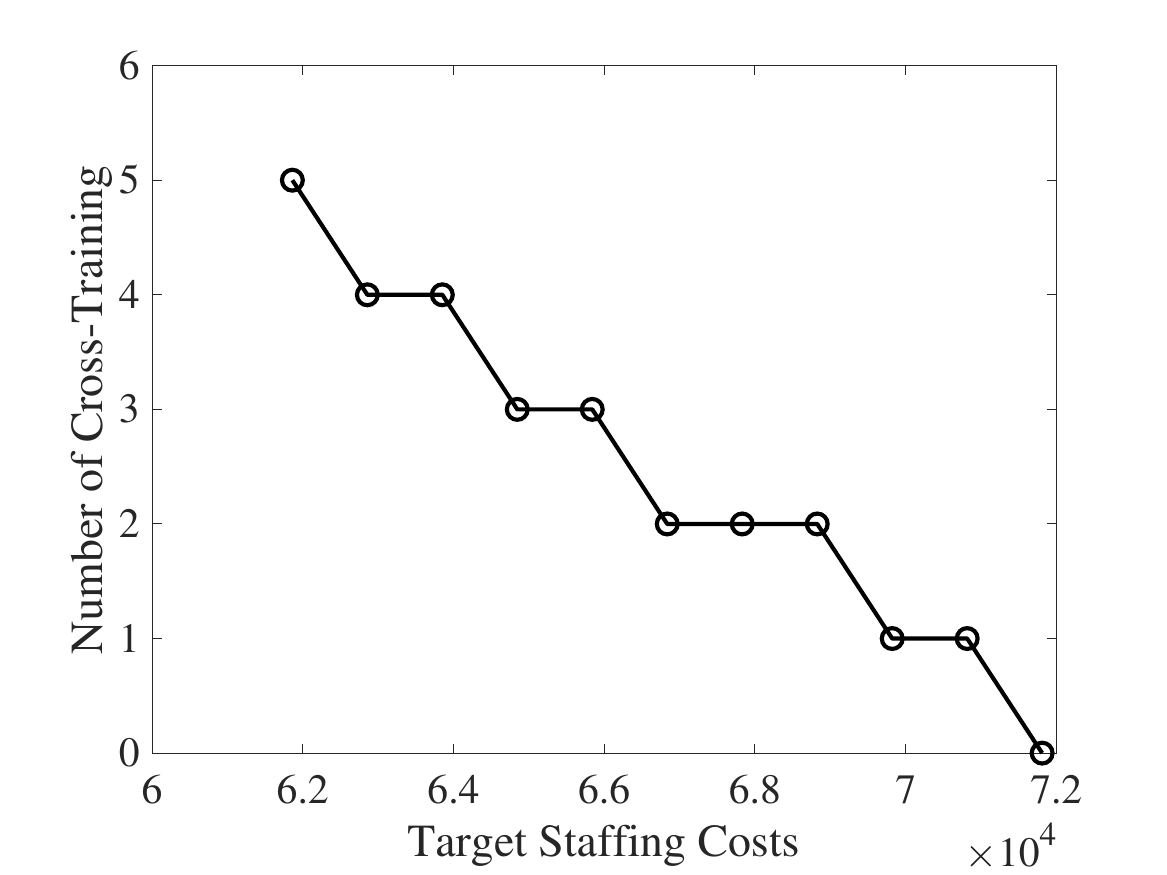}
		\caption{10-unit system}
	\end{subfigure}
    \caption{Amount of cross-training as a function of target operational cost.} \label{Fig:Number of Cross}
\end{figure}

\subsection{Sparse nurse pool design} \label{sec:results_poolstructure}

{\color{black}In several application domains (e.g., production systems~\cite{jordan1995principles,chen2015optimal}), it has been observed that a sparse design can already harvest most of the operational flexibility.} To verify this intuition in nurse staffing, {\color{black}we conduct sensitivity analysis on the target staffing cost $T$ in (OPD).} Specifically, we pick ten values of $T$ between {\color{black}$Z^{\star}_{0}$} (i.e., the optimal value of (DRNS) with no nurse pools) and {\color{black}$Z^{\star}_{1}$} (i.e., the optimal value of (DRNS) under Structure 1) uniformly. 
For each value of $T$, we solve (OPD) to obtain the minimum amount of cross-training \#($T$) that guarantees that the {\color{black} worst-case expected staffing cost} is no larger than $T$. 
We report \#($T$) as a function of $T$ in Figure~\ref{Fig:Number of Cross}, 
{\color{black}where the left figure depicts \#($T$) of the 5-unit system and the right figure depicts \#($T$) of the 10-unit system. We observe that in order to achieve $T=Z^{\star}_{1}$ (i.e., harvesting all operational flexibility of having nurse pools), the numbers of cross-training required by the 5-unit and 10-unit systems are as small as $4$ and $5$, respectively. This demonstrates the intuition that even sparse nurse pools lead to high flexibility.}

\subsection{Patterns of the optimal nurse pool design}\label{sec:results_ONPD}
{\color{black}Motivated by the sparseness intuition, we study the patterns of the optimal nurse pool design. To this end, we generate three sets of random instances of (OPD) based on the 10-unit system. In each set, we randomly select 4 units to have higher level of variability (denoted by $H$) and the remaining 6 units to have lower variability (denoted by $L$), as displayed in Table~\ref{Table:3 cases}.

\begin{table}[H]
	\caption{Units with different levels of variability} \label{Table:3 cases}
	\centering
	\begin{tabular}{|c||c|c||c|c|}
		\hline
		& \multicolumn{2}{c||}{Units in $H$} & \multicolumn{2}{c|}{Units in $L$}\\ \hline
		& \makecell{Standard deviation \\ of nurse demand} & Absence rate & \makecell{Standard deviation \\ of nurse demand} & Absence rate \\ \hline
		Set 1	& High & High  & Low & High \\ \hline
		Set 2	& Low  & High  & Low & Low  \\ \hline
		Set 3	& High & High  & Low & Low  \\ \hline
	\end{tabular}
\end{table}
\noindent We specify what ``High'' and ``Low'' in Table~\ref{Table:3 cases} represent as follows:
\begin{itemize}
\item For low standard deviation of nurse demand, we use the sd$_j$ value in Table \ref{Table:5Unit-Parameters}; and for high standard deviation of nurse demand, we randomly generate sd$_j$ from the interval $[3, 4]$.
\item For low absenteeism rate, we set $f_j(w_j) := (1-A^{\text{\tiny u}}_j) w_j$, where $A^{\text{\tiny u}}_j$ is randomly extracted from the interval $[0.10, 0.16]$; and for high absenteeism rate, we use $f_j(w_j)$ in Table \ref{Table:5Unit-Parameters}.
\end{itemize}
In addition, we set $T := Z^{\star}_1$ in (OPD), i.e., we search for sparse pool structures that produce the same operational flexibility as Structure~\ref{assump:one-pool}.
\begin{figure}[H]
	\centering 
	\includegraphics[scale=0.4]{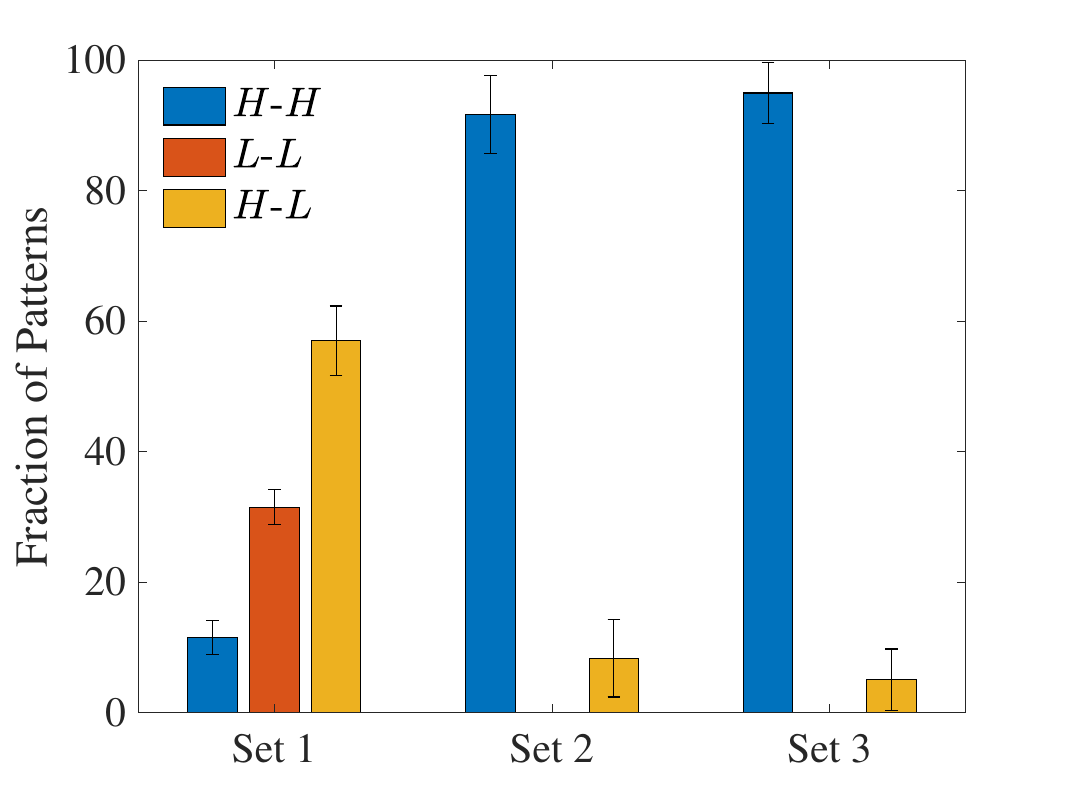}
	\caption{Patterns of optimal pool design. The height of a bar represents the fraction of a pattern appearing in an optimal pool design, and the error bar represents the corresponding 90\% confidence interval.} \label{Fig:CI}
\end{figure}
We generate 40 instances for each set and report the patterns of optimal pool design in Figure~\ref{Fig:CI}, in which we call a pool ``$H$-$H$'' if all of its units come from $H$, ``$L$-$L$'' if all of its units come from $L$, and ``$H$-$L$'' if the units come from both $H$ and $L$. From this figure, we observe that $L$-$L$ pools and $H$-$L$ pools appear frequently in an optimal pool design in instance set 1, while $H$-$H$ pools become dominant in instance sets 2 and 3. Comparing the results of instance set 1 with those of set 3, we observe that $L$-$L$ pools vanish once we decrease the absenteeism rate in $L$-units. This indicates that (OPD) tends to pool together units with higher variability. In addition, comparing the results of sets 2 and 3, we observe that $H$-$H$ pools remain dominant regardless of the standard deviation of nurse demand. This indicates that the absenteeism rate plays a more important role than nurse demand in deciding the pattern of optimal pool design. This experiment suggests that we should prioritize pooling together the units with higher variability, especially those with higher absenteeism rates.}

\section{Conclusions} \label{sec:conclusion}
We studied a two-stage (DRNS) model for nurse staffing under both exogenous demand uncertainty and endogenous absenteeism uncertainty. We derived a min-max reformulation for (DRNS) under general nurse pool structures, leading to a separation algorithm that provably finds a globally optimal solution within a finite number of iterations. Under practical pool structures including one pool, disjoint pools, and chained pools, we {\color{black}derived MILP} reformulations for (DRNS) and significantly improved the computational efficacy. Via numerical case studies, we found that modeling absenteeism improves the out-of-sample performance of staffing decisions, and such improvement is positively correlated with the value of operational flexibility. For nurse pool design, we found that sparse pool structures can already harvest most of the operational flexibility. More importantly, it is particularly effective to pool together the units with higher nurse absence rates.

\newpage

\begin{appendices}

\section{Proof of Lemma \ref{prop:sns-sub-tu}} \label{apx-prop:sns-sub-tu}
\proof We rewrite formulation \eqref{sns-sub-obj}--\eqref{sns-sub-con-integer} as
\begin{subequations}
\begin{align}
V(\tilde{w}, \tilde{y}, \tilde{d}) \ = \ \min_{z, x} \ & \ \sum_{j=1}^J c^{\mbox{\tiny x}}_j x_j  \nonumber \\
\mbox{s.t.} \ & \ x_j + \sum_{i \in [I]: j \in P_i} z_{ij} \geq \tilde{d}_j - \tilde{w}_j, \ \ \forall j \in [J], \label{apx-sns-sub-con-demand} \\
& \ -\sum_{j \in P_i} z_{ij} \geq -\tilde{y}_i, \ \ \forall i \in [I], \label{apx-sns-sub-con-pool} \\
& \ x_j \in \mathbb{Z}_+, \ \ \forall j \in [J], \ \ z_{ij} \in \mathbb{Z}_+, \ \ \forall i \in [I], \ \forall j \in P_i. \nonumber
\end{align}
\end{subequations}
We note that the constraint matrix of the above formulation is totally unimodular (TU), and so the conclusion follows. 
To see the TU property, we consider the following constraint matrix:
\begin{equation*}
\begin{bmatrix}
\left(x_j + \sum_{i \in [I]: j \in P_i} z_{ij}\right) \\[0.15cm]
\left(-\sum_{j \in P_i} z_{ij}\right)
\end{bmatrix}.
\end{equation*}
It follows that (a) each entry of this matrix is $-1$, $0$, or $1$, (b) this matrix has at most two nonzero entries in each column, and (c) the entries sum up to be zero for any column containing two nonzero entries. Hence, the constraint matrix is TU based on Proposition 2.6 in~\cite{nemhauser1999integer}. The conclusion follows because $\tilde{d}_j - \tilde{w}_j$ and $-\tilde{y}_i$ are integers for all $j \in [J]$ and for all $i \in [I]$, respectively. \qed

\section{Verifying Assumption \ref{assump:technical}} \label{apx-prop:non-empty}
We present necessary and sufficient conditions for Assumption \ref{assump:technical} in the following proposition.
\begin{proposition} \label{prop:non-empty}
For any given $w$ and $y$, $\mathcal{D}(w,y)$ is non-empty if and only if the following three conditions are satisfied:
\begin{enumerate}[1.]
\item $f_j(w_j) \in [0, w_j]$ for all $j \in [J]$;
\item $g_i(y_i) \in [0, y_i]$ for all $i \in [I]$;
\item For all $j \in [J]$, the optimal value of the following linear program is non-positive:
\begin{subequations}
\label{apx-assump:technical-note-1}
\begin{align}
\min_{p_j \geq 0, \tau \geq 0} \ & \ \sum_{q=1}^Q (\tau^+_q + \tau^-_q) \\
\mbox{s.t.} \ & \ \sum_{k=d^{\mbox{\tiny L}}_j}^{d^{\mbox{\tiny U}}_j} k^q p_{jk} + \tau^+_q - \tau^-_q = \mu_{jq}, \ \ \forall q \in [Q], \\
& \ \sum_{k=d^{\mbox{\tiny L}}_j}^{d^{\mbox{\tiny U}}_j} p_{jk} = 1.
\end{align}
\end{subequations}
\end{enumerate}
\end{proposition}
\proof ({\bf Necessity}) Suppose that $\mathcal{D}(w,y) \neq \emptyset$. Then, there exists a $\mathbb{P} \in \mathcal{P}(\Xi)$ such that $\mathbb{E}_{\mathbb{P}}[\tilde{d}_j^q] = \mu_{jq}$ for all $j \in [J]$ and $q \in [Q]$, $\mathbb{E}_{\mathbb{P}}[\tilde{w}_j] = f_j(w_j)$ for all $j \in [J]$, and $\mathbb{E}_{\mathbb{P}}[\tilde{y}_i] = g_i(y_i)$ for all $i \in [I]$. It follows that, for all $j \in [J]$, we have $f_j(w_j) \leq \mbox{esssup}_{\Xi}\{\tilde{w}_j\} \leq w_j$ and $f_j(w_j) \geq \mbox{essinf}_{\Xi}\{\tilde{w}_j\} \geq 0$, leading to $f_j(w_j) \in [0, w_j]$. Likewise, it holds that $g_i(y_i) \in [0, y_i]$ for all $i \in [I]$. In addition, for all $j \in [J]$ and $k \in [d^{\mbox{\tiny L}}_j, d^{\mbox{\tiny U}}_j]_{\mathbb{Z}}$, we let $\bar{p}_{jk} = \mathbb{P}\{\tilde{d}_j = k\}$. It follows that, for all $q \in [Q]$, $\sum_{k=d^{\mbox{\tiny L}}_j}^{d^{\mbox{\tiny U}}_j} \bar{p}_{jk} = \sum_{k=d^{\mbox{\tiny L}}_j}^{d^{\mbox{\tiny U}}_j} \mathbb{P}\{\tilde{d}_j = k\} = 1$. Moreover,
\begin{equation*}
\sum_{k=d^{\mbox{\tiny L}}_j}^{d^{\mbox{\tiny U}}_j} k^q \bar{p}_{jk} \ = \ \sum_{k=d^{\mbox{\tiny L}}_j}^{d^{\mbox{\tiny U}}_j} k^q \mathbb{P}\{\tilde{d}_j = k\} \ = \ \mathbb{E}_{\mathbb{P}}[\tilde{d}_j^q] \ = \ \mu_{jq}.
\end{equation*}
Hence, together with $\tau^+_q = \tau^-_q = 0$, $\bar{p}_{jk}$ constitutes a feasible solution to linear program \eqref{apx-assump:technical-note-1} with an objective value being zero. As zero is also a lower bound of the objective value, $\bar{p}_{jk}$ is optimal to \eqref{apx-assump:technical-note-1} and accordingly the optimal value of this linear program equals zero. This holds for all $j \in [J]$ and proves the necessity of the three conditions.

\noindent ({\bf Sufficiency}) Suppose that the three conditions are satisfied. For all $j \in [J]$, as $f_j(w_j) \in [0, w_j] \equiv \mbox{conv}([0, w_j]_{\mathbb{Z}})$ by condition 1, there exists a $\mathbb{P}_{\tilde{w}_j} \in \mathcal{P}([0, w_j]_{\mathbb{Z}})$ such that $f_j(w_j) = \mathbb{E}_{\mathbb{P}_{\tilde{w}_j}}[\tilde{w}_j]$. Likewise, for all $i \in [I]$, there exists a $\mathbb{P}_{\tilde{y}_i} \in \mathcal{P}([0, y_i]_{\mathbb{Z}})$ such that $g_i(y_i) = \mathbb{E}_{\mathbb{P}_{\tilde{y}_i}}[\tilde{y}_i]$. In addition, as the optimal value of \eqref{apx-assump:technical-note-1} is non-positive and $\tau^+_q, \tau^-_q \geq 0$ for all $q \in [Q]$, the optimal value of \eqref{apx-assump:technical-note-1} equals zero. It follows that, for all $j \in [J]$, there exist $p_{jk}$ such that $\sum_{k=d^{\mbox{\tiny L}}_j}^{d^{\mbox{\tiny U}}_j} k^q p_{jk} = \mu_{jq}$ for all $q \in [Q]$ and $\sum_{k=d^{\mbox{\tiny L}}_j}^{d^{\mbox{\tiny U}}_j} p_{jk} = 1$. Defining $\mathbb{P}_{\tilde{d}_j} \in \mathcal{P}([d^{\mbox{\tiny L}}_j, d^{\mbox{\tiny U}}_j]_{\mathbb{Z}})$ such that $\mathbb{P}_{\tilde{d}_j}\{\tilde{d}_j = k\} = p_{jk}$ for all $k \in [d^{\mbox{\tiny L}}_j, d^{\mbox{\tiny U}}_j]_{\mathbb{Z}}$, we have $\mathbb{E}_{\mathbb{P}_{\tilde{d}_j}}[\tilde{d}^q_j] = \mu_{jq}$. Therefore, the probability distribution
$$
\mathbb{P} \ := \ \Pi_{j =1}^J \mathbb{P}_{\tilde{w}_j} \times \Pi_{i=1}^I \mathbb{P}_{\tilde{y}_i} \times \Pi_{j =1}^J \mathbb{P}_{\tilde{d}_j}
$$
satisfies constraints \eqref{demand-moments}--\eqref{nurse-linear-pool} and hence $\mathbb{P} \in \mathcal{D}(w,y)$. It follows that $\mathcal{D}(w,y) \neq \emptyset$ and the proof is completed. \qed

\section{Proof of Proposition \ref{prop:ref}} \label{apx-prop:ref}
\proof First, denoting $\tilde{\xi} := (\tilde{w}, \tilde{y}, \tilde{d})$, we present $\sup_{\mathbb{P} \in \mathcal{D}(w,y)}\mathbb{E}_{\mathbb{P}}[ V(\tilde{w}, \tilde{y}, \tilde{d})]$ as the following optimization problem:
\begin{subequations}
\begin{align}
\max_{p \geq 0} \ & \ \sum_{\tilde{\xi} \in \Xi} p_{\tilde{\xi}} V(\tilde{\xi}) \nonumber \\
\mbox{s.t.} \ & \ \sum_{\tilde{\xi} \in \Xi} p_{\tilde{\xi}} \tilde{w}_j = f_j(w_j), \ \ \forall j \in [J], \label{ref-note-1} \\
& \ \sum_{\tilde{\xi} \in \Xi} p_{\tilde{\xi}} \tilde{y}_i = g_i(y_i), \ \ \forall i \in [I], \label{ref-note-3} \\
& \ \sum_{\tilde{\xi} \in \Xi} p_{\tilde{\xi}} \tilde{d}^q_j = \mu_{jq}, \ \ \forall j \in [J], \ \forall q \in [Q], \label{ref-note-5} \\
& \ \sum_{\tilde{\xi} \in \Xi} p_{\tilde{\xi}} = 1, \label{ref-note-6}
\end{align}
\end{subequations}
where decision variables $p_{\tilde{\xi}}$ represent the probability of the random variables being realized as $\tilde{\xi}$, and constraints \eqref{ref-note-1}--\eqref{ref-note-6} describe the ambiguity set $\mathcal{D}(w,y)$ defined in \eqref{demand-moments}--\eqref{nurse-linear-pool}. The dual of this formulation is
\begin{subequations}
\begin{align}
\min_{\gamma, \lambda, \rho, \theta} \ & \ \sum_{j=1}^J \sum_{q=1}^Q \mu_{jq} \rho_{jq} + \sum_{j=1}^J f_j(w_j) \gamma_j + \sum_{i=1}^I g_i(y_i) \lambda_i + \theta \label{ref-note-obj} \\
\mbox{s.t.} \ & \ \theta + \sum_{j=1}^J \sum_{q=1}^Q \rho_{jq} \tilde{d}^q_j + \sum_{j=1}^J \gamma_j \tilde{w}_j + \sum_{i=1}^I \lambda_i \tilde{y}_i \geq V(\tilde{\xi}), \ \ \forall \tilde{\xi} \in \Xi, \label{ref-note-7}
\end{align}
\end{subequations}
where dual variables $\gamma_j$, $\lambda_i$, $\rho_{jq}$, and $\theta$ are associated with primal constraints \eqref{ref-note-1}--\eqref{ref-note-6}, respectively, and dual constraints \eqref{ref-note-7} are associated with primal variables $p_{\tilde{\xi}}$. By Assumption \ref{assump:technical}, strong duality holds between the primal and dual formulations because they are both linear programs. As the objective function aims to minimize the value of $\theta$, we observe by constraints \eqref{ref-note-7} that $\theta = \sup_{\tilde{\xi} \in \Omega}\{V(\tilde{\xi}) - \sum_{j=1}^J \sum_{q=1}^Q \rho_{jq} \tilde{d}^q_j - \sum_{j=1}^J \gamma_j \tilde{w}_j - \sum_{i=1}^I \lambda_i \tilde{y}_i\}$. Hence, $\sup_{\mathbb{P} \in \mathcal{D}(w,y)}\mathbb{E}_{\mathbb{P}}[ V(\tilde{\xi})]$ equals the optimal value of the following min-max optimization problem:
\begin{subequations}
\begin{align}
& \min_{\gamma, \lambda, \rho} \max_{\tilde{\xi} \in \Xi} \Biggl\{ V(\tilde{\xi}) - \sum_{j=1}^J\Biggl[ \sum_{q=1}^Q \rho_{jq} \tilde{d}^q_j + \gamma_j \tilde{w}_j \Biggr] - \sum_{i=1}^I \lambda_i \tilde{y}_i \Biggr\} \nonumber \\
& + \sum_{j=1}^J \Biggl[ \sum_{q=1}^Q \mu_{jq} \rho_{jq} + f_j(w_j) \gamma_j \Biggr] + \sum_{i=1}^I g_i(y_i) \lambda_i. \label{ref-note-42}
\end{align}
\end{subequations}
Second, in view of the dual formulation \eqref{sns-sub-dual-obj}--\eqref{sns-sub-dual-con-bd} of $V(\tilde{\xi})$, we rewrite the maximum term in \eqref{ref-note-42} as
\begin{subequations}
\begin{align}
& \max_{(\tilde{w}, \tilde{y}, \tilde{d}) \in \Xi} \max_{(\alpha, \beta) \in \Lambda} \left\{ \sum_{j=1}^J (\tilde{d}_j - \tilde{w}_j) \alpha_j + \sum_{i=1}^I \tilde{y}_i \beta_i - \sum_{j=1}^J\Biggl[ \sum_{q=1}^Q \rho_{jq} \tilde{d}^q_j + \gamma_j \tilde{w}_j \Biggr] - \sum_{i=1}^I \lambda_i \tilde{y}_i \right\} \nonumber \\
= \ & \max_{(\alpha, \beta) \in \Lambda} \max_{(\tilde{w}, \tilde{y}, \tilde{d}) \in \Xi} \left\{ \sum_{j=1}^J (\tilde{d}_j - \tilde{w}_j) \alpha_j + \sum_{i=1}^I \tilde{y}_i \beta_i - \sum_{j=1}^J\Biggl[ \sum_{q=1}^Q \rho_{jq} \tilde{d}^q_j + \gamma_j \tilde{w}_j \Biggr] - \sum_{i=1}^I \lambda_i \tilde{y}_i \right\} \nonumber \\
= \ & \max_{(\alpha, \beta) \in \Lambda} \Biggl\{ \sum_{j=1}^J \max_{\tilde{w}_j \in [0, w_j]_{\mathbb{Z}}}\Bigl\{(- \alpha_j - \gamma_j) \tilde{w}_j\Bigr\} + \sum_{i=1}^I \max_{\tilde{y}_i \in [0, y_i]_{\mathbb{Z}}} \Bigl\{ (\beta_i - \lambda_i) \tilde{y}_i \Bigr\} + \sum_{j=1}^J \max_{\tilde{d}_j \in [d^{\tinyl}_j, d^{\tinyu}_j]_{\mathbb{Z}}} \Bigl\{\alpha_j \tilde{d}_j - \sum_{q=1}^Q\rho_{jq}\tilde{d}^q_j \Bigr\} \Biggr\}. \nonumber
\end{align}
\end{subequations}
Finally, as $(- \alpha_j - \gamma_j) \tilde{w}_j$ is linear in $\tilde{w}_j$, we have
\begin{equation*}
\max_{\tilde{w}_j \in [0, w_j]_{\mathbb{Z}}}\Bigl\{(- \alpha_j - \gamma_j) \tilde{w}_j\Bigr\} = \max\Bigl\{0, (- \alpha_j - \gamma_j) w_j\Bigr\} = \Bigl[ (- \alpha_j - \gamma_j) w_j \Bigr]_+.
\end{equation*}
Similarly, we have $\max_{\tilde{y}_i \in [0, y_i]_{\mathbb{Z}}} \{ (\beta_i - \lambda_i) \tilde{y}_i \} = [(\beta_i - \lambda_i) y_i]_+$. This completes the proof. \qed

\section{Proof of Lemma \ref{thm:extreme-points}} \label{apx-thm:extreme-points}
\proof As $\max_{(\alpha, \beta) \in \Lambda} F(\alpha, \beta)$ is to maximize a convex function over a polyhedron, we only need to analyze the extreme directions and extreme points of $\Lambda$.

First, the extreme directions of $\Lambda$ are $(\alpha, \beta) = (0, - e_i)$ for all $i \in [I]$, where $e_i$ represents the $i^{\mbox{\tiny th}}$ standard basis vector. As $\tilde{y}_i \geq 0$, moving along any of these extreme directions (i.e., decreasing the value of any $\beta_i$) does not increase the value of $F(\alpha, \beta)$. Hence, we can omit these extreme directions in the attempt of maximizing $F(\alpha, \beta)$ and accordingly $\bar{\beta}_i = \min\{ -\bar{\alpha}_j: j \in P_i\} = - \max\{\bar{\alpha}_j: j \in P_i\} $ without loss of optimality. This proves property (b) in the claim. In addition, there exists an extreme point of $\Lambda$ that is optimal to $\max_{(\alpha, \beta) \in \Lambda} F(\alpha, \beta)$.

Second, we prove, by contradiction, that each extreme point of $\Lambda$ satisfies property (a) in the claim. Suppose that there exists an extreme point $(\bar{\alpha}, \bar{\beta})$ such that property (a) fails, i.e., $\bar{\alpha}_{j^*} \notin \{ 0, c^{\mbox{\tiny x}}_1, \ldots, c^{\mbox{\tiny x}}_{j^*} \}$ for some $j^* \in [J]$. Consider the set $\mathcal{I}(j^*) := \{i \in [I]: - \bar{\beta}_i = \bar{\alpha}_{j^*}\}$. We discuss the following two cases. In each case, we shall construct two points in $\Lambda$ such that their midpoint is $(\bar{\alpha}, \bar{\beta})$, which provides a desired contradiction.
\begin{enumerate}[1.]
\item If $\mathcal{I}(j^*) = \emptyset$, then $- \bar{\beta}_i > \bar{\alpha}_{j^*}$ for all $i$ such that $j^* \in P_i$. 
Defining $\epsilon := (1/2) \min \Big\{   - \bar{\beta}_i - \bar{\alpha}_{j^*}, \forall i \in [I]: j^* \in P_i, \ \bar{\alpha}_{j^*}, \  \min \big\{ |\bar{\alpha}_{j^*} - c^{\mbox{\tiny x}}_{\ell}| : \ell \in [j^*] \big\} \Big\} > 0$, we construct two points $(\bar{\alpha}^{+}, \bar{\beta})$ and $(\bar{\alpha}^{-}, \bar{\beta})$ such that $\bar{\alpha}^{+}_{j^*} = \bar{\alpha}_{j^*} + \epsilon$, $\bar{\alpha}^{-}_{j^*} = \bar{\alpha}_{j^*} - \epsilon$, and $\bar{\alpha}^{+}_j = \bar{\alpha}^{-}_j = \bar{\alpha}_j$ for all $j \neq j^*$. Then, it is clear that $(\bar{\alpha}^{+}, \bar{\beta}), (\bar{\alpha}^{-}, \bar{\beta}) \in \Lambda$. But $(\bar{\alpha}, \bar{\beta}) = (1/2)(\bar{\alpha}^{+}, \bar{\beta}) + (1/2)(\bar{\alpha}^{-}, \bar{\beta})$, which contradicts the fact that $(\bar{\alpha}, \bar{\beta})$ is an extreme point of $\Lambda$.
\item If $\mathcal{I}(j^*) \neq \emptyset$, then we define $\mathcal{J}(j^*) := \bigcup_{i \in \mathcal{I}(j^*)} \{j \in P_i: \bar{\alpha}_j = - \bar{\beta}_i\}$. It follows that $\bar{\alpha}_j = \bar{\alpha}_{j^*}$ for all $j \in \mathcal{J}(j^*)$. Hence, for each $i \in \mathcal{I}(j^*)$, $\bar{\alpha}_j = \bar{\alpha}_{j^*}$ for all $j \in P_i\cap\mathcal{J}(j^*)$ and $\bar{\alpha}_j < \bar{\alpha}_{j^*}$ for all $j \in P_i \setminus \mathcal{J}(j^*)$. 
We define $\epsilon := (1/2)\min\Bigl\{\min \big\{\bar{\alpha}_{j^*} - \bar{\alpha}_j: i \in \mathcal{I}(j^*), j \in P_i\setminus \mathcal{J}(j^*) \big\}, \ \min \big\{- \bar{\beta}_i - \bar{\alpha}_{j^*}: i \notin \mathcal{I}(j^*), \ - \bar{\beta}_i > \bar{\alpha}_{j^*} \big\}, \
\bar{\alpha}_{j^*}, \  \min \big\{ |\bar{\alpha}_{j^*} - c^{\mbox{\tiny x}}_{\ell}| : \ell \in [j^*] \big\}  \Bigr\}$. Then $\epsilon > 0$ because it is the minimum of a finite number of positive reals.\footnote{Here we adopt the convention that $\min\{a: a \in A\} = \infty$ if $A = \emptyset$. For example, if there does not exist an $i \notin \mathcal{I}(j^*)$ such that $- \bar{\beta}_i > \bar{\alpha}_{j^*}$, then $\min\{- \bar{\beta}_i - \bar{\alpha}_{j^*}: i \notin \mathcal{I}(j^*), \ - \bar{\beta}_i > \bar{\alpha}_{j^*}\} = \infty$.} We construct two points $(\bar{\alpha}^{+}, \bar{\beta}^{+})$ and $(\bar{\alpha}^{-}, \bar{\beta}^{-})$ such that
\begin{align*}
& \bar{\alpha}^{+}_j = \begin{cases}
\bar{\alpha}_{j^*} + \epsilon & \forall j \in \mathcal{J}(j^*) \\
\bar{\alpha}_j & \mbox{otherwise}
\end{cases}, \ \
\bar{\alpha}^{-}_j = \begin{cases}
\bar{\alpha}_{j^*} - \epsilon & \forall j \in \mathcal{J}(j^*) \\
\bar{\alpha}_j & \mbox{otherwise}
\end{cases}, \\[0.2cm]
& \bar{\beta}^{+}_i = \begin{cases}
- (\bar{\alpha}_{j^*} + \epsilon) & \forall i \in \mathcal{I}(j^*) \\
\bar{\beta}_i & \mbox{otherwise}
\end{cases}, \ \
\bar{\beta}^{-}_i = \begin{cases}
- (\bar{\alpha}_{j^*} - \epsilon) & \forall i \in \mathcal{I}(j^*) \\
\bar{\beta}_i & \mbox{otherwise}
\end{cases}.
\end{align*}
It is clear that $(\bar{\alpha}, \bar{\beta}) = (1/2)(\bar{\alpha}^{+}, \bar{\beta}^{+}) + (1/2)(\bar{\alpha}^{-}, \bar{\beta}^{-})$. To finish the proof, it remains to show that $(\bar{\alpha}^{+}, \bar{\beta}^{+}), (\bar{\alpha}^{-}, \bar{\beta}^{-}) \in \Lambda$. To see this, we check constraints \eqref{sns-sub-dual-con-max} and \eqref{sns-sub-dual-con-bd}. For constraints \eqref{sns-sub-dual-con-bd}, we have $\bar{\alpha}^{+}_j, \bar{\alpha}^{-}_j \in (0, c^{\mbox{\tiny x}}_j)$ for all $j \in \mathcal{J}(j^*)$ by the definition of $\epsilon$. Additionally, for all $j \notin \mathcal{J}(j^*)$, we have $\bar{\alpha}^{+}_j = \bar{\alpha}^{-}_j = \bar{\alpha}_j \in [0, c^{\mbox{\tiny x}}_j]$. Hence, constraints \eqref{sns-sub-dual-con-bd} are indeed satisfied and it remains to check constraints \eqref{sns-sub-dual-con-max}. For each $i \in \mathcal{I}(j^*)$, $- \bar{\beta}^{+}_i = \bar{\alpha}_{j^*} + \epsilon = \bar{\alpha}^{+}_j$ for all $j \in P_i\cap\mathcal{J}(j^*)$ and $- \bar{\beta}^{+}_i = \bar{\alpha}_{j^*} + \epsilon \geq \bar{\alpha}_{j^*} \geq \bar{\alpha}_j = \bar{\alpha}^{+}_j$ for all $j \in P_i\setminus\mathcal{J}(j^*)$, where the first inequality is because $\epsilon > 0$, and the second inequality follows from the definition of $\mathcal{J}(j^*)$. Meanwhile, $- \bar{\beta}^{-}_i = \bar{\alpha}_{j^*} - \epsilon = \bar{\alpha}^{-}_j$ for all $j \in P_i\cap\mathcal{J}(j^*)$, and $- \bar{\beta}^{-}_i = \bar{\alpha}_{j^*} - \epsilon \geq \bar{\alpha}_j = \bar{\alpha}^{-}_j$ for all $j \in P_i\setminus\mathcal{J}(j^*)$, where the inequality follows from the definition of $\epsilon$ and the last equality is because $j \notin \mathcal{J}(j^*)$. It follows that constraints \eqref{sns-sub-dual-con-max} are indeed satisfied for all $i \in \mathcal{I}(j^*)$. For each $i \notin \mathcal{I}(j^*)$, $\bar{\beta}^{+}_i = \bar{\beta}^{-}_i = \bar{\beta}_i$ and $-\bar{\beta}_i \neq \bar{\alpha}_{j^*}$. We discuss the following two sub-cases to complete the proof.
\begin{enumerate}[(a)]
\item If $- \bar{\beta}_i > \bar{\alpha}_{j^*}$, then $-\bar{\beta}^{+}_i = -\bar{\beta}_i \geq \bar{\alpha}_{j^*} + \epsilon \geq \bar{\alpha}^{+}_j$, where the first inequality follows from the definition of $\epsilon$. In addition, by construction $-\bar{\beta}^{-}_i = -\bar{\beta}_i > \bar{\alpha}_{j^*} \geq \bar{\alpha}^{-}_j$ for all $j \in P_i$.
\item If $- \bar{\beta}_i < \bar{\alpha}_{j^*}$, then $j \notin \mathcal{J}(j^*)$ for all $j \in P_i$ because otherwise $- \bar{\beta}_i \geq \bar{\alpha}_j = \bar{\alpha}_{j^*}$. It follows that $\bar{\alpha}^{+}_j = \bar{\alpha}^{-}_j = \bar{\alpha}_j$ and so $- \bar{\beta}^{+}_i = - \bar{\beta}_i \geq \bar{\alpha}_j = \bar{\alpha}^{+}_j$ and $- \bar{\beta}^{-}_i = - \bar{\beta}_i \geq \bar{\alpha}_j = \bar{\alpha}^{-}_j$. \qed
\end{enumerate}
\end{enumerate}

\section{Proof of Theorem \ref{thm:ip}} \label{apx-thm:ip}
\proof  First, pick any $(\alpha, \beta) \in \Lambda$ that satisfies the optimality conditions (a)--(b) stated in Lemma \ref{thm:extreme-points}. We shall show that there exists a feasible solution $(t, s, r, p)$ to formulation \eqref{IntegerProgram} that attains the same objective function value as $F(\alpha, \beta)$. 

To this end, for all $j \in [J]$ and $k \in [j]$, we let $t_{jk} = 1$ if $\alpha_j= c^{\text{\tiny x}}_k$, and $t_{jk}=0$ otherwise.
In addition, for all $i \in [I]$, if $\alpha_j = 0$ for all $j \in P_i$, then we let $s_{ik}=0$ for all $k \in [J(i)]$; and otherwise, i.e. there exists a $k^* \in [J(i)]$ such that $\max \{ \alpha_j : j \in P_i\} = c^{\text{\tiny x}}_{k^*}$, we let $s_{i k^*} = 1$ and all other $s_{ik} = 0$. Also, we define $r$ and $p$ as in \eqref{IntegerProgram-1} and \eqref{IntegerProgram-3}, respectively. 
By construction, $(t,s,r,p)$ satisfies \eqref{IntegerProgram-1}--\eqref{IntegerProgram-3}.
It follows that the objective function value of $(t,s,r,p)$ equals
\begin{align*}
& \sum_{j=1}^J \bigg( c ^{\text{\tiny r}}_j r_j + \sum_{k =1}^j c^{\text{\tiny t}}_{jk} t_{jk} \bigg) + \sum_{i=1}^I \bigg( c^{\text{\tiny p}}_i p_i +  \sum_{k =1}^{J(i)} c^{\text{\tiny s}}_{ik} s_{ik} \bigg)  \\
= \ 
& \sum_{j=1}^J \Bigg[
\mathbbm{1}_{\{ 0 \}}\bigl(  \alpha_j \bigr) \Bigg\{ [  - \gamma_j w_j ]_+ + \sup_{\tilde{d}_j \in [d^{\tinyl}_j, d^{\tinyu}_j]_{\mathbb{Z}}} \bigg\{  - \sum_{q=1}^Q \rho_{jq} \tilde{d}^q_j \bigg\} \Bigg\} + \\ 
& \sum_{k=1}^j \mathbbm{1}_{\{c^{\mbox{\tiny x}}_k\}}\bigl(  \alpha_j \bigr) \Bigg\{ [ (-c^{\text{\tiny x}}_k - \gamma_j ) w_j ]_+ + \sup_{\tilde{d}_j \in [d^{\tinyl}_j, d^{\tinyu}_j]_{\mathbb{Z}}} \bigg\{ c^{\text{\tiny x}}_k \tilde{d}_j - \sum_{q=1}^Q \rho_{jq} \tilde{d}^q_j \bigg\} \Bigg\} \Bigg] +
\\
& \sum_{i=1}^I \Bigg[ \Bigg\{ 
\mathbbm{1}_{\{0\}} \Bigl(  \max_{j \in P_i} \{\alpha_j \}  \Bigr) [ - \lambda_i  y_i  ]_+ \Bigg\} + 
\sum_{k=1}^{J(i)} \Bigg\{ \mathbbm{1}_{\{ c^{\mbox{\tiny x}}_k\}} \Bigl( \max_{j \in P_i} \{\alpha_j \}  \Bigr) [ (-c^{\text{\tiny x}}_k - \lambda_i ) y_i  ]_+ \Bigg\} \Bigg] \\
= \
& \sum_{j=1}^J \Biggl\{ \Bigl[(- \alpha_j - \gamma_j) w_j\Bigr]_+ + \sup_{\tilde{d}_j \in [d^{\tinyl}_j, d^{\tinyu}_j]_{\mathbb{Z}}} \Bigl\{\alpha_j \tilde{d}_j - \sum_{q=1}^Q\rho_{jq}\tilde{d}^q_j \Bigr\} \Biggr\} + \sum_{i=1}^I \Bigl[ (\beta_i - \lambda_i) y_i \Bigr]_+ \\
= \
& F(\alpha, \beta),
\end{align*}
where the first equality follows from the definition of $(t, s, r, p)$ and the second equality follows from the optimality conditions stated in Lemma \ref{thm:extreme-points}.

Second, pick any feasible solution $(t, s, r, p)$ to formulation \eqref{IntegerProgram}. 
We construct an $(\alpha, \beta) \in \Lambda$ that satisfies the optimality conditions (a)--(b) in Lemma \ref{thm:extreme-points} and $F(\alpha, \beta)$ equals the objective function value \eqref{IntegerProgram-0} of $(t, s, r, p)$. 
Specifically, for all $j \in [J]$, we let $\alpha_j =  \sum_{k=1}^j c^{\mbox{\tiny x}}_k t_{jk}$ and, for all $i \in [I]$, $\beta_i = - \sum_{k=1}^{J(i)} c^{\text{\tiny x}}_k s_{ik}$. 
Then, for all $i \in [I]$ and $j \in P_i$,
\begin{align*}
\beta_i + \alpha_j & = - \sum_{k=1}^{J(i)} c^{\text{\tiny x}}_k s_{ik} + \sum_{k=1}^j c^{\text{\tiny x}}_k t_{jk} \leq  - \sum_{k=1}^{J(i)} c^{\text{\tiny x}}_k (\sum_{\ell \in P_i: \ell \geq k} t_{\ell k}) + \sum_{k=1}^j c^{\text{\tiny x}}_k t_{jk} \\
& =  - \sum_{\ell \in P_i} \sum_{k=1}^{\ell} c^{\text{\tiny x}}_k  t_{\ell k} + \sum_{k=1}^j c^{\text{\tiny x}}_k t_{jk} =  - \sum_{\ell \in P_i \setminus \{j\}} \sum_{k=1}^{\ell} c^{\text{\tiny x}}_k  t_{\ell k} \leq 0,
\end{align*}
where the first inequality is due to constraints \eqref{encode-2}, and the last inequality is due to $c^{\text{\tiny x}}_{k}, t_{\ell k} \geq 0$.
Next, we have $\alpha_j \in \{0, c^{\textmd{\mbox{\tiny x}}}_1, \ldots, c^{\textmd{\mbox{\tiny x}}}_j \}$ for all $j \in [J]$ due to constraints \eqref{IntegerProgram-1}, and $\beta_i \in \{0, -c^{\text{\tiny x}}_1, \ldots, -c^{\text{\tiny x}}_{J(i)} \}$ for all $i \in [I]$ due to constraints \eqref{IntegerProgram-3}.
Hence, $(\alpha, \beta) \in \Lambda$ and satisfies optimality condition (a).
It remains to show that the constructed $(\alpha,\beta)$ satisfies optimality condition (b).

For all $i \in [I]$, if $\sum_{k=1}^{J(i)} s_{ik} = 0$, i.e., $\beta_i = 0$, then $t_{jk} = 0$ for all $j \in P_i$ and $k \in [j]$ due to constraints \eqref{encode-4}. It follows that $\alpha_j = 0$ for all $j \in P_i$ and so $\beta_i = - \max\{\alpha_j: j \in P_i\} = 0$. 
On the other hand, if $\sum_{k=1}^{J(i)} s_{ik} = 1$, there exists a $k^* \in [J(i)]$ with $s_{ik^*} = 1$, i.e., $\beta_i = -c^{\text{\tiny x}}_{k^*}$.
By constraints \eqref{encode-2} and \eqref{encode-3}, $\sum_{j \in P_i: j \geq k^*} t_{jk^*} \geq 1$ and  $t_{j \ell} = 0$ for all $j \in P_i: j \geq k^*+1$ and $\ell \in [k^*+1, j]_{\mathbb{Z}\textsl{}}$. It follows that $\max\{\alpha_j: j \in P_i\} = c^{\textmd{\mbox{\tiny x}}}_{k^*}$ and so $\beta_i = - \max\{\alpha_j: j \in P_i\} = - c^{\text{\tiny x}}_{k^*}$. Hence, $(\alpha, \beta)$ satisfies the optimality conditions (b). Finally,
\begin{align*}
F(\alpha, \beta) = \ & \ \sum_{j=1}^J \Biggl\{ \Bigl[( - \sum_{k=1}^j c^{\textmd{\mbox{\tiny x}}}_k t_{jk} - \gamma_j) w_j\Bigr]_+ + \sup_{\tilde{d}_j \in [d^{\tinyl}_j, d^{\tinyu}_j]_{\mathbb{Z}}} \Bigl\{\bigl( \sum_{k=1}^j c^{\textmd{\mbox{\tiny x}}}_k t_{jk} \bigr) \tilde{d}_j - \sum_{q=1}^Q\rho_{jq}\tilde{d}^q_j \Bigr\} \Biggr\} \\
& \ + \sum_{i=1}^I \Bigl[ \Bigl( - \sum_{k=1}^{J(i)} c^{\textmd{\mbox{\tiny x}}}_k s_{ik} - \lambda_i\Bigr) y_i \Bigr]_+ \\
= \ & \sum_{j=1}^J \bigg( c ^{\text{\tiny r}}_j r_j + \sum_{k =1}^j c^{\text{\tiny t}}_{jk} t_{jk} \bigg) + \sum_{i=1}^I \bigg( c^{\text{\tiny p}}_i p_i +  \sum_{k =1}^{J(i)} c^{\text{\tiny s}}_{ik} s_{ik} \bigg) 
\end{align*}
by the definition of $(\alpha, \beta)$ and constraints \eqref{IntegerProgram-1}--\eqref{IntegerProgram-3}. This completes the proof. \qed

\section{Proof of Proposition \ref{prop:big-M}} \label{apx-prop:big-M}
\proof Let $G(\gamma, \lambda, \rho)$ be the objective function of problem \eqref{ref-note-8}--\eqref{ref-note-con}, $(\hat{\gamma}, \hat{\lambda}, \hat{\rho})$ be any feasible solution, and $S^*$ be the set of optimal solution to problem $\max_{(\alpha, \beta) \in \Lambda}F(\alpha, \beta)$ for the given $(\hat{\gamma}, \hat{\lambda}, \hat{\rho})$. 
Suppose that there exists a $j \in [J]$ such that $\hat{\gamma}_j < - c^{\mbox{\tiny x}}_j$. Then, $- \hat{\gamma}_j - \alpha^*_j > 0$ and $[(- \hat{\gamma}_j - \alpha^*_j)w_j]_+ = (- \hat{\gamma}_j - \alpha^*_j)w_j$ for all $(\alpha^*, \beta^*) \in S^*$ {\color{black}because $\alpha^*_j \leq c^{\mbox{\tiny x}}_j$ by \eqref{sns-sub-dual-con-bd}}.
Additionally, due to Lemma \ref{thm:extreme-points}, we can replace polyhedron $\Lambda$ by the (compact) set of its extreme points $\mbox{ex}(\Lambda)$ without loss of optimality, i.e., $\max_{(\alpha, \beta) \in \Lambda}F(\alpha, \beta) = \max_{(\alpha, \beta) \in \mbox{\tiny ex}(\Lambda)}F(\alpha, \beta)$. It then follows from Theorem 2.87 in~\cite{ruszczynski2006nonlinear} that, for all subgradient $\varpi \in \partial G(\hat{\gamma}, \hat{\lambda}, \hat{\rho})$, the entry in $\varpi$ with regard to variable $\gamma_j$ at $(\hat{\gamma}, \hat{\lambda}, \hat{\rho})$ equals $f_j(w_j) - w_j$, i.e., $\varpi(\gamma_j)|_{(\hat{\gamma}, \hat{\lambda}, \hat{\rho})} = f_j(w_j) - w_j \leq 0$. Noting that $\varpi(\gamma_j)|_{(\hat{\gamma}, \hat{\lambda}, \hat{\rho})} \leq 0$ holds valid whenever $\hat{\gamma}_j < - c^{\mbox{\tiny x}}_j$, we can increase $\hat{\gamma}_j$ to $- c^{\mbox{\tiny x}}_j$ without any loss of optimality.

Now suppose that $\hat{\gamma}_j > 0$. Then, we have $-\hat{\gamma}_j - \alpha^*_j < 0$ and $[(- \hat{\gamma}_j - \alpha^*_j)w_j]_+ = 0$ for all $(\alpha^*, \beta^*) \in S^*$ {\color{black}because $\alpha^*_j \geq 0$ by \eqref{sns-sub-dual-con-bd}}. It follows from a similar implication as in the previous paragraph that, for all subgradient $\varpi \in \partial G(\hat{\gamma}, \hat{\lambda}, \hat{\rho})$, we have $\varpi(\gamma_j)|_{(\hat{\gamma}, \hat{\lambda}, \hat{\rho})} = f_j(w_j) \geq 0$. Noting that this holds valid whenever $\hat{\gamma}_j > 0$, we can decrease $\hat{\gamma}_j$ to $0$ without any loss of optimality. Therefore, there exists an optimal solution $(\gamma^*, \lambda^*, \rho^*)$ to problem \eqref{ref-note-8}--\eqref{ref-note-con} such that $\gamma^*_j \in [-c^{\mbox{\tiny x}}_j, 0]$ for all $j \in [J]$.

Following a similar proof, we can show that there exists an optimal solution $(\gamma^*, \lambda^*, \rho^*)$ to problem \eqref{ref-note-8}--\eqref{ref-note-con} such that $\lambda^*_i \in [-c^{\mbox{\tiny x}}_{J(i)}, 0]$ for all $i \in [I]$. We omit the details for the sake of saving space. \qed

\section{Proof of Theorem \ref{thm:finite}} \label{apx-thm:finite}
\proof In each iteration of Algorithm \ref{algo-sep}, we solve a relaxation of the (DRNS) reformulation \eqref{ref-linear-1}--\eqref{ref-linear-4}. It follows that, if the algorithm stops in an iteration and returns a solution $(u^*, v^*)$ then $(u^*, v^*)$ satisfies all the constraints \eqref{ref-linear-3} because of Step 5. Then, $(u^*, v^*)$ is feasible to formulation \eqref{ref-linear-1}--\eqref{ref-linear-4} and meanwhile optimal to its relaxation. Hence, $(u^*, v^*)$ is optimal to formulation \eqref{ref-linear-1}--\eqref{ref-linear-4}, i.e., optimal to (DRNS).

It remains to show that Algorithm \ref{algo-sep} stops in a finite number of iterations. To see this, we notice that the set $\mathcal{H}$ contains a finite number of elements. Indeed, binary variables $t$ and $s$ only have a finite number of possible values. Although $r$ and $p$ are continuous variables, they also only have a finite number of possible values due to constraints \eqref{IntegerProgram-1}--\eqref{IntegerProgram-3}. \qed

\section{Proof of Lemma \ref{lem:valid}} \label{apx-lem:valid}
\proof
Pick any $i \in [I]$. We note that $\sum_{\ell = 1}^{J(i)} s_{i \ell} \in \{0,1\}$ due to constraints \eqref{encode-1} and discuss the following three cases.
First, if $\sum_{\ell = 1}^{J(i)} s_{i \ell} = 0$, i.e., $s_{i \ell} = 0$ for all $\ell \in [J(i)]$, then $t_{jk}=0$ for all $j \in P_i$ and $k \in [j]$ by constraints \eqref{encode-4}. 
In this case, inequalities \eqref{VI-general} are valid because both left-hand and right-hand sides of \eqref{VI-general} are $0$.
Second, if $\sum_{\ell = 1}^{J(i)} s_{i \ell} = 1$ and $\sum_{\ell = k}^{J(i)} s_{i \ell} = 1$, then inequalities \eqref{VI-general} reduce to $\sum_{\ell = k}^j t_{j \ell} \leq 1$ for all $j \in P_i$ with $j \geq k$, which are valid due to constraints \eqref{encode-0}.
Third, if $\sum_{\ell = 1}^{J(i)} s_{i \ell} = 1$ and $\sum_{\ell = k}^{J(i)} s_{i \ell} = 0$, then there exists a $k^* < k$ such that $s_{ik^*} = 1$. By constraints \eqref{encode-3}, $t_{j \ell} = 0$ for all $j \in P_i$ with $j \geq k^* + 1$ and for all $\ell \in [k^*+1,j]_{\mathbb{Z}}$. It follows that both left-hand and right-hand sides of \eqref{VI-general} are $0$, and so inequalities \eqref{VI-general} are valid.

Finally, inequalities \eqref{VI-general} imply constraints \eqref{encode-4} because
\begin{align*}
& t_{jk} \leq \sum_{\ell = k }^j t_{j \ell} \leq \sum_{\ell = k}^{J(i)} s_{i \ell} \leq \sum_{\ell = 1}^{J(i)} s_{i \ell}, \ \forall i \in [I], \forall j \in P_i, \forall k \in [j].
\end{align*}
\qed

\section{Proof of Proposition \ref{prop:onepool-convex}} \label{apx-prop:onepool-convex}

We first show that $\overline{\mathcal{H}}_{\text{\tiny 1}} \subseteq \text{conv} (\mathcal{H}_{\text{\tiny 1}})$.
To prove this, we pick any fractional solution $(\hat{s}, \hat{t})$ that satisfies constraints \eqref{H-onepool-1}--\eqref{H-onepool-3}, \eqref{onepool-VI}, and \eqref{onepool-nonneg} in $\overline{\mathcal{H}}_{\text{\tiny 1}}$, and find a finite number of points $\{(s^1,t^1), \ldots, (s^n, t^n)\}$ such that: (i) each $(s^i,t^i)$ is binary-valued and satisfies constraints \eqref{H-onepool-1}--\eqref{H-onepool-6} in $\mathcal{H}_{\text{\tiny 1}}$ and (ii) their convex combination produces $(\hat{s}, \hat{t})$, i.e., we find nonnegative weights $\{\pi^i: i \in [n]\}$ with $\sum_{i=1}^n \pi^i = 1$ such that $\hat{s}_k = \sum_{i=1}^n \pi^i s^i_k$ for all $k \in [J]$ and $\hat{t}_{jk} = \sum_{i=1}^n \pi^i t^i_{jk}$ for all $j \in [J]$ and $k \in [j]$.

To this end, we develop Algorithm \ref{algo:convex}, which iteratively constructs binary points $\{(s^i,t^i): i \in [n]\}$ and weights $\{\pi^i: i\in [n]\}$. This algorithm initializes by setting an incumbent $(\hat{s}^1, \hat{t}^1) := (\hat{s},\hat{t})$ and, in each iteration $i \in [n]$, update the incumbent $(\hat{s}^{i+1}, \hat{t}^{i+1}) := (\hat{s}^{i}- \pi^i s^{i}, \ \hat{t}^{i} - \pi^i t^{i})$ as it constructs a new binary point $(s^i, t^i)$ with weight $\pi^i$, till the incumbent becomes $(0, 0)$. For ease of exposition, we assume in Algorithm \ref{algo:convex} that optimizing over an empty set yields zero, i.e., $\max\{f(x): x \in \emptyset\} = 0$. In addition, we provide an example of implementing Algorithm \ref{algo:convex} at the end of this section. We show the correctness of Algorithm \ref{algo:convex} by proving its properties in Theorem \ref{thm:onepool-Algorithm}.

\begin{theorem} \label{thm:onepool-Algorithm}
	For each binary point $(s^i, t^i)$ produced in Algorithm \ref{algo:convex}, the following properties hold:
	\begin{enumerate}
		\item[1.] $(s^i, t^i)$ satisfies constraints \eqref{H-onepool-1}--\eqref{H-onepool-6} in $\mathcal{H}_{\text{\tiny 1}}$.
		\item[2.] In line 6, the weight $\pi^i \in (0,1]$.	
		\item[3.] In line 7, the updated fractional solution $(\hat{s}^{i+1}, \hat{t}^{i+1})$ satisfies constraints \eqref{H-onepool-1}--\eqref{H-onepool-3}, \eqref{onepool-VI}, and \eqref{onepool-nonneg} in $\overline{\mathcal{H}}_{\text{\tiny 1}}$.
	\end{enumerate}
	In a finite number of iterations, Algorithm \ref{algo:convex} terminates with the following properties:
	\begin{enumerate}
		\item[4.] After we execute line 7 for the last time, $\hat{s}^n_k = 0$ for all $k \in [J]$ and $\hat{t}^n_{jk} = 0$ for all $j \in [J]$ and $k \in [j]$.
		\item[5.] $\pi^n \in [0,1)$.
	\end{enumerate}
	Therefore, $(\hat{s}, \hat{t})$ is a convex combination of $\{(s^1, t^1), \ldots, (s^n, t^n)\}$ with weights $\pi^1, \ldots, \pi^n$, respectively.	
\end{theorem}
\proof See Appendix \ref{apx-thm:onepool-Algorithm}.  \qed

\begin{algorithm}[H]
	\caption{Finding binary points $\{(s^1, t^1), \ldots , (s^n, t^n) \}$ and their weights $\{\pi^1, \ldots, \pi^n \}$} 
	\label{algo:convex}
	\begin{algorithmic}[1]	
		\STATE \textbf{Initialization:} $i = 1$, and $\hat{s}^1_k = \hat{s}_k, \ \forall k \in [J], \ \hat{t}^1_{jk} = \hat{t}_{jk}, \ \forall j \in [J], \forall k \in [j]$.
		\FOR{$m = 1:J$} 
		\WHILE{$\hat{s}^i_{m} > 0$}
		\STATE \textbf{Find:} 
		\begin{small}
			\begin{align*}
			& j^* = \min \big\{  j \in [m,J]_{\mathbb{Z}} : \ \hat{t}^i_{j m} > 0 \big\}, \\
			& A_{m} = \big\{j^* \big\} \cup \left\{ j \in [m,J]_{\mathbb{Z}}: \ \sum_{\ell=m}^j \hat{t}^i_{j \ell} = \sum_{\ell=m}^J \hat{s}^i_{\ell} \right\}, \ \ B_{m} = [m,J]_{\mathbb{Z}} \setminus A_{m}, \\
			& \text{Index}(j) =
			\begin{cases}
			m & \forall j \in A_{m}, \\
			\max \big\{k \in [ \min \{j,m-1\} ]: \hat{t}^i_{jk} > 0 \big\} & \forall j \in B_{m} \cup [m-1].
			\end{cases}
			\end{align*}
		\end{small}			
		\STATE \textbf{Construct a binary point $(s^{i},t^{i})$:}		
		\begin{small}
			\begin{align*}
			s^i_k = 
			\begin{cases}
			1 & \text{ if } k = m \\
			0 & \text{ otherwise }
			\end{cases}
			,\ \forall k \in [J], \ \text{ and } \
			t^i_{jk} =
			\begin{cases}
			1 & \text{ if } k = \text{Index}(j) \neq 0 \\
			0 & \text{ otherwise }
			\end{cases}		
			,\ \forall j \in [J], \forall k \in [j].
			\end{align*}
		\end{small}
		\STATE \textbf{Construct a weight $\pi^i$:} 	
		\vspace{-3mm}
		\begin{small}
			\begin{align*}
			& \pi^i = \min \bigg\{ \sum_{\ell=m}^J \hat{s}^i_{\ell} - M_{m}, \ \ \hat{s}^i_{m}, \ \ \hat{t}^i_{j, \text{Index}(j)}, \ \forall j \in [J] :	\text{Index}(j) \neq 0 \bigg\}, \\ 
			& \text{ where } M_{m} = \max \bigg\{ \sum_{\ell=m}^j \hat{t}^i_{j \ell}:j \in B_m \bigg\}.
			\end{align*}
		\end{small}
		\vspace{-3mm}
		\STATE \textbf{Update the fractional solution $(\hat{s}^{i+1}, \hat{t}^{i+1})$:}
		\begin{small}
			\begin{align*}
			& \hat{s}^{i+1}_{k} = \hat{s}^i_{k} - \pi^i s^i_{k}, \ \forall k \in [J] \ \text{ and } \
			\hat{t}^{i+1}_{jk} = \hat{t}^i_{jk} - \pi^i t^i_{jk}, \ \forall j \in [J], \forall k \in [j].
			\end{align*}
		\end{small}
		\STATE \textbf{Update:} $i = i + 1$.		
		\ENDWHILE
		\ENDFOR	
		\STATE \textbf{Set:} $n = i$, $\pi^n = 1 -  (\pi^1+\cdots+\pi^{n-1})$, $s^n_k=0, \ \forall k \in [J]$, and $t^n_{jk} = 0, \ \forall j \in [J], \forall k \in [j]$.
	\end{algorithmic}
\end{algorithm}

\noindent It follows from Theorem \ref{thm:onepool-Algorithm} that $\overline{\mathcal{H}}_{\text{\tiny 1}} \subseteq \text{conv} (\mathcal{H}_{\text{\tiny 1}})$.
On the other hand, $\text{conv} (\mathcal{H}_{\text{\tiny 1}}) \subseteq \overline{\mathcal{H}}_{\text{\tiny 1}}$ because $\mathcal{H}_{\text{\tiny 1}} \subseteq \overline{\mathcal{H}}_{\text{\tiny 1}}$ by construction.
Therefore, $\overline{\mathcal{H}}_{\text{\tiny 1}} = \text{conv} (\mathcal{H}_{\text{\tiny 1}})$ and this completes the proof. \qed

\section{Proof of Theorem \ref{thm:onepool-Algorithm}} \label{apx-thm:onepool-Algorithm}
\proof
We first prove Property 3 and then use it to prove all other properties. \\
\noindent\textbf{(Property 3)}
We prove by mathematical induction.
It is clear that $(\hat{s}^1, \hat{t}^1) \equiv (\hat{s}, \hat{t})$ satisfies Property 3. It remains to show that $(\hat{s}^{i+1}, \hat{t}^{i+1})$ satisfies constraints \eqref{H-onepool-1}--\eqref{H-onepool-3}, \eqref{onepool-VI}, and \eqref{onepool-nonneg} under the assumption that they are satisfied by $(\hat{s}^{i}, \hat{t}^{i})$.
Before getting into details, we rewrite the construction of $(\hat{s}^{i+1}, \hat{t}^{i+1})$ based on line 5 as follows:
\begin{align} \label{updated_fractional} 
\hat{s}^{i+1}_k = 
\begin{cases}
\hat{s}^i_k 		& \text{if } k \in [m-1] \\
\hat{s}^i_m - \pi^i 	& \text{if } k = m \\
\hat{s}^i_k 		& \text{if } k \in [m+1,J]_{\mathbb{Z}}
\end{cases}
, \
\hat{t}^{i+1}_{jk} = 
\begin{cases}
\hat{t}^{i}_{jk} 		& \text{if } k \in [m-1],  j \in [k,J]_{\mathbb{Z}} : k < \text{Index(j)}  \\
\hat{t}^{i}_{jk} - \pi^i  & \text{if } k \in [m-1],  j \in [k,J]_{\mathbb{Z}} : k = \text{Index(j)} \\
0 						& \text{if } k \in [m-1],  j \in [k,J]_{\mathbb{Z}} : k > \text{Index(j)} \\
\hat{t}^{i}_{jm} - \pi^i  & \text{if } k  = m,  j \in A_m \\
\hat{t}^{i}_{jm} 		& \text{if } k  = m,  j \in B_m \\
\hat{t}^{i}_{jk} 		& \text{if } k \in [m+1,J]_{\mathbb{Z}},  j \in [k,J]_{\mathbb{Z}} \\
\end{cases}
.
\end{align}

First, we consider constraints \eqref{onepool-nonneg}. We notice that $\hat{s}^{i+1}_k \geq 0$ for all $k \in [J]$ because $\pi^i \leq \hat{s}^i_m$ by construction. Likewise, $\hat{t}^{i+1}_{jk} \geq 0$ for all $j \in [J]$ and $k \in [j]$. 

Second, we consider constraint \eqref{H-onepool-2}. By \eqref{updated_fractional}, we have $\sum_{k=1}^J \hat{s}^{i+1}_k = \sum_{k=1}^J \hat{s}^{i}_k - \pi^i \leq \sum_{k=1}^J \hat{s}^{i}_k \leq 1$ because $\pi^i \geq 0$ by construction.

Third, we consider constraints \eqref{H-onepool-1}. By \eqref{updated_fractional}, for all $j \in [J]$, we have $\sum_{k=1}^j \hat{t}^{i+1}_{jk} = \sum_{k=1}^j \hat{t}^{i}_{jk} - \mathbbm{1}(\text{Index}(j) \neq 0) \pi^i \leq \sum_{k=1}^j \hat{t}^{i}_{jk} \leq 1$ because $\pi^i \geq 0$ by construction.

Fourth, we consider constraints \eqref{onepool-VI}.
For fixed $k \in [J]$ and $j \in [k, J]_{\mathbb{Z}}$, we rewrite each term in \eqref{onepool-VI} as follows:
\begin{subequations}
	\begin{align} 
	& \sum_{\ell=k}^J \hat{s}^{i+1}_{\ell} = 
	\begin{cases}
	\sum_{\ell=k}^J \hat{s}^i_{\ell} - \pi^i  & \text{if } k \in [m] \\
	\sum_{\ell=k}^J \hat{s}^i_{\ell} 		& \text{if } k \in [m+1,J]_{\mathbb{Z}}
	\end{cases}
	, \\
	& \sum_{\ell=k}^j \hat{t}^{i+1}_{j \ell} =
	\begin{cases}
	\sum_{\ell=k}^j \hat{t}^i_{j\ell} - \pi^i   & \text{if } k \in [m-1], j \in [k, m-1]_{\mathbb{Z}}: k \leq \text{Index}(j) \\
	0                                         & \text{if } k \in [m-1], j \in [k, m-1]_{\mathbb{Z}}: k > \text{Index}(j) \\
	\sum_{\ell=k}^j \hat{t}^i_{j\ell} - \pi^i   & \text{if } j \in A_m, k \in [m] \\
	\sum_{\ell=k}^j \hat{t}^i_{j\ell}         & \text{if } j \in A_m, k \in [m+1, j]_{\mathbb{Z}} \\
	\sum_{\ell=k}^j \hat{t}^i_{j\ell} - \pi^i   & \text{if } j \in B_m, k \in [j]: k \leq \text{Index}(j) \\
	\sum_{\ell=k}^j \hat{t}^i_{j\ell}         & \text{if } j \in B_m, k \in [j]: k > \text{Index}(j) \\
	\end{cases}
	. \label{cases}
	\end{align}
\end{subequations}
We prove by enumerating all $6$ cases of $\sum_{\ell=k}^j \hat{t}^{i+1}_{j \ell}$ in \eqref{cases}. 
In cases 1 and 3, \eqref{onepool-VI} is satisfied because $\sum_{\ell=k}^j \hat{t}^i_{j\ell} - \pi^i \leq \sum_{\ell=k}^J \hat{s}^i_\ell - \pi^i$, which follows from the induction assumption that $(\hat{s}^i, \hat{t}^i)$ satisfies \eqref{onepool-VI}. 
Likewise, \eqref{onepool-VI} is satisfied in cases 4 and 5. 
In case 2, \eqref{onepool-VI} is satisfied because $\pi^i \leq \hat{s}^i_m \leq \sum_{\ell=k}^J \hat{s}^i_\ell$. 
In case 6, we discuss the following two sub-cases:
(i) If $k \geq m + 1$, then \eqref{onepool-VI} is satisfied because $\sum_{\ell=k}^j \hat{t}^i_{j\ell} \leq \sum_{\ell=k}^J \hat{s}^i_\ell$ by the induction assumption.
(ii) If $\text{Index}(j) < k \leq m$, then $\sum_{\ell=k}^j \hat{t}^i_{j\ell} = \sum_{\ell=m}^j \hat{t}^i_{j\ell}$ by definition of $\text{Index}(j)$. It follows that $\pi^i \leq \sum_{\ell=m}^J \hat{s}^i_\ell - M_m \leq \sum_{\ell=k}^J \hat{s}^i_\ell - \sum_{\ell=m}^j \hat{t}^i_{j\ell}$, where the first inequality is due to the definition of $\pi^i$ and the second inequality is due to that of $M_m$. This proves that $(\hat{s}^{i+1}, \hat{t}^{i+1})$ satisfies \eqref{onepool-VI}.

Finally, we consider constraints \eqref{H-onepool-3}. 
By \eqref{updated_fractional}, we rewrite the right-hand side of constraints \eqref{H-onepool-3} as follows:
\begin{align}
& \sum_{j=k}^J \hat{t}^{i+1}_{jk} =
\begin{cases}
\sum_{j=k}^J (\hat{t}^{i}_{jk} - \mathbbm{1}_{\text{Index}(j)}(k) \pi^i) & \text{if } k \in [m-1] \\
\sum_{j=m}^J \hat{t}^{i}_{jm} -  \pi^i 									& \text{if } k = m \text{ and } A_m = \{j^*\} \\
\sum_{j=m}^J \hat{t}^{i}_{jm} - \sum_{j \in A_m} \pi^i 					& \text{if } k = m \text{ and } A_m \supset \{j^*\} \\
\sum_{j=k}^J \hat{t}^{i}_{jk} 											& \text{if } k \in [m+1,J]_{\mathbb{Z}}
\end{cases}
, \label{cases-2}
\end{align}
where $\mathbbm{1}_{\text{Index}(j)}(k) = 1$ if $k=\text{Index}(j)$ and $\mathbbm{1}_{\text{Index}(j)}(k) = 0$ otherwise.
We prove by enumerating all $4$ cases of $\sum_{j=k}^J \hat{t}^{i+1}_{jk}$ in \eqref{cases-2}.
In case 1, we notice that (i) $\hat{s}^i_{k} = 0, \ \forall k \in [m-1]$ due to line 3, (ii) $\hat{t}^i_{jk} \geq 0$ for all $j \in [k, J]_{\mathbb{Z}}$, and (iii) $\pi^i \leq \hat{t}^i_{j, \text{Index}(j)}, \ \forall j \in [J]:\text{Index}(j) \neq 0$ by the definition of $\pi^i$. Therefore, \eqref{H-onepool-3} is satisfied:
\begin{align*}
& \hat{s}^{i+1}_k = \hat{s}^i_k = 0 \leq \sum_{j=k}^J (\hat{t}^{i}_{jk} - \mathbbm{1}_{\text{Index}(j)}(k) \pi^i) = \sum_{j=k}^J \hat{t}^{i+1}_{jk}.
\end{align*}
In case 2, \eqref{H-onepool-3} is satisfied because $\hat{s}^i_m - \pi^i \leq \sum_{j=m}^J \hat{t}^i_{jm} - \pi^i$, which follows from the induction assumption that $(\hat{s}^i, \hat{t}^i)$ satisfies \eqref{H-onepool-3}.
In case 3, we notice that $\sum_{\ell=m}^j \hat{t}^i_{j \ell} = \sum_{\ell = m}^J \hat{s}^i_{\ell}$  for all $j \in A_m \setminus \{j^*\} \neq \emptyset$.
As $\sum_{\ell=m+1}^j \hat{t}^i_{j \ell} \leq \sum_{\ell = m+1}^J \hat{s}^i_{\ell}$ because $(\hat{s}^i, \hat{t}^i)$ satisfies \eqref{onepool-VI} by the induction assumption, we have $\hat{t}^i_{jm} \geq \hat{s}^i_m, \ \forall j \in A_m \setminus \{j^*\}$.
Then, $\hat{s}^i_m - \pi^i \leq \sum_{j \in A_m \setminus \{j^*\}} (\hat{t}^i_{jm} - \pi^i) + (\hat{t}^i_{j^* m} - \pi^i)$, where $\hat{t}^i_{jm} - \pi^i \geq 0$ for all $j \in A_m \setminus \{j^*\}$ and $\hat{t}^i_{j^* m} - \pi^i \geq 0$ by the definition of $\pi^i$. It follows that \eqref{H-onepool-3} is satisfied:
\begin{align*}
& \hat{s}^{i+1}_k = \hat{s}^i_m - \pi^i \leq \sum_{j \in A_m} (\hat{t}^i_{jm}-\pi^i) \leq \sum_{j=m}^J \hat{t}^{i}_{jm} - \sum_{j \in A_m} \pi^i= \sum_{j=k}^J \hat{t}^{i+1}_{jk}.
\end{align*}
In case 4, \eqref{H-onepool-3} is satisfied because $\hat{s}^{i}_k \leq \sum_{j=k}^J \hat{t}^i_{jk}$, which follows from the induction assumption that $(\hat{s}^i, \hat{t}^i)$ satisfies \eqref{H-onepool-3}.

\noindent
\textbf{(Property 1)}
By construction, $(s^i, t^i)$ produced in line 5 satisfies constraints \eqref{H-onepool-1}, \eqref{H-onepool-2}, and \eqref{H-onepool-6}.
It remains to show constraints \eqref{H-onepool-3}--\eqref{H-onepool-5}, which reduce to the following inequalities as we set $s^i_m = 1$:
\begin{subequations}
	\begin{align}
	& 1 \leq \sum_{j=m}^J t^i_{j m}, \label{integer_t-1} \\
	& t^i_{j \ell} \leq 0, \ \forall \ell \in [m+1, J]_{\mathbb{Z}}, \forall j \in [\ell,J]_{\mathbb{Z}}. \label{integer_t-2}  
	\end{align}
\end{subequations}
To that end, we note from line 4 that for all $j \in [J]$, $\text{Index}(j) \in [0,m]_{\mathbb{Z}}$. 
This implies that, in line 5, $t^i_{j\ell} = 0$ for all $\ell \in [m+1, J]_{\mathbb{Z}}$ and $j \in [\ell,J]_{\mathbb{Z}}$. This proves \eqref{integer_t-2}.
To see \eqref{integer_t-1}, we notice that the construction of $(s^i, t^i)$ implies that $\hat{s}^i_m > 0$ in line 3. Then, Property 3 implies that $\sum_{j=m}^J \hat{t}^i_{jm} \geq \hat{s}^i_m > 0$. It follows that, in line 5, we are able to find a non-empty $A_m$ and so there exists a $j \in [m,J]_{\mathbb{Z}}$ such that $t^i_{jm} = 1$. This proves \eqref{integer_t-1}.

\noindent
\textbf{(Property 2)}
To show that weight $\pi^i > 0$, it suffices to check that each term in the definition of $\pi^i$ in line 6 is positive.
First, by Property 3 and the definition of $B_m$ in line 4, we have $\sum_{\ell=m}^j \hat{t}^i_{j \ell} < \sum_{\ell=m}^J \hat{s}^i_{\ell}, \ \forall j \in B_m$, which implies that $\sum_{\ell=m}^J \hat{s}^i_{\ell} - M_m = \sum_{\ell=m}^J \hat{s}^i_{\ell} - \max_{j \in B_m} \{ \sum_{\ell=m}^j \hat{t}^i_{j \ell} \} > 0$.
Second, $\hat{s}^i_m > 0$ due to line 3. 
It remains to show that $\hat{t}^i_{j,\text{Index}(j)} > 0, \ \forall j \in [J]:\text{Index}(j) \neq 0$.
By noting that $[J] = A_m \cup B_m \cup [m-1]$, we consider the following three cases: (i) $j \in B_m \cup [m-1]$, (ii) $j = j^*$, and (iii) $j \in A_m \setminus \{j^*\} \neq \emptyset$.
In case (i), $\text{Index}(j) \neq 0$ implies that the set $\{k \in [\min\{j,m-1\}] : \hat{t}^i_{jk} > 0 \} \neq \emptyset$, and thus $\hat{t}^i_{j,\text{Index}(j)} > 0$.
In case (ii), we have $\hat{t}^i_{j^*,\text{Index}(j^*)} = \hat{t}^i_{j^* m} > 0$ by the definition of $j^*$.
In case (iii), $\sum_{\ell=m}^j \hat{t}^i_{j \ell} = \sum_{\ell=m}^J \hat{s}^i_{\ell}$ as $j \in A_m \setminus \{j^*\}$.
Since $\sum_{\ell=m+1}^j \hat{t}^i_{j \ell} \leq \sum_{\ell=m+1}^J \hat{s}^i_{\ell}$ by Property 3, $\hat{t}^i_{jm} \geq \hat{s}^i_m > 0$. 
This indicates that $\hat{t}^i_{j,\text{Index}(j)} = \hat{t}^i_{jm} > 0$.
Finally, $\pi^i \leq 1$ because each term in the definition of $\pi^i$ is less than or equal to $1$ by construction.

\noindent
\textbf{(Property 4)}
We claim that, for all $m \in [J]$, the while loop in line 3 terminates in a finite number of iterations, i.e., $\hat{s}^i_m$ becomes zero after a finite number of updates in line 7. It follows from this claim that Algorithm \ref{algo:convex} terminates in a finite number of iterations. Accordingly, as $\hat{s}^n_k = 0, \forall k \in [J]$ at termination, we have $\hat{t}^n_{jk} = 0, \forall j \in [J], \forall k \in [j]$ because $(\hat{s}^n, \hat{t}^n)$ satisfies constraints \eqref{onepool-VI} by Property 3.

To see the claim, we consider the value of $\pi^i$ returned by the minimum comparison in line 6. 
In particular, we distinguish between two cases: (A) $\pi^i=\hat{s}^i_m$ or $\pi^i = \hat{t}^i_{j,\text{Index}(j)}$ for a $j \in [J]$ with $\text{Index}(j) \neq 0$, and (B) $\pi^i = \sum_{\ell=m}^J \hat{s}^i_{\ell} - M_m$, and in addition, $\pi^i < \hat{s}^i_m$ and $\pi^i < \hat{t}^i_{j,\text{Index}(j)}$ for all $j \in [J]$ with $\text{Index}(j) \neq 0$. We notice that $\hat{s}^i_m$ will be updated to zero in line 7 after case (A) takes place for a finite number of times. This is because we only have a finite number of positive values among $\{\hat{s}^i_k: k \in [m]\}$ and $\{\hat{t}^i_{jk}: k \in [m], j \in [k, J]_{\mathbb{Z}}\}$.
Now suppose that case (B) takes place. In this case, as $\pi^i = \sum_{\ell=m}^J \hat{s}^i_{\ell} - M_m < \hat{s}^i_m$, we know that $M_m > 0$ and so there exists a $j \in B_m$ such that $M_m = \sum_{\ell=m}^j \hat{t}^i_{j\ell}$. Then, after the update in line 7, we have $\sum_{\ell=m}^j \hat{t}^{i+1}_{j\ell} = \sum_{\ell=m}^j \hat{t}^i_{j\ell} = M_m$ and $\sum_{\ell=m}^J \hat{s}^{i+1}_{\ell} = \sum_{\ell=m}^J \hat{s}^i_{\ell} - \pi^i = M_m$. It follows that $\sum_{\ell=m}^j \hat{t}^{i+1}_{j\ell} = \sum_{\ell=m}^J \hat{s}^{i+1}_{\ell}$, increasing the cardinality of $A_m$ by $1$. 
Since $|A_m| \leq J-m+1$, this indicates that case (B) can take place for at most $(J-m+1)$ times before case (A) takes place. Therefore, the claim holds valid.

\noindent
\textbf{(Property 5)}
By construction of the binary point in line 5, vectors $\hat{s}^{i+1}$ and $\hat{s}^i$ are only different in the $m^{\mbox{\tiny th}}$ entry, i.e., $\hat{s}^{i+1}_m = \hat{s}^i_m - \pi^i$ and $\hat{s}^{i+1}_k = \hat{s}^i_k$ for all $k \neq m$. It follows that in line 10, i.e., when the for loop terminates, $0 < \sum_{i=1}^{n-1} \pi^i = \sum_{m=1}^J \hat{s}_m \leq 1$, where the first inequality is due to Property 2 and the last inequality is due to Property 3. It follows that $\pi^n \equiv 1 - \sum_{i=1}^{n-1} \pi^i \in [0, 1)$. \qed

\noindent We conclude this section with the following example.
\begin{example}
Consider an example of $1$ pool that covers $4$ units, i.e., $I=1$, $P_1=[4]$, and $J(1)=4$.
Suppose that we are given the following fractional solution $(\hat{s},\hat{t})$:
\begin{center}
	\begin{tabular}{|l|l|l|l|}
		\hline	
		$\hat{t}^1_{11} = 0.6$ &    &   &       \\ \hline
		$\hat{t}^1_{21} = 0$ & $\hat{t}^1_{22} = 0.2$   &   &       \\ \hline
		$\hat{t}^1_{31} = 0.2$ & $\hat{t}^1_{32} = 0.1$   & $\hat{t}^1_{33} = 0.3$  &    \\ \hline
		$\hat{t}^1_{41} = 0$ & $\hat{t}^1_{42} = 0.5$   & $\hat{t}^1_{43} = 0$  &  $\hat{t}^1_{44} = 0.1$    \\ \hline
		$\hat{s}^1_1= 0$ & $\hat{s}^1_2= 0.2$   & $\hat{s}^1_3= 0.3$  &  $\hat{s}^1_4 = 0.1$   \\ \hline
	\end{tabular}
\end{center}
Then, Algorithm \ref{algo:convex} finds $3$ binary points $(s^1,t^1), (s^2,t^2), (s^3,t^3)$ and their weights $\pi^1$, $\pi^2$, $\pi^3$ as follows:\\
\noindent
\textbf{First iteration:} $i = 1$, $m = 2$, $j^*=2$, $A_2= \{2,4 \}$, $B_2=\{3\}$, $M_2=0.4$, \\
$\text{Index}(1)=1,\text{Index}(2)=2,\text{Index}(3)=1,\text{Index}(4)=2$, \\
$s^1_{2} =t^1_{11}=t^1_{22}=t^1_{31}=t^1_{42} = 1$, and \\ 
$\pi^1 = \{ 0.2, 0.2, 0.6, 0.2, 0.2, 0.5  \} = 0.2$.
\begin{center}
	\begin{tabular}{|l|l|l|l|}
		\hline	
		$\hat{t}^2_{11} = 0.4$ 	&    &   &       \\ \hline
		$\hat{t}^2_{21} = 0$ 	& $\hat{t}^2_{22} = 0$   &   &       \\ \hline
		$\hat{t}^2_{31} = 0$ 	& $\hat{t}^2_{32} = 0.1$   	& $\hat{t}^2_{33} = 0.3$  &    \\ \hline
		$\hat{t}^2_{41} = 0$ 	& $\hat{t}^2_{42} = 0.3$   	& $\hat{t}^2_{43} = 0$  	&  $\hat{t}^2_{44} = 0.1$    \\ \hline
		$\hat{s}^2_1= 0$ 		& $\hat{s}^2_2= 0$   		& $\hat{s}^2_3= 0.3$  		&  $\hat{s}^2_4 = 0.1$   \\ \hline
	\end{tabular}
\end{center}
\textbf{Second iteration:} $i = 2$, $m = 3$, $j^*=3$, $A_3= \{3 \}$, $B_3=\{4\}$, $M_3=0.1$, \\
$\text{Index}(1)=1,\text{Index}(2)=0,\text{Index}(3)=3,\text{Index}(4)=2$, \\
$s^2_{3}=t^2_{11}=t^2_{33}=t^2_{42}= 1$, and \\ 
$\pi^2 = \{ 0.3, 0.3, 0.4, 0.3, 0.3  \} = 0.3$.
\begin{center}
	\begin{tabular}{|l|l|l|l|}
		\hline	
		$\hat{t}^3_{11} = 0.1$ 	&    &   &       \\ \hline
		$\hat{t}^3_{21} = 0$ 	& $\hat{t}^3_{22} = 0$   &   &       \\ \hline
		$\hat{t}^3_{31} = 0$ 	& $\hat{t}^3_{32} = 0.1$   	& $\hat{t}^3_{33} = 0$  &    \\ \hline
		$\hat{t}^3_{41} = 0$ 	& $\hat{t}^3_{42} = 0$   	& $\hat{t}^3_{43} = 0$  	&  $\hat{t}^3_{44} = 0.1$    \\ \hline
		$\hat{s}^3_1= 0$ 		& $\hat{s}^3_2= 0$   		& $\hat{s}^3_3= 0$  		&  $\hat{s}^3_4 = 0.1$   \\ \hline
	\end{tabular}
\end{center}
\textbf{Third iteration:} $i = 3$, $m = 4$, $j^*=4$, $A_4= \{4 \}$, $B_4=\emptyset$, $M_4=0$, \\
$\text{Index}(1)=1,\text{Index}(2)=0,\text{Index}(3)=2,\text{Index}(4)=4$, \\
$s^3_{4} =t^3_{11} = t^3_{32}=t^3_{44} = 1$, and \\ 
$\pi^3 = \{ 0.1, 0.1, 0.1, 0.1, 0.1  \} = 0.1$.
\begin{center}
	\begin{tabular}{|l|l|l|l|}
		\hline	
		$\hat{t}^4_{11} = 0$ 	&    &   &       \\ \hline
		$\hat{t}^4_{21} = 0$ 	& $\hat{t}^4_{22} = 0$   &   &       \\ \hline
		$\hat{t}^4_{31} = 0$ 	& $\hat{t}^4_{32} = 0$   	& $\hat{t}^4_{33} = 0$  &    \\ \hline
		$\hat{t}^4_{41} = 0$ 	& $\hat{t}^4_{42} = 0$   	& $\hat{t}^4_{43} = 0$  	&  $\hat{t}^4_{44} = 0$    \\ \hline
		$\hat{s}^4_1= 0$ 		& $\hat{s}^4_2= 0$   		& $\hat{s}^4_3= 0$  		&  $\hat{s}^4_4 = 0$   \\ \hline
	\end{tabular}
\end{center}
$\pi^4 = 1 - 0.2 - 0.3 - 0.1 = 0.4$.
Therefore, $(\hat{s}, \hat{t}) = \sum_{i=1}^4 \pi^i (s^i, t^i)$. \qed
\end{example}


\section{Proof of Theorem \ref{thm:integrality-one}} \label{apx-thm:integrality-one}
\proof 
Under Structure \ref{assump:one-pool}, $\max_{(\alpha, \beta) \in \Lambda} F(\alpha, \beta)$ can be reformulated as the following linear program:
\begin{align}
& \max_{(\alpha, \beta) \in \Lambda} F(\alpha, \beta) = \max_{(t,s,r,p) \in \overline{\mathcal{H}}_{\text{\tiny 1}}} \ \sum_{j=1}^J \bigg( c^{\text{\tiny r}}_j r_j + \sum_{k =1}^j c^{\text{\tiny t}}_{jk} t_{jk} \bigg) +  c^{\text{\tiny p}}_1 p + \sum_{k = 1}^{J} c^{\text{\tiny s}}_{1k} s_{k}   \nonumber \\
& = c^{\text{\tiny p}}_1 + \sum_{j=1}^J c^{\text{\tiny r}}_j + \max_{s \geq 0 : \eqref{H-onepool-2} } \bigg[ \sum_{k=1}^J (c^{\text{\tiny s}}_{1k} - c^{\text{\tiny p}}_1) s_k + \bigg\{ \max_{t \geq 0} \sum_{j=1}^{J} \sum_{k=1}^j (c^{\text{\tiny t}}_{jk} - c^{\text{\tiny r}}_j) t_{jk} \bigg\} \bigg] \label{onepool-linearprogram} \\
& \hspace{79mm} \mbox{s.t.} \ \eqref{H-onepool-1}, \eqref{H-onepool-3}, \eqref{onepool-VI}. \nonumber
\end{align}
An optimal solution $(s^*,t^*)$ to \eqref{onepool-linearprogram} is integral as $\overline{\mathcal{H}}_{\text{\tiny 1}}=\text{conv} ({\mathcal{H}}_{\text{\tiny 1}})$ is an integral polyhedron. Specifically, by constraints \eqref{H-onepool-2}, we have $s^* \in \{0\}\cup\{e_k: k \in [J]\}$, where $e_k$ is the $k^{\text{\tiny th}}$ standard basis vector, and for given $s^*$, the inner maximization problem in \eqref{onepool-linearprogram} finds an optimal $t^*$, which is also integral.
We rewrite \eqref{onepool-linearprogram} by discussing two cases of $s^*$: (i) $s^*=0$ and (ii) $s^* \in \{e_k: k \in [J] \}$.

In case (i), we have $t^*_{jk} = 0$ for all $j \in [J]$ and $k \in [j]$ by constraints \eqref{onepool-VI}.
Then, \eqref{onepool-linearprogram} reduces to $c^{\text{\tiny p}}_1 + \sum_{j=1}^J c^{\text{\tiny r}}_j$.

In case (ii), $s_{k^*} = 1$ for some $k^* \in [J]$ and $s_{k} = 0$ for all $k \in [J] \setminus \{k^*\}$.
By constraints \eqref{onepool-VI}, we have $t_{jk} = 0$ if $j,k \geq k^*+1$. 
Then, we can rewrite \eqref{onepool-linearprogram} as follows:
\begin{align*}
& c^{\text{\tiny p}}_1 + \sum_{j=1}^J c^{\text{\tiny r}}_j + (c^{\text{\tiny s}}_{1 k^*} - c^{\text{\tiny p}}_1) + \max_{t \geq 0} \ \sum_{j=1}^{J} \sum_{k=1}^{\min \{j, k^* \} } (c^{\text{\tiny t}}_{jk} - c^{\text{\tiny r}}_j) t_{jk} \\
& \hspace{47mm} \mbox{s.t.} \ \sum_{k=1}^{ \min \{ j, k^* \} } t_{jk} \leq 1, \ \forall j \in [J], \ \ \ 1 \leq \sum_{j = k^*}^J t_{j k^*}.
\end{align*}
By constraint \eqref{H-onepool-1}, $t_{j^*k^*}=1$ for some $j^* \in [k^*,J]_{\mathbb{Z}}$ and $t_{j^* k} = 0$ for all $k \in [j^*] \setminus \{k^*\}$. Then, we recast \eqref{onepool-linearprogram} as
\begin{align*}
& c^{\text{\tiny p}}_1 + \sum_{j=1}^J c^{\text{\tiny r}}_j + (c^{\text{\tiny s}}_{1k^*} - c^{\text{\tiny p}}_1) +  (c^{\text{\tiny t}}_{j^* k^*} - c^{\text{\tiny r}}_{j^*}) + \sum_{j \in [J] \setminus \{j^*\} } \max_{t \geq 0} \bigg\{ \sum_{k=1}^{\min \{j, k^* \} } (c^{\text{\tiny t}}_{jk} - c^{\text{\tiny r}}_j) t_{jk}\bigg| \ \sum_{k=1}^{ \min \{j, k^*\} } t_{jk} \leq 1 \bigg\}  \\
& = \ \sum_{j=1}^J c^{\text{\tiny r}}_j + c^{\text{\tiny s}}_{1k^*}  +  (c^{\text{\tiny t}}_{j^* k^*} - c^{\text{\tiny r}}_{j^*}) + \sum_{j \in [J] \setminus \{j^*\} } \max_{k \in [ \min \{j, k^* \} ] } \ \Big\{ c^{\text{\tiny t}}_{jk} - c^{\text{\tiny r}}_j \Big\}.
\end{align*}
Therefore, we have
\begin{align*}
& \max_{(\alpha,\beta) \in \Lambda} F(\alpha, \beta) =  \\
& \max \Bigg\{ c^{\text{\tiny p}}_1 + \sum_{j=1}^J c^{\text{\tiny r}}_j, \ \ \max_{\substack{j^* \in [J], \\ k^* \in [j^*] }} \bigg\{ \sum_{j=1}^J c^{\text{\tiny r}}_j + c^{\text{\tiny s}}_{1k^*}  +  (c^{\text{\tiny t}}_{j^* k^*} - c^{\text{\tiny r}}_{j^*}) + \sum_{j \in [J] \setminus \{j^*\} } \max_{k \in [\min \{j, k^* \} ] } \ \Big\{ c^{\text{\tiny t}}_{jk} - c^{\text{\tiny r}}_j \Big\} \bigg\}
\Bigg\}.
\end{align*}
This completes the proof. \qed

\section{Proof of Proposition \ref{prop:one-pool-ref}} \label{apx-prop:one-pool-ref}
\proof 
By Theorem \ref{thm:integrality-one}, constraints \eqref{ref-linear-3} are equivalent to
\begin{align*}
\theta \geq & 
\max \Bigg\{
c^{\text{\tiny p}}_1 + \sum_{j=1}^J c^{\text{\tiny r}}_j, \ \max_{\substack{j^* \in [J], \\  k^* \in [j^*]}} \bigg\{ \sum_{j=1}^J c^{\text{\tiny r}}_j + c^{\text{\tiny s}}_{1k^*}  +  (c^{\text{\tiny t}}_{j^* k^*} - c^{\text{\tiny r}}_{j^*}) + \sum_{j \in [J] \setminus \{j^*\} } \bigg( \max_{k \in [\min \{j, k^* \} ] } \ \Big\{ (c^{\text{\tiny t}}_{jk} - c^{\text{\tiny r}}_j) \Big\} \bigg) \bigg\}
\Bigg\}.
\end{align*}
where $c^{\text{\tiny r}}_j$, $c^{\text{\tiny t}}_{jk}$, $c^{\text{\tiny p}}_1$, and $c^{\text{\tiny s}}_{1k}$ are computed by \eqref{ref-linear-5}--\eqref{ref-linear-8}.
Let $\chi_{jk^*} :=  \max_{k \in [ \min \{ j, k^* \} ] } \ \Big\{ (c^{\text{\tiny t}}_{jk} - c^{\text{\tiny r}}_j) \Big\}$ for all $j \in [J]$ and $k^* \in [J]$.
Then, constraints \eqref{ref-linear-3} can be rewritten as follows:
\begin{subequations}
\label{ref-linear-3-1}
\begin{align}
& \theta \geq c^{\text{\tiny p}}_1 + \sum_{j=1}^J c^{\text{\tiny r}}_j, \\
& \theta \geq \sum_{j=1}^J c^{\text{\tiny r}}_j + c^{\text{\tiny s}}_{1k^*} +  (c^{\text{\tiny t}}_{j^* k^*} - c^{\text{\tiny r}}_{j^*}) + \sum_{j \in  [J] \setminus \{j^*\}} \chi_{jk^*}, \ \forall j^* \in [J], \forall k^* \in [j^*],  \\
& \chi_{jk^*} \geq  (c^{\text{\tiny t}}_{jk} - c^{\text{\tiny r}}_j), \ \forall j \in [J], \forall k^* \in [J], \forall k \in [ \min \{j,k^* \}].
\end{align}
\end{subequations}
Moreover, we define the following auxiliary variables:
\begin{subequations}
\label{onepool-aux}
\begin{align}
& \hspace{-0.15cm} \left.\begin{array}{l}
\zeta^{\mbox{\tiny e}}_j \ := \ \bigl[ - \gamma_j w^{\tinyl}_j - \sum_{\ell=1}^{w^{\mbox{\tinyu}}_j - w^{\mbox{\tinyl}}_j}  \varphi_{j\ell} \bigr]_+ \\[0.1cm]
\eta^{\mbox{\tiny e}}_j \ := \ \sup_{\tilde{d}_j \in [d^{\tinyl}_j, d^{\tinyu}_j]_{\mathbb{Z}}} \bigl\{  - \sum_{q=1}^Q\tilde{d}^q_j\rho_{jq} \bigr\}
\end{array}\right\} \ \ \forall j \in [J], \\
& \hspace{-0.15cm} \left.\begin{array}{l}
\zeta^{\mbox{\tiny x}}_{jk} \ := \ \bigl[(- c^{\mbox{\tiny x}}_k - \gamma_j) w^{\tinyl}_j - \sum_{\ell=1}^{w^{\mbox{\tinyu}}_j - w^{\mbox{\tinyl}}_j} ( c^{\mbox{\tiny x}}_k u_{j \ell} + \varphi_{j \ell})\bigr]_+ \\[0.1cm]
\eta^{\mbox{\tiny x}}_{jk} \ := \ \sup_{\tilde{d}_j \in [d^{\tinyl}_j, d^{\tinyu}_j]_{\mathbb{Z}}} \bigl\{ c^{\mbox{\tiny x}}_k \tilde{d}_j - \sum_{q=1}^Q\tilde{d}^q_j\rho_{jq} \bigr\} \\[0.1cm]
\end{array}\right\} \ \ \forall j \in [J], \forall k \in [j], \\
& \ \phi^{\mbox{\tiny e}}_1 \ := \ \bigl[  - \lambda_1 y^{\tinyl}_1 - \sum_{\ell=1}^{y^{\mbox{\tinyu}}_1 - y^{\mbox{\tinyl}}_1} \nu_{1\ell} \bigr]_+, \\
& \ \phi^{\mbox{\tiny x}}_{1k} \ := \ \Bigl[(- c^{\mbox{\tiny x}}_k - \lambda_1) y^{\tinyl}_1 - \sum_{\ell=1}^{y^{\mbox{\tinyu}}_1 - y^{\mbox{\tinyl}}_1} (c^{\mbox{\tiny x}}_k v_{1\ell} + \nu_{1\ell} ) \Bigr]_+, \ \forall k \in [J],
\end{align}
\end{subequations}
Since
$c^{\text{\tiny r}}_j = \zeta^{\mbox{\tiny e}}_j + \eta^{\mbox{\tiny e}}_j$, 
$c^{\text{\tiny t}}_{jk} = \zeta^{\mbox{\tiny x}}_{jk} + \eta^{\mbox{\tiny x}}_{jk}$, 
$c^{\text{\tiny p}}_1 = \phi^{\mbox{\tiny e}}_1$, and 
$c^{\text{\tiny s}}_{1k} = \phi^{\mbox{\tiny x}}_{1k}$,
constraints \eqref{ref-linear-3-1} can be rewritten as follows:
\begin{subequations}  
	\label{ref-linear-3-2}
	\begin{align}
	& \theta \geq \phi^{\mbox{\tiny e}}_1 + \sum_{j=1}^J (\zeta^{\mbox{\tiny e}}_j + \eta^{\mbox{\tiny e}}_j), \\
	& \theta \geq \sum_{j=1}^J (\zeta^{\mbox{\tiny e}}_j + \eta^{\mbox{\tiny e}}_j) + \phi^{\mbox{\tiny x}}_{1k^*} + (\zeta^{\mbox{\tiny x}}_{j^* k^*} + \eta^{\mbox{\tiny x}}_{j^* k^*}) - (\zeta^{\mbox{\tiny e}}_{j^*} + \eta^{\mbox{\tiny e}}_{j^*}) + \sum_{j \in  [J] \setminus \{j^*\}} \chi_{jk^*}, \ \forall j^* \in [J], \forall k^* \in [j^*], \\
	& \chi_{jk^*} \geq  (\zeta^{\mbox{\tiny x}}_{jk} + \eta^{\mbox{\tiny x}}_{jk}) - (\zeta^{\mbox{\tiny e}}_j + \eta^{\mbox{\tiny e}}_j), \ \forall j \in [J], \forall k^* \in [j], \forall k \in [\min\{j,k^*\}].
	\end{align}
\end{subequations}
Replacing constraints \eqref{ref-linear-3} with \eqref{onepool-aux}--\eqref{ref-linear-3-2} in the formulation \eqref{ref-linear-1}--\eqref{ref-linear-4} leads to the claimed reformulation of (DRNS). This completes the proof. \qed

\section{Proof of Proposition \ref{prop:exact-disjoint}} \label{apx-prop:exact-disjoint}
We start by proving the following technical lemma.
\begin{lemma} \label{lem:technical}
Consider sets $A_i \subseteq \mathbb{R}^{k_i}$ for all $i \in [I]$ and let $A := \Pi_{i=1}^I A_i$. Then, $\mbox{conv}(A) = \Pi_{i=1}^I \mbox{conv}(A_i)$.
\end{lemma}
\proof First, as $A = \Pi_{i=1}^I A_i$ and $A_i \subseteq \mbox{conv}(A_i)$, we have $A \subseteq \Pi_{i=1}^I \mbox{conv}(A_i)$ and hence $\mbox{conv}(A) \subseteq \Pi_{i=1}^I \mbox{conv}(A_i)$. Second, to show that $\Pi_{i=1}^I \mbox{conv}(A_i) \subseteq \mbox{conv}(A)$, we pick any $a \in \Pi_{i=1}^I \mbox{conv}(A_i)$ and prove that $a \in \mbox{conv}(A)$. To this end, we denote $a := [a_1, \ldots, a_I]^{\top}$, where $a_i \in \mbox{conv}(A_i)$ for all $i \in [I]$. Then, for all $i \in [I]$, there exist $\{\lambda_i^{n_i}\}_{n_i = 1}^{N_i}$ and $\{a_i^{n_i}\}_{n_i = 1}^{N_i}$ such that each $\lambda_i^{n_i} \geq 0$, each $a_i^{n_i} \in A_i$, $\sum_{n_i=1}^{N_i} \lambda_i^{n_i} = 1$, and $\sum_{n_i=1}^{N_i} \lambda_i^{n_i} a_i^{n_i} = a_i$ for all $i \in [I]$.

Denote set $\mathcal{N} := \{(n_1, \ldots, n_I): n_i \in [N_i], \ \forall i \in [I]\}$, vector $a^n := [a_1^{n_1}, \ldots, a_I^{n_I}]^{\top}$ for all $n := [n_1, \ldots, n_I]^{\top} \in \mathcal{N}$, and scalar $\lambda^n := \Pi_{i=1}^I \lambda_i^{n_i}$ for all $n \in \mathcal{N}$. Then, $\lambda^n \geq 0$ and $a^n \in A$ for all $n \in \mathcal{N}$. In addition, $\sum_{n \in \mathcal{N}} \lambda^n = \sum_{n \in \mathcal{N}} \Pi_{i=1}^I \lambda_i^{n_i} = (\lambda_1^1 + \cdots + \lambda_1^{N_1}) (\lambda_2^1 + \cdots + \lambda_2^{N_2}) \cdots (\lambda_I^1 + \cdots + \lambda_I^{N_I}) = 1$. Furthermore, for all $i \in [I]$, we have
\begin{align*}
\sum_{n \in \mathcal{N}} \lambda^n a_i^{n_i} \ = & \ \sum_{m_i = 1}^{N_i} \sum_{n \in \mathcal{N}: n_i = m_i} \lambda^n a_i^{n_i} \\
= \ & \ \sum_{m_i = 1}^{N_i} a_i^{m_i} \sum_{n \in \mathcal{N}: n_i = m_i} \lambda^n \\
= \ & \ \sum_{m_i = 1}^{N_i} a_i^{m_i} \lambda_i^{m_i} \\
= \ & \ a_i,
\end{align*}
where third equality is because, for fixed $i \in [I]$ and $m_i \in [N_i]$, $\sum_{n \in \mathcal{N}: n_i = m_i} \lambda^n = (\lambda_1^1 + \cdots + \lambda_1^{N_1}) \cdots (\lambda_{m_i-1}^1 + \cdots + \lambda_{m_i-1}^{N_{m_i-1}}) \lambda_i^{m_i} (\lambda_{m_i+1}^1 + \cdots + \lambda_{m_i+1}^{N_{m_i+1}}) \cdots (\lambda_I^1 + \cdots + \lambda_I^{N_I}) = \lambda_i^{m_i}$, and the last inequality follows from the definitions of $\{\lambda_i^{m_i}\}_{m_i = 1}^{N_i}$ and $\{a_i^{m_i}\}_{m_i = 1}^{N_i}$. It follows that $a \equiv [a_1, \ldots, a_I]^{\top} = \sum_{n \in \mathcal{N}} \lambda^n [a_1^{n_1}, \ldots, a_I^{n_I}]^{\top} \equiv \sum_{n \in \mathcal{N}} \lambda^n a^n$ and hence $a \in \mbox{conv}(A)$. This completes the proof. \qed

\noindent We are now ready to present the main proof of this section.\\
\noindent{\em Proof of Proposition \ref{prop:exact-disjoint}:} 
Since $\mathcal{H}_{\text{\tiny D}}$ is separable in index $i$, we have $\mathcal{H}_{\text{\tiny D}} = \prod_{i=1}^I \mathcal{H}^{\text{\tiny D}}_i$, where each $\mathcal{H}^{\text{\tiny D}}_i$ is defined as
\begin{subequations}
	\begin{align}
	\mathcal{H}^{\text{\tiny D}}_i := \bigg\{ (t,s_i,r,p_i): \
	& r_j + \sum_{k=1}^j t_{jk} = 1, \ \forall j \in P_i, \ \ \ p_i + \sum_{k=1}^{J(i)} s_{ik} = 1, \label{H-disjoint-0} \\
	& \sum_{k = 1}^j t_{jk} \leq 1, \ \forall j \in P_i, \label{H-disjoint-1} \\
	& \sum_{k=1}^{J(i)} s_{ik} \leq 1, \label{H-disjoint-2} \\
	& s_{ik} \leq \sum_{j \in P_i: j \geq k} t_{jk}, \ \forall k \in [J(i)], \label{H-disjoint-3} \\
	& s_{ik} + t_{j \ell} \leq 1, \ \forall k \in [J(i)], \forall j \in P_i: j \geq k+1, \forall \ell \in [k+1,J(i)]_{\mathbb{Z}}, \nonumber \\
	& t_{jk} \leq \sum_{\ell = 1}^{J(i)} s_{i \ell}, \  \forall j \in P_i, \forall k \in [j], \nonumber \\
	& s_{ik} \in \mathbb{B}, \  \forall k \in [J(i)], \ \ t_{jk} \in \mathbb{B}, \ \forall j \in P_i, \forall k \in [j] \nonumber
	\bigg\}.
	\end{align}	
\end{subequations}
Following a similar proof in Section \ref{subsec:one-pool},
we can show that incorporating inequalities \eqref{VI-general} produces the convex hull of $\mathcal{H}^{\text{\tiny D}}_i$, i.e., 
\begin{align}
\overline{\mathcal{H}}^{\text{\tiny D}}_i := \bigg\{ (t,s_i,r,p_i): \
& \eqref{H-disjoint-0}-\eqref{H-disjoint-3}, \nonumber \\
& \sum_{\ell = k }^j t_{j \ell} \leq \sum_{\ell = k}^{J(i)} s_{i \ell}, \ \forall j \in P_i, \forall k \in [j], \label{H-disjoint-VI} \\
& s_{ik} \in \mathbb{R}_+, \forall k \in [J(i)], \ \ t_{jk} \in \mathbb{R}_+, \ \forall j \in P_i, \forall k \in [j] \nonumber
\bigg\}.
\end{align}
Let $\overline{\mathcal{H}}_{\text{\tiny D}} = \prod_{i=1}^I \overline{\mathcal{H}}^{\text{\tiny D}}_i$.
Then, it follows from Lemma \ref{lem:technical} that $\overline{\mathcal{H}}_{\text{\tiny D}} = \mbox{conv} ( \mathcal{H}_{\text{\tiny D}})$.
Following a similar proof as that of Theorem \ref{thm:integrality-one}, we have
\begin{align}
& \max_{(\alpha, \beta) \in \Lambda} F(\alpha, \beta) = \max_{(t,s,r,p) \in \overline{\mathcal{H}}_{\text{\tiny D}}} \ \sum_{j=1}^J \bigg( c^{\text{\tiny r}}_j r_j + \sum_{k = 1}^j c^{\text{\tiny t}}_{jk} t_{jk} \bigg) + \sum_{i=1}^I \bigg( c^{\text{\tiny p}}_i p_i + \sum_{k = 1}^{J(i)} c^{\text{\tiny s}}_{ik} s_{ik} \bigg)  \nonumber \\
& = \sum_{i=1}^I \bigg[ \max_{(t,s_i,r,p_i) \in \overline{\mathcal{H}}^{\text{\tiny D}}_i} \ \sum_{j \in P_i} \bigg( c^{\text{\tiny r}}_j r_j + \sum_{k =1}^j c^{\text{\tiny t}}_{jk} t_{jk} \bigg) +  \bigg( c^{\text{\tiny p}}_i p_i + \sum_{k = 1}^{J(i)} c^{\text{\tiny s}}_{ik} s_{ik} \bigg) \bigg]  \nonumber \\
& = \sum_{i=1}^I \bigg[  c^{\text{\tiny p}}_i + \sum_{j \in P_i} c^{\text{\tiny r}}_j + \max_{s_i \geq 0: \eqref{H-disjoint-2}} \ \sum_{k = 1}^{J(i)} (c^{\text{\tiny s}}_{ik}-c^{\text{\tiny p}}_i) s_{ik} + \bigg( \max_{t \geq 0} \sum_{j \in P_i} \sum_{k =1}^j (c^{\text{\tiny t}}_{jk} - c^{\text{\tiny r}}_j) t_{jk} \bigg)  \bigg]  \label{disjoint-LP} \\
& \hspace{88mm} \mbox{s.t.} \ \eqref{H-disjoint-1}, \eqref{H-disjoint-3}, \eqref{H-disjoint-VI}. \nonumber
\end{align}
For every $i \in [I]$, an optimal solution $(s^*_i,t^*)$ to \eqref{disjoint-LP} is integral as $\overline{\mathcal{H}}^{\text{\tiny D}}_i=\text{conv} ({\mathcal{H}}^{\text{\tiny D}}_i)$ is an integral polyhedron.
Following a similar proof in Theorem \ref{thm:integrality-one}, we have
\begin{align*}
& \max_{(\alpha,\beta) \in \Lambda} F(\alpha, \beta) = \sum_{i=1}^I \Bigg[ \max \Bigg\{ c^{\text{\tiny p}}_i + \sum_{j \in P_i} c^{\text{\tiny r}}_j, \\ 
& \hspace{15mm} \max_{j^* \in P_i, k^* \in [j^*] } \bigg\{ \sum_{j \in P_i} c^{\text{\tiny r}}_j + c^{\text{\tiny s}}_{ik^*}  +  (c^{\text{\tiny t}}_{j^* k^*} - c^{\text{\tiny r}}_{j^*}) + \sum_{j \in P_i \setminus \{j^*\} } \bigg( \max_{k \in [\min \{j, k^* \} ] } \ \Big\{ c^{\text{\tiny t}}_{jk} - c^{\text{\tiny r}}_j \Big\} \bigg) \bigg\}
\Bigg\}
\Bigg].
\end{align*}
We define $\chi_{jk^*} :=  \max_{k \in [ \min \{ j, k^* \} ] } \ \Big\{ (c^{\text{\tiny t}}_{jk} - c^{\text{\tiny r}}_j) \Big\}$ for all $j \in [J]$ and $k^* \in [J]$.
Then, we rewrite \eqref{disjoint-LP} as follows:
\begin{align*}
\min \ & \sum_{i=1}^I \theta_i \\
\mbox{s.t.} \ 
& \theta_i \geq c^{\text{\tiny p}}_i + \sum_{j \in P_i} c^{\text{\tiny r}}_j, \ \forall i \in [I], \\
& \theta_i \geq \sum_{j \in P_i} c^{\text{\tiny r}}_j + c^{\text{\tiny s}}_{ik^*} + (c^{\text{\tiny t}}_{j^* k^*} - c^{\text{\tiny r}}_{j^*}) + \sum_{j \in P_i \setminus \{j^*\} } \chi_{jk^*}, \ \forall i \in [I], \forall j^* \in P_i, \forall k^* \in [j^*], \\
& \chi_{jk^*} \geq  (c^{\text{\tiny t}}_{jk} - c^{\text{\tiny r}}_j), \ \forall j \in [J], \forall k^* \in [J], \forall k \in [\min \{ j, k^* \} ].
\end{align*}
Moreover, we define the following auxiliary variables:
\begin{subequations}
\label{disjoint-aux}
\begin{align}
& \hspace{-0.15cm} \left.\begin{array}{l}
\zeta^{\mbox{\tiny e}}_j \ := \ \Bigl[ - \gamma_j w^{\tinyl}_j - \sum_{\ell=1}^{w^{\mbox{\tinyu}}_j - w^{\mbox{\tinyl}}_j} \varphi_{j\ell} \Bigr]_+ \\[0.1cm]
\eta^{\mbox{\tiny e}}_j \ := \ \sup_{\tilde{d}_j \in [d^{\tinyl}_j, d^{\tinyu}_j]_{\mathbb{Z}}} \bigl\{ - \sum_{q=1}^Q\tilde{d}^q_j\rho_{jq} \bigr\}
\end{array}\right\} \ \ \forall j \in [J], \\
& \hspace{-0.15cm} \left.\begin{array}{l}
\zeta^{\mbox{\tiny x}}_{jk} \ := \ \Bigl[(- c^{\mbox{\tiny x}}_k - \gamma_j) w^{\tinyl}_j - \sum_{\ell=1}^{w^{\mbox{\tinyu}}_j - w^{\mbox{\tinyl}}_j} ( c^{\mbox{\tiny x}}_k u_{j \ell} + \varphi_{j \ell})\Bigr]_+ \\[0.1cm]
\eta^{\mbox{\tiny x}}_{jk} \ := \ \sup_{\tilde{d}_j \in [d^{\tinyl}_j, d^{\tinyu}_j]_{\mathbb{Z}}} \bigl\{ c^{\mbox{\tiny x}}_k \tilde{d}_j - \sum_{q=1}^Q\tilde{d}^q_j\rho_{jq} \bigr\} \\[0.1cm]
\end{array}\right\} \ \ \forall j \in [J], \forall k \in [j], \\
& \ \phi^{\mbox{\tiny e}}_i \ := \ \Bigl[  - \lambda_i  y^{\tinyl}_i - \sum_{\ell=1}^{y^{\mbox{\tinyu}}_i - y^{\mbox{\tinyl}}_i} \nu_{i\ell} \Bigr]_+, \ \forall i \in [I], \\
& \ \phi^{\mbox{\tiny x}}_{ik} \ := \ \Bigl[(- c^{\mbox{\tiny x}}_k - \lambda_i) y^{\tinyl}_i - \sum_{\ell=1}^{y^{\mbox{\tinyu}}_i - y^{\mbox{\tinyl}}_i} (c^{\mbox{\tiny x}}_k v_{i\ell} + \nu_{i\ell} ) \Bigr]_+, \ \forall i \in [I], \forall k \in [J(i)],
\end{align}
\end{subequations}
Since
$c^{\text{\tiny r}}_j = \zeta^{\mbox{\tiny e}}_j + \eta^{\mbox{\tiny e}}_j$, 
$c^{\text{\tiny t}}_{jk} = \zeta^{\mbox{\tiny x}}_{jk} + \eta^{\mbox{\tiny x}}_{jk}$, 
$c^{\text{\tiny p}}_i = \phi^{\mbox{\tiny e}}_i$, and 
$c^{\text{\tiny s}}_{ik} = \phi^{\mbox{\tiny x}}_{ik}$, we can rewrite \eqref{disjoint-LP} as follows:
\begin{align*}
\min \ & \sum_{i=1}^I \theta_i \\
\mbox{s.t.} \ 
& \eqref{disjoint-aux}, \\
& \theta_i \geq \phi^{\mbox{\tiny e}}_i + \sum_{j \in P_i} (\zeta^{\mbox{\tiny e}}_j + \eta^{\mbox{\tiny e}}_j), \ \forall i \in [I], \\
& \theta_i \geq \sum_{j \in P_i} (\zeta^{\mbox{\tiny e}}_j + \eta^{\mbox{\tiny e}}_j) + \phi^{\mbox{\tiny x}}_{ik^*} + (\zeta^{\mbox{\tiny x}}_{j^* k^*} + \eta^{\mbox{\tiny x}}_{j^* k^*}) - (\zeta^{\mbox{\tiny e}}_{j^*} + \eta^{\mbox{\tiny e}}_{j^*}) + \sum_{j \in P_i \setminus \{j^*\} } \chi_{jk^*}, \ \forall i \in [I], \forall j^* \in P_i, \forall k^* \in [j^*], \\
& \chi_{jk^*} \geq (\zeta^{\mbox{\tiny x}}_{jk} + \eta^{\mbox{\tiny x}}_{jk}) - (\zeta^{\mbox{\tiny e}}_j + \eta^{\mbox{\tiny e}}_j), \ \forall j \in [J], \forall k^* \in [J], \forall k \in [\min\{j,k^*\}].
\end{align*}
This completes the proof. \qed

\section{Example for $\mathcal{H}$ being not integral under the Chained-Pool Structure} \label{apx:example-chained-pool}
{\color{black}
\begin{example}
Consider an example of chained pools with $I=J=3$, $P_1=\{1,3\}$, $P_2=\{2,3\}$, and $P_3=\{1,2\}$ (see Figure~\ref{Fig:chained-pool}). The polyhedron $\overline{\mathcal{H}}$, which is obtained by incorporating valid inequalities \eqref{VI-general} and relaxing binary restrictions in $\mathcal{H}$, reads
\begin{subequations}
\begin{align}
	\overline{\mathcal{H}} = \big\{ & (t,s,r,p): \nonumber \\
	& r_1 + t_{11} = r_2 + t_{21} + t_{22} = r_3 + t_{31} + t_{32} + t_{33} = 1, \nonumber \\
	& p_1 + s_{11} + t_{32} + t_{33} = p_2 + s_{21} + s_{22} + t_{33} = p_3 + s_{31} + t_{22} = 1, \nonumber \\
	& s_{11} + t_{32} + t_{33} \leq 1, \ \ s_{21} + s_{22} + t_{33} \leq 1, \ \ s_{31} + t_{22} \leq 1, \label{binding-1} \\ 
	& s_{11} \leq t_{11} + t_{31}, \ s_{21} \leq t_{21} + t_{31}, \ s_{22} \leq t_{22} + t_{32}, \ s_{31} \leq t_{11} + t_{21}, \label{binding-2}\\
	& t_{11} \leq s_{11} + t_{32} + t_{33}, \  t_{31} \leq s_{11}, \nonumber \\ 
	& t_{21} + t_{22} \leq s_{21} + s_{22} + t_{33}, \ t_{31} + t_{32} \leq s_{21} + s_{22}, \ t_{22} \leq s_{22} + t_{33}, \ \ t_{32}  \leq s_{22}, \nonumber \\
	& t_{11} \leq s_{31} + t_{22}, \ t_{21}  \leq s_{31}, \nonumber \\
	& t_{11}, t_{21}, t_{22}, t_{31}, t_{32}, t_{33}, s_{11}, s_{21}, s_{22}, s_{31} \geq 0 \nonumber
	\big\}.
\end{align}
We observe that $\overline{\mathcal{H}}$ is $10$-dimensional. Hence, replacing inequalities in constraints \eqref{binding-1} and \eqref{binding-2} with equalities yields the following extreme point:
\begin{align*}
(s_{11}, s_{21}, s_{22}, s_{31}, t_{11}, t_{21}, t_{22}, t_{31}, t_{32}, t_{33}, r_1,r_2,r_3,p_1,p_2,p_3) = \bigg(\frac{1}{2}, 0,\frac{1}{2},\frac{1}{2},\frac{1}{2}, 0,\frac{1}{2},0,0,\frac{1}{2},\frac{1}{2},\frac{1}{2},\frac{1}{2},0,0,0 \bigg),
\end{align*}
which is fractional. Therefore, $\mathcal{H}$ is not integral, i.e., $\overline{\mathcal{H}} \neq \text{conv} (\mathcal{H})$.
\end{subequations}
\end{example}
}

\section{Proof of Theorem \ref{thm:longest-path}} \label{apx-thm:longest-path}
\proof 
Each trajectory of states $t_{[1]}, (t_{[1]}, t_{[2]}), \ldots, (t_{[1]}, t_{[I]})$ in the DP corresponds to an $\texttt{S}$-$\texttt{T}$ path in the network $(\mathcal{N}, \mathcal{A})$, where the objective function value of the trajectory $V_{I} (t_{[1]}, t_{[I]}) + \sum_{k\preceq 1 \curlyvee I} (c^{\text{\tiny s}}_{Ik}  - c^{\text{\tiny p}}_I) h_k (t_{[I]}, t_{[1]} )$ equals the length of the $\texttt{S}$-$\texttt{T}$ path by definition of the arc lengths $c_{[m,n]}$. 
Likewise, each $\texttt{S}$-$\texttt{T}$ path in $(\mathcal{N}, \mathcal{A})$ corresponds to a trajectory of states in the DP and the length of the path equals the objective function value of the trajectory. 
This proves that the length of the longest path in $(\mathcal{N}, \mathcal{A})$ equals $\max_{(\alpha, \beta) \in \Lambda}F(\alpha, \beta)$ and completes the proof. \qed

\section{Proof of Proposition \ref{prop:two-longest}} \label{apx-prop:two-longest}
\proof Taking the dual of the longest-path formulation yields
\begin{align*}
\min_{\pi} \ & \pi_{\texttt{S}} - \pi_{\texttt{T}} + \sum_{i=1}^I (c^{\text{\tiny r}}_i + c^{\text{\tiny p}}_i) \\
\mbox{s.t.} \
& \pi_{\texttt{S}} - \pi_{t_{[1]}} \geq (c^{\text{\tiny t}}_{11} - c^{\text{\tiny r}}_1) t_{11}, \ \forall t_{[1]} \in \mathcal{B}_1, \\
& \pi_{t_{[1]}} - \pi_{(t_{[1]}, e_{2\ell})} \geq 
\begin{cases}
( c^{\text{\tiny s}}_{11} - c^{\text{\tiny p}}_{1} ) t_{11} & \text{if } \ell = 0 \\
(c^{\text{\tiny t}}_{2\ell} - c^{\text{\tiny r}}_{2} ) + (c^{\text{\tiny s}}_{1\ell} - c^{\text{\tiny p}}_{1}) & \text{if } \ell \neq 0  \\
\end{cases} 
\quad \forall t_{[1]} \in \mathcal{B}_1, \forall \ell \preceq 2, \\
& \pi_{(t_{[1]}, e_{(i-1)k})} - \pi_{(t_{[1]}, e_{i\ell})} \geq
\begin{cases}
0, & \text{if } k = 0, \ell = 0\\
c^{\text{\tiny s}}_{(i-1)k} - c^{\text{\tiny p}}_{i-1} , & \text{if } k \neq 0, \ell = 0 \\
(c^{\text{\tiny t}}_{i\ell} - c^{\text{\tiny r}}_{i}) + (c^{\text{\tiny s}}_{(i-1) (k\curlyvee\ell)} - c^{\text{\tiny p}}_{i-1}), & \text{if } \ell \neq 0,  \\
\end{cases} 
\\
& \hspace{43mm}
\forall t_{[1]} \in \mathcal{B}_1, \forall i \in [3,I]_{\mathbb{Z}}, \forall k \preceq i-1, \forall \ell \preceq i, \\
& \pi_{ (t_{[1]}, e_{Ik})}- \pi_{\texttt{T}} \geq
\begin{cases}
(c^{\text{\tiny s}}_{I1} - c^{\text{\tiny p}}_{I}) t_{11} & \text{if } k = 0 \\
c^{\text{\tiny s}}_{Ik} - c^{\text{\tiny p}}_{I} & \text{if } k \neq 0
\end{cases}
\quad \forall t_{[1]} \in \mathcal{B}_1, \forall k \preceq I,
\end{align*}
where dual variables $\pi$ are associated with the (primal) flow balance constraints and all dual constraints are associated with the primal variables $x$. The strong duality holds valid because the longest-path formulation is finitely optimal. The claimed reformulation of (DRNS) then follows from a similar proof as that of Proposition \ref{prop:one-pool-ref}. \qed

\section{Proof of Proposition \ref{prop:opd}} \label{apx-prop:opd}
\proof We linearize the bilinear terms in constraints \eqref{onpd-con-target}--\eqref{onpd-con-3}. 
First, for constraints \eqref{onpd-con-1}--\eqref{onpd-con-3}, we define auxiliary variables
\begin{subequations}
\begin{align}
& \hspace{-0.15cm} \left.
\begin{array}{l}
\zeta^{\text{\tiny e}}_{ij} := \zeta^{\text{\tiny e}}_j a_{ij}, \ \ \eta^{\mbox{\tiny e}}_{ij} := \eta^{\mbox{\tiny e}}_j a_{ij},
\\ [0.1cm]
\zeta^{\text{\tiny x}}_{ijk} := \zeta^{\text{\tiny x}}_{jk} a_{ij}, \ \ \eta^{\mbox{\tiny x}}_{ijk} := \eta^{\mbox{\tiny x}}_{jk} a_{ij}, \ \ \forall k \in [j], \\ [0.1cm] 
\end{array}\right\} \ \ \forall i \in [I+1], \forall j \in [J], \label{apx-opd-note-1} \\
& \hspace{-0.15cm} \left.
\begin{array}{l}
\phi^{\mbox{\tiny x}}_{ijk} := \phi^{\mbox{\tiny x}}_{ik} a_{ij}, \ \ \forall k \in [j], \\[0.1cm]
\lambda^{\text{\tiny a}}_{ij} := \lambda_i a_{ij},
\\ [0.1cm]
v^{\text{\tiny a}}_{ij \ell} := v_{i \ell} a_{ij}, \ \ \nu^{\text{\tiny a}}_{ij \ell} := \nu_{i \ell} a_{ij} , \ \forall \ell \in [y^{\text{\tiny U}}_i - y^{\text{\tiny L}}_i], \\[0.1cm]
\end{array}\right\} \ \ \forall i \in [I+1], \forall j \in [J]. \label{apx-opd-note-1-1}
\end{align}
We equivalently linearize these bilinear equalities as
\begin{align}
& \hspace{-0.15cm} \left.
\begin{array}{l}
\zeta^{\text{\tiny e}}_j = \sum_{i=1}^{I+1} \zeta^{\text{\tiny e}}_{ij}, \ \ 
\eta^{\text{\tiny e}}_j = \sum_{i=1}^{I+1} \eta^{\text{\tiny e}}_{ij},  \\ [0.1cm]
\zeta^{\text{\tiny x}}_{jk} = \sum_{i=1}^{I+1} \zeta^{\text{\tiny x}}_{ijk}, \ \ 
\eta^{\text{\tiny x}}_{jk} = \sum_{i=1}^{I+1} \eta^{\text{\tiny x}}_{ijk}, \ \forall k \in [j], \\[0.1cm]
\end{array}\right\} \ \ \forall j \in [J], \label{apx-opd-note-2} \\
& \hspace{-0.15cm} \left.
\begin{array}{l}
0 \leq \zeta^{\text{\tiny e}}_{ij} \leq K a_{ij}, \ \ -K a_{ij} \leq \eta^{\text{\tiny e}}_{ij} \leq K a_{ij},  \\ [0.1cm]
0 \leq \zeta^{\text{\tiny x}}_{ijk} \leq K a_{ij}, \ \ -K a_{ij} \leq \eta^{\text{\tiny x}}_{ijk} \leq K a_{ij}, \ \forall k \in [j], \\[0.1cm]
\end{array}\right\} \ \ \forall i \in [I+1], \forall j \in [J], \label{apx-opd-note-3} \\
& \hspace{-0.15cm} \left.
\begin{array}{l}
0 \leq \phi^{\text{\tiny x}}_{ijk} \leq K a_{ij}, \ \ \phi^{\text{\tiny x}}_{ik} - K (1-a_{ij}) \leq \phi^{\text{\tiny x}}_{ijk} \leq \phi^{\text{\tiny x}}_{ik} + K (1-a_{ij}), \ \forall k \in [j], \\[0.1cm]
- c^{\text{\tiny x}}_{J} a_{ij} \leq \lambda^{\text{\tiny a}}_{ij} \leq 0, \ \ \lambda_i \leq \lambda^{\text{\tiny a}}_{ij} \leq \lambda_i + c^{\text{\tiny x}}_{J} (1- a_{ij}), \\ [0.1cm]
0 \leq v^{\text{\tiny a}}_{ij \ell} \leq a_{ij}, \ \ v_{i \ell} - (1-a_{ij}) \leq v^{\text{\tiny a}}_{ij \ell} \leq v_{i \ell} + (1-a_{ij}), \ \forall \ell \in [y^{\text{\tiny U}}_i - y^{\text{\tiny L}}_i], \\ [0.1cm]
-c^{\text{\tiny x}}_{J} a_{ij} \leq \nu^{\text{\tiny a}}_{ij \ell} \leq c^{\text{\tiny x}}_{J} a_{ij}, \ \forall \ell \in [y^{\text{\tiny U}}_i - y^{\text{\tiny L}}_i], \\ [0.1cm]
\nu_{i \ell} - c^{\text{\tiny x}}_{J}(1-a_{ij}) \leq \nu^{\text{\tiny a}}_{ij \ell} \leq \nu_{i \ell} + c^{\text{\tiny x}}_{J}(1-a_{ij}), \ \forall \ell \in [y^{\text{\tiny U}}_i - y^{\text{\tiny L}}_i], \\ [0.1cm]
\end{array}\right\} \ \ \forall i \in [I+1], \forall j \in [J], \label{apx-opd-note-3-1} 
\end{align}
To see the equivalence, on the one hand, we notice that constraints \eqref{apx-opd-note-2} follow from \eqref{apx-opd-note-1} and \eqref{onpd-note-1}. Similarly, constraints \eqref{apx-opd-note-3} follow from \eqref{apx-opd-note-1} and the facts that $a_{ij}$ are binary and $\zeta^{\text{\tiny e}}_j, \zeta^{\text{\tiny x}}_{jk} \geq 0$. On the other hand, constraints \eqref{apx-opd-note-2} and \eqref{onpd-note-1} imply that $\zeta^{\text{\tiny e}}_{ij} = \zeta^{\text{\tiny e}}_j$ if $a_{ij} = 1$, and constraints \eqref{apx-opd-note-3} imply that $\zeta^{\text{\tiny e}}_{ij} = 0$ if $a_{ij} = 0$. We hence have $\zeta^{\text{\tiny e}}_{ij} = \zeta^{\text{\tiny e}}_j a_{ij}$. Likewise, we establish $\eta^{\mbox{\tiny e}}_{ij} = \eta^{\mbox{\tiny e}}_j a_{ij}$, $\zeta^{\mbox{\tiny x}}_{ijk} = \zeta^{\mbox{\tiny x}}_{jk} a_{ij}$, and $\eta^{\mbox{\tiny x}}_{ijk} = \eta^{\mbox{\tiny x}}_{jk} a_{ij}$, and hence constraints \eqref{apx-opd-note-1}. 
Moreover, constraints \eqref{apx-opd-note-3-1} are the McCormick inequalities that equivalently express the bilinear terms in  \eqref{apx-opd-note-1-1}.
It follows that constraints \eqref{onpd-con-1}--\eqref{onpd-con-3} can be recast as
\begin{align*}
& \ \theta_i \geq  \phi^{\mbox{\tiny e}}_i + \sum_{j =1}^J (\zeta^{\mbox{\tiny e}}_{ij} + \eta^{\mbox{\tiny e}}_{ij}), \ \forall i \in [I], \\
& \ \theta_i \geq  \sum_{j =1}^J (\zeta^{\mbox{\tiny e}}_{ij} + \eta^{\mbox{\tiny e}}_{ij}) + \Big\{ \phi^{\mbox{\tiny x}}_{ij^*k^*} + (\zeta^{\mbox{\tiny x}}_{ij^* k^*} + \eta^{\mbox{\tiny x}}_{ij^* k^*}) - (\zeta^{\mbox{\tiny e}}_{ij^*} + \eta^{\mbox{\tiny e}}_{ij^*}) \Big\} +  \sum_{j \in [J] \setminus \{j^*\} } \chi_{ijk^*}, \\
& \hspace{9mm} \forall i \in [I], \forall j^* \in [J], \forall k^* \in [j^*], \\
& \ \chi_{ijk^*} \geq (\zeta^{\mbox{\tiny x}}_{ij k} + \eta^{\mbox{\tiny x}}_{ij k}) - (\zeta^{\mbox{\tiny e}}_{ij} + \eta^{\mbox{\tiny e}}_{ij}), \ \forall i \in [I], \forall j \in [J], \forall k^* \in [J], \forall k \in [\min \{j,k^*\}],  \\
& \hspace{-0.15cm} \left.
\begin{array}{l}
\phi^{\mbox{\tiny x}}_{ijk} + ( c^{\mbox{\tiny x}}_k + \lambda^{\text{\tiny a}}_{ij}) y^{\tinyl}_i + \sum_{\ell=1}^{y^{\mbox{\tinyu}}_i - y^{\mbox{\tinyl}}_i}  (c^{\mbox{\tiny x}}_k v^{\text{\tiny a}}_{ij\ell} + \nu^{\text{\tiny a}}_{ij\ell} )  \geq 0,  \\ [0.1cm]
\phi^{\mbox{\tiny x}}_{ijk} \geq 0, 
\end{array}\right\} \ \ \forall i \in [I+1], \forall j \in [J],\forall k\in [j], \\
& \ \eqref{apx-opd-note-2}, \eqref{apx-opd-note-3}, \eqref{apx-opd-note-3-1}.
\end{align*}

Second, we linearize constraint \eqref{onpd-con-target} by claiming that
\begin{align}
& \ \Bigl( \theta_i + c^{\mbox{\tiny y}}_i y^{\tinyl}_i + g_i(y^{\tinyl}_i) \lambda_i + \sum_{\ell=1}^{y^{\mbox{\tinyu}}_i - y^{\mbox{\tinyl}}_i} \bigl(c^{\mbox{\tiny y}}_i v_{i\ell} + \delta_{i\ell} \nu_{i\ell} \bigr) \Bigr) o_i \nonumber \\
= \ & \ \theta_i + c^{\mbox{\tiny y}}_i y^{\tinyl}_i o_i + g_i(y^{\tinyl}_i) \lambda_i + \sum_{\ell=1}^{y^{\mbox{\tinyu}}_i - y^{\mbox{\tinyl}}_i} \bigl(c^{\mbox{\tiny y}}_i v_{i\ell} + \delta_{i\ell} \nu_{i\ell} \bigr), \label{apx-opd-note-7}
\end{align}
which holds valid if $\theta_i = 0$, $\lambda_i = 0$, and $v_{i\ell} = 0$ for all $\ell \in [y^{\mbox{\tinyu}}_i - y^{\mbox{\tinyl}}_i]$ whenever $o_i = 0$ (note that variables $\nu_{i\ell}$ also vanish in this case because $\nu_{i\ell} = \lambda_i v_{i\ell}$). To this end, we incorporate constraints
\begin{align}
-c^{\mbox{\tiny x}}_{J} o_i \leq \lambda_i \leq 0 \label{apx-opd-note-4}
\end{align}
because $\lambda_i \in [-c^{\mbox{\tiny x}}_{J(i)}, 0]$ without loss of optimality by Proposition \ref{prop:big-M}. This guarantees that $o_i = 0$ implies $\lambda_i = 0$. Additionally, we incorporate constraints
\begin{align}
v_{i1} \leq o_i. \label{apx-opd-note-5}
\end{align}
\end{subequations}
Then, $o_i = 0$ implies $v_{i1} = 0$ and hence $v_{i\ell}=0$ for all $\ell \in [y^{\mbox{\tinyu}}_i - y^{\mbox{\tinyl}}_i]$ due to constraints \eqref{ref-note-33}. Furthermore, $o_i = 0$ implies that $a_{ij} = 0$ for all $j \in [J]$ by constraints \eqref{onpd-note-2}. It follows that $\zeta^{\text{\tiny e}}_{ij} = \eta^{\text{\tiny e}}_{ij} = \zeta^{\text{\tiny x}}_{ijk} = \eta^{\text{\tiny x}}_{ijk} = \phi^{\text{\tiny x}}_{ijk} = 0$ for all $j \in [J]$ and $k \in [j]$. 
It remains to ensure that $\phi^{\text{\tiny e}}_i = 0$ whenever $o_i = 0$, which is true due to constraints \eqref{disjoint-ref-1}.

Indeed, if $o_i = 0$ then $\lambda_i = 0$ by constraints \eqref{apx-opd-note-4} and $\nu_{i\ell} = v_{i\ell} = 0$ for all $\ell \in [y^{\mbox{\tinyu}}_i - y^{\mbox{\tinyl}}_i]$ by constraints \eqref{apx-opd-note-5}. Therefore, constraint \eqref{onpd-con-target} is equivalently linearized through equality \eqref{apx-opd-note-7} and incorporating constraints \eqref{apx-opd-note-4}--\eqref{apx-opd-note-5}. This completes the proof. \qed

\section{Symmetry Breaking Inequalities for the (OPD) Model} \label{apx-symmetry-breaking}
We consider two types of symmetry among integer solutions. First, suppose that there are 2 pools and 4 units. The following two unit assignments lead to symmetric integer solutions: (i) assigning all units to pool 1 and no unit to pool 2 (i.e., $a_{1j} = 1$ and $a_{2j} = 0$ for all $j \in [4]$) and (ii) assigning all units to pool 2 and no unit to pool 1 (i.e., $a_{1j} = 0$ and $a_{2j} = 1$ for all $j \in [4]$). We call this ``pool symmetry.'' To break this symmetry, we designate that all open pools have smaller indices than the closed ones. This designation breaks the pool symmetry because the above case (ii) is now prohibited. Accordingly, we add the following inequalities to the (OPD) formulation:
\begin{align*}
& o_i \geq o_{i+1}, \ \ \forall i \in [I-1].
\end{align*}
Second, the following two unit assignments also lead to symmetric integer solutions: (iii) assigning units 1 and 3 to pool 1 and units 2 and 4 to pool 2 (i.e., $a_{11} = 1 - a_{12} = a_{13} = 1 - a_{14} = 1$ and $1 - a_{21} = a_{22} = 1 - a_{23} = a_{24} = 1$) and (iv) assigning units 1 and 3 to pool 2 and units 2 and 4 to pool 1 (i.e., $1 - a_{11} = a_{12} = 1 - a_{13} = a_{14} = 1$ and $a_{21} = 1 - a_{22} = a_{23} = 1 - a_{24} = 1$). We call this ``unit symmetry.'' To break this symmetry, we rank the pools based on the smallest unit index in each pool. That is, we designate that the smallest unit index in pool $i$ is smaller than that in pool $i+1$ for all $i \in [I-1]$, if both pools are opened. This designation breaks the unit symmetry because the above case (iv) is now prohibited. Accordingly, we add the following inequalities to the (OPD) formulation:
\begin{align*}
& \sum_{\ell = 1}^{j-1} a_{i\ell} \geq a_{(i+1)j}, \ \ \forall i \in [I-1], \ \forall j \in [J].
\end{align*}

\begin{table}[H] 
	\centering
	\caption{Nurse demand and attendance entries of the 5-unit system for 4 years.}	
	\label{Table:Data-5unit}
	\begin{tabular}{|c||c|c|c|c|c||c|c|c|c|c|}
	\hline
			& \multicolumn{5}{c||}{Demand} & \multicolumn{5}{c|}{Attendance}  \\ \hline
	Entries & $\widehat{d}_1$ & $\widehat{d}_2$ & $\widehat{d}_3$ & $\widehat{d}_4$ & $\widehat{d}_5$ & $(w_1,\widehat{w}_1)$ & $(w_2,\widehat{w}_2)$ & $(w_3,\widehat{w}_3)$ & $(w_4,\widehat{w}_4)$ & $(w_5,\widehat{w}_5)$    \\ \hline
	$1$  & 12&	14&	18&	13&	16&  (11,8) & (12,9) & (16,12) & (11,8) & (14,11) \\ \hline
	$2$ & 12&	14&	18&	14&	15& (10,8) & (11,9) & (15,12) & (10,8) & (13,10) \\ \hline
	$3$ & 12&	13&	17&	12&	14& (11,7) & (11,9) & (15,12) & (10,8) & (14,11) \\ \hline
	$\vdots$ & $\vdots$ & $\vdots$ & $\vdots$ & $\vdots$  & $\vdots$ & $\vdots$ & $\vdots$ & $\vdots$ & $\vdots$ & $\vdots$  \\ \hline
	\end{tabular}
\end{table}

\begin{figure}[H]
	\centering	
\begin{subfigure}[b]{0.38\textwidth}
	\includegraphics[width=\textwidth]{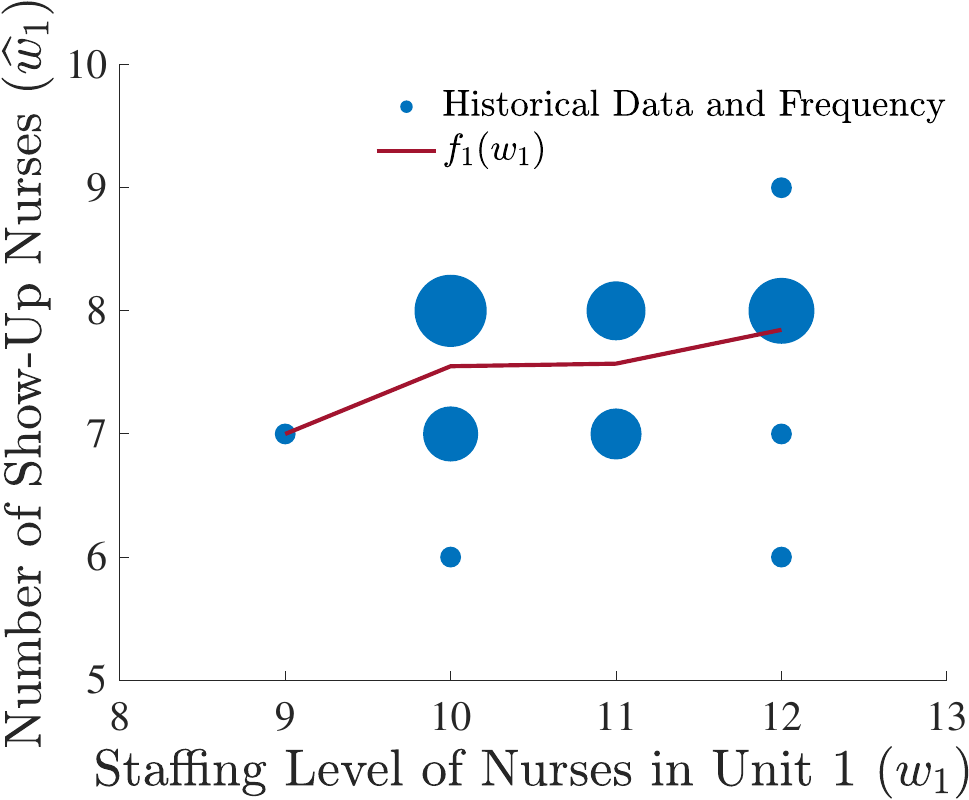}
	\includegraphics[width=\textwidth]{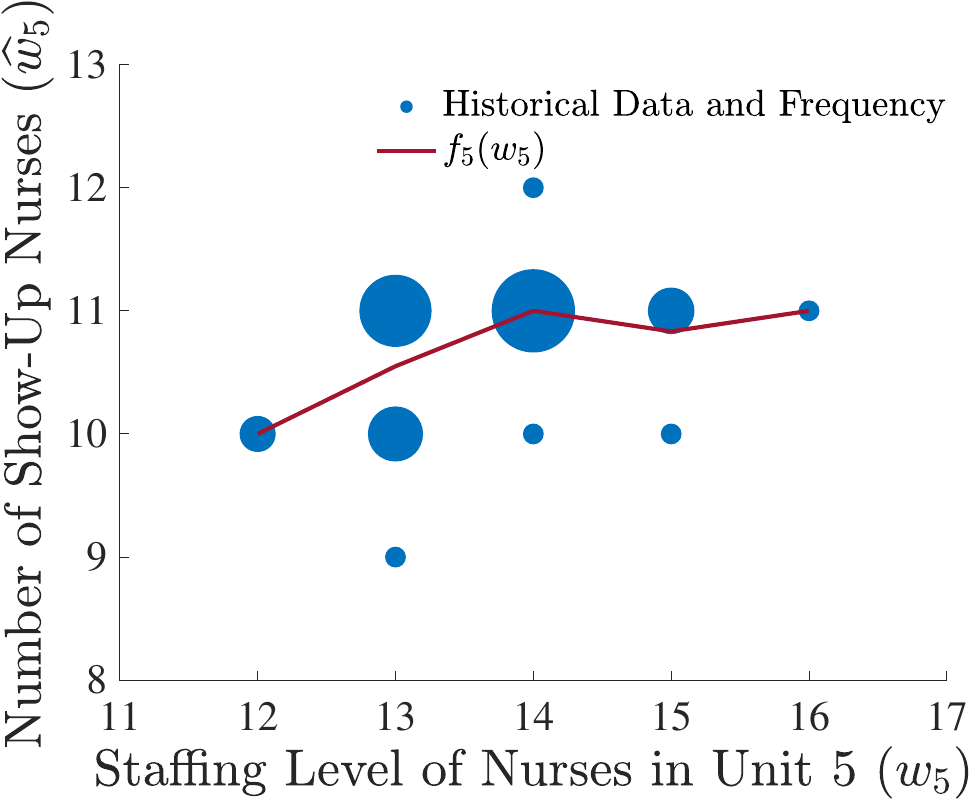}		
	\caption{}
	\label{Fig:Showupnurses}
\end{subfigure} \hspace{-0.0cm}
\begin{subfigure}[b]{0.4\textwidth}
	\includegraphics[width=\textwidth]{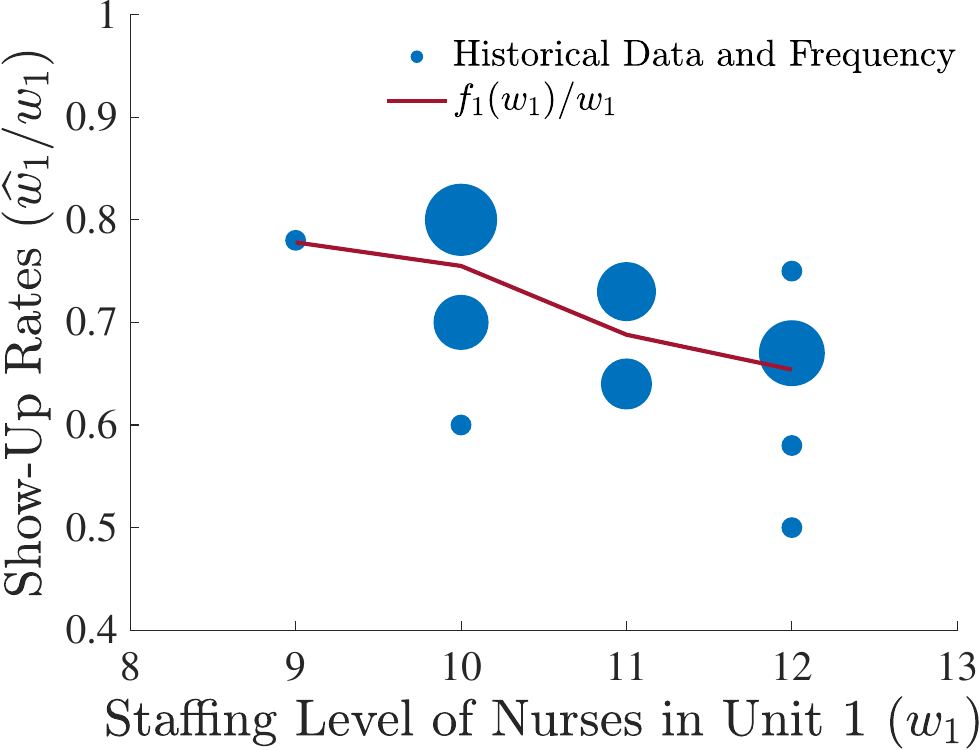}
	\includegraphics[width=\textwidth]{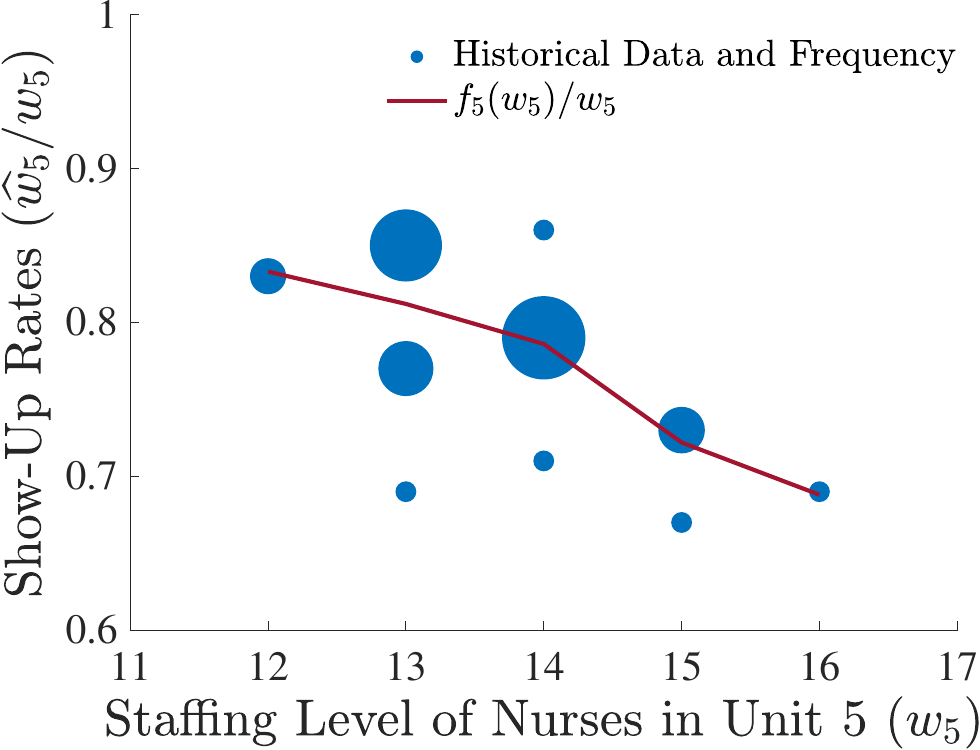}
	\caption{}	
	\label{Fig:Showuprates}
\end{subfigure}	
\caption{The number of show-up nurses (Figure \eqref{Fig:Showupnurses}) and the show-up rates (Figure \eqref{Fig:Showuprates}) with respect to the staffing level of nurses in units CVC5 and 7B.}
\label{Fig:Attendance}
\end{figure}

\section{A 5-Unit System} \label{apx-5-unit}
{\color{black}
Our test instances are based on 5 units, abbreviated as CVC5, 4B, 4C, 7A1, and 7B in our collaborating hospital. We obtain the nurse demand and attendance entries of these 5 units during 4 years from July-2010 to June-2014. We display the format of this data in Table \ref{Table:Data-5unit}. For example, the first entry states that, in the first day, $w_1=11$ unit nurses were assigned to CVC5 and $\widehat{w}_1=8$ of them showed up for work, while the nurse demand was $\widehat{d}_1=12$ for the unit. For this particular day, the nurse show-up rate is then $\widehat{w}_1/w_1 = 72.7\%$. In Figure \ref{Fig:Attendance}, we depict the number of show-up nurses, as well as the show-up rates, as functions of staffing levels in CVC5 and 7B based on all data entries. From this figure, we observe that not only the number of show-up nurses, but also the show-up rates, depend on the corresponding staffing levels.

\begin{table}[H]
	\centering
	\captionsetup{width=15cm}
	\caption{Parameters of $\mathcal{D}(w, y)$ in the 5-unit system.}
	\label{Table:5Unit-Parameters}
	\begin{tabular}{|c||c|c|c|c||c|c|l|}
		\hline
		& \multicolumn{4}{c||}{Nurse demand} & \multicolumn{3}{c|}{Endogenous nurse show-up} \\ \hline
		Unit $j$ & $d^{\text{\tiny L}}_j$ & $d^{\text{\tiny U}}_j$ & $\mu_{j1}$ & $\text{sd}_{j}$ & $w^{\text{\tiny L}}_j$ &  $w^{\text{\tiny U}}_j$ & \multicolumn{1}{c|}{$\{(w_j, f_j(w_j)): w_j \in [w^{\text{\tiny L}}_j, w^{\text{\tiny U}}_j]_{\mathbb{Z}} \}$}  \\ \hline
		$1$ & $4$ & $12$ & $10.39$ & $1.75$ & $9$ & $12$  & $\{ (9, 7),  (10, 7.55),  (11, 7.57),  (12, 7.85) \}$ \\
		$2$ & $6$ & $14$ & $12.91$ & $1.34$ & $10$ & $13$ & $\{ (10, 8.5),  (11, 8.89), (12, 9),  (13, 9)  \}$ \\	
		$3$ & $8$ & $18$ & $16.60$ & $1.55$ & $14$ & $17$ & $\{ (14,11.64),  (15, 11.68),  (16,11.86),  (17, 12) \}$ \\			
		$4$ & $5$ & $14$ & $12.08$ & $1.48$ & $8$ & $13$  & $\{ (8,7),  (9,7.5),  (10, 8),  (11,8.75),  (12, 8.82),  (13,9) \}$ \\	
		$5$ & $6$ & $16$ & $14.58$ & $1.42$ & $12$ & $16$ & $\{ (12,10),  (13, 10.55),  (14,11), (15, 10.83),  (16,11) \}$ \\	\hline
	\end{tabular}
\end{table}

To calibrate the ambiguity set $\mathcal{D}(w, y)$, we randomly split the data into a training data set (taking up 80\% of all data entries) and a testing data set (20\%). We set $Q := 2$, i.e., we consider mean value $\mu_{j1}$ and standard deviation $\text{sd}_j$ of the nurse demand, and we set $(d^{\text{\tiny L}}_j, d^{\text{\tiny U}}_j, \mu_{j1}, \text{sd}_j)$ as the empirical min, max, mean, and standard deviation of the training data, respectively. Likewise, we set $(w^{\text{\tiny L}}_j, w^{\text{\tiny U}}_j, f_j(w_j))$ as the empirical min, max, and mean of the training data. For pool nurses, we assume a constant absence rate such that, for all $i \in [I]$, $g_i(y_i) := (1-A^{\mbox{\tiny p}}_i) y_i$, where $A^{\mbox{\tiny p}}_i$ is randomly extracted from the interval $[0, 0.02]$. Additionally, we set $y^{\tinyl}_i := 0$ and $y^{\tinyu}_i := \sum_{j \in P_i} (d^{\text{\tiny U}}_j - w^{\text{\tiny L}}_j)$. We summarize the parameters of $\mathcal{D}(w, y)$ in Table \ref{Table:5Unit-Parameters}.

Given an annual salary ranging from $\$60,000$ to $\$100,000$, we set the cost of staffing a unit nurse during a specific shift as $c^{\text{\tiny w}}_j = \$250$ for all $j \in [J]$ and that of staffing a pool nurse during the shift as $c^{\text{\tiny y}}_i = \$416$ for all $i \in [I]$. In addition, we assume that staffing temp nurses costs $c^{\text{\tiny x}}_j= \$1000 + \tilde{\alpha}_j$ for all $j \in [J]$, where $\tilde{\alpha}_j$ is uniformly extracted from the interval $[0, 100]$.

To evaluate given staffing levels, denoted by $(w^*, y^*)$, we conduct an out-of-sample simulation based on the testing data. Specifically, we compute an average staffing cost $\sum_{j=1}^J c^{\mbox{\tiny w}}_j w^*_j + \sum_{i=1}^I c^{\mbox{\tiny y}}_i y^*_i + \mathbb{E}_{\mathbb{P}(w^*)}\big[ V(\tilde{w}, y^*, \tilde{d}) \big]$, where $\mathbb{P}(w^*)$ represents the product measure of the empirical distributions of $\tilde{w}_j$ (with regard to the fixed staffing level $w^*_j$) and the empirical distribution of $\tilde{d}$, all based on the testing data.
}

\section{Computational Efficacy of Sep and Sep$_{\text{\tiny VI}}$} \label{apx-seperation}
{\color{black} We compare the computational efficacy of Sep, Sep$_{\text{\tiny VI}}$, and directly solving the exponentially-sized MILP reformulation \eqref{drns-milp-exp} (denoted by MILP$_{\text{\tiny Exp}}$). To this end, we construct $J$-unit systems, where $J \in [5, 9]_{\mathbb{Z}}$, by duplicating units from the 5-unit system introduced in Appendix \ref{apx-5-unit}. In addition, for fixed $J$, we construct a general pool with $P_1 = \{ j \in [J] : j \text{ is even} \}$, $P_2 = \{ j \in [J] : j \text{ is odd} \}$, $P_3 = \{ j \in [J] : j \leq J/2 \}$, and $P_4 = \{ j \in [J] : j > J/2 \}$. In Table \ref{Table:Arbitrary}, we report the time for constructing each test instance in GUROBI and that for solving the instance.
\begin{table}[H]	
	\centering
	\caption{Wall-clock seconds for solving (DRNS) using various approaches.}
	\label{Table:Arbitrary}
	\begin{tabular}{|c||r|r|r|}	
	\hline
	 \multicolumn{4}{|c|}{Total time (construction time, solution time)} \\ \hline
	$J$ &  \multicolumn{1}{c|}{MILP$_{\text{\tiny Exp}}$} &  \multicolumn{1}{c|}{Sep}	& \multicolumn{1}{c|}{Sep$_{\text{\tiny VI}}$} \\ \hline  
	$5$ & $0.39 \ (0.13, 0.26)$ & $0.24 \ (0.01,0.23)$ & $0.17 \ (0.01, 0.16)$  \\
	$6$ & $3.13 \ (1.02, 2.11)$ & $0.64 \ (0.02, 0.62)$ & $0.25 \ (0.02, 0.23)$  \\		
	$7$ & $30.35 \ (10.75, 19.60)$ & $0.79 \ (0.02, 0.77)$ & $0.28 \ (0.02, 0.26)$ \\
	$8$ & $377.84 \ (121.40, 256.44)$ & $1.77 \ (0.02, 1.75)$ & $0.40 \ (0.02, 0.38)$ \\ 
	$9$ & $6330.28 \ (1649.20, 4681.08)$ & $1.62 \ (0.03, 1.59)$ & $0.44 \ (0.03, 0.41)$ \\ \hline	
	\end{tabular}	
\end{table}
From the above table, we observe that the time it takes to construct a test instance in GUROBI can be non-trivial. For example, the construction time for formulation \eqref{drns-milp-exp} increases drastically as $J$ increases, mainly because of the exponentially-sized constraints \eqref{ref-linear-3}. In fact, the construction time exceeds 1 hour once $J > 9$, which is why we stop this comparison at $J = 9$. It is worth mentioning that this observation remains when using other commercial solvers (e.g., CPLEX). This demonstrates the computational efficacy of Sep and Sep$_{\text{\tiny VI}}$, which incorporate constraints \eqref{ref-linear-3} on-the-fly. In addition, the solution times for Sep and Sep$_{\text{\tiny VI}}$ are significantly more scalable than that of MILP$_{\text{\tiny Exp}}$.}

\section{Stochastic Programming Formulation} \label{apx-SPformulation}
In Section~\ref{sec:results_SP_DRO}, we compare DRNS with the following stochastic programming (SP) model:
\begin{subequations}
	\label{SP-formulation1}
	\begin{align}
	\min_{w,y} \ & \sum_{j=1}^J c^{\text{\tiny w}}_j w_j + \sum_{i=1}^I c^{\text{\tiny y}}_i y_i + \mathbb{E}_{\mathbb{P}} \Big[ V(w, y, \tilde{u}, \tilde{p}, \tilde{d} ) \Big] \\
	\mbox{s.t.} \
	& \mbox{\eqref{sns-con-w-bd}--\eqref{sns-con-integer}}, \nonumber
	\end{align}
	where ($\tilde{u}, \tilde{p}, \tilde{d}$) represent the random show-up rate of unit nurses, show-up rate of pool nurses, and nurse demand, respectively. $V(w, y, \tilde{u}, \tilde{p}, \tilde{d})$ equals the optimal value of the following linear program:
	\begin{align}
	V(w, y, \tilde{u}, \tilde{p}, \tilde{d} ) := 
	\max_{(\alpha, \beta) \in \Lambda} \ & \Big\{\sum_{j=1}^J (\tilde{d}_j -  \tilde{u}_j w_j  ) \alpha_j + \sum_{i=1}^I (  \tilde{p}_i y_i ) \beta_i \Big\},
	\end{align}
\end{subequations}
where $\Lambda = \{(\alpha, \beta): \mbox{\eqref{sns-sub-dual-con-max}--\eqref{sns-sub-dual-con-bd}}\}$. We remark that, unlike DRNS, this SP model assumes \emph{exogenous} uncertainties, i.e., the distribution of $\tilde{u}, \tilde{p}$ does not depend on the staffing levels $(w, y)$. To solve this SP model, we adopt the sample average approximation method using the data described in Appendix~\ref{apx-5-unit}.

\section{Capacitated Temporary Nurses} \label{apx-capacitated_temp_nurses}
In this section, we extend the (DRNS) model with a finite supply of the temporary nurses. In case of nurse shortage, we assume that there is a capacity $x^{\text{\tiny U}}_j$ of temporary nurses in unit $j$ for all $j \in [J]$, and if there are remaining shortage after hiring all available nurses, a unit penalty cost $c^{\text{\tiny u}}$ will be applied. This leads to a revised formulation \eqref{sns-sub} as follows:
\begin{subequations}
	\begin{align}
	V(\tilde{w}, \tilde{y}, \tilde{d}) \ = \ \min_{z, x, u, e} \ & \ \sum_{j=1}^J  c^{\mbox{\tiny x}}_j x_j + c^{\text{\tiny u}} u_j \label{sns-sub-obj-u} \\
	\mbox{s.t.} \ & \ \sum_{i \in [I]:\ j \in P_i} z_{ij} + x_j + u_j - e_j = \tilde{d}_j - \tilde{w}_j, \ \ \forall j \in [J], \label{sns-sub-con-demand-u} \\
	& \ \sum_{j \in P_i} z_{ij} \leq \tilde{y}_i, \ \ \forall i \in [I], \label{sns-sub-con-pool-u} \\
	& \ x_j \leq x^{\text{\tiny U}}_j, \ \forall j \in [J], \label{sns-sub-con-xu} \\
	& \ x_j, u_j, e_j \in \mathbb{Z}_+, \ \ \forall j \in [J], \ \ z_{ij} \in \mathbb{Z}_+, \ \ \forall i \in [I], \ \forall j \in P_i, \label{sns-sub-con-integer-u}
	\end{align}
\end{subequations}
where variables $u_j$ represent the number of nurse shortage in unit $j$. In what follows, we show that the main theoretical results in Sections \ref{sec: Solution Approach} and \ref{sec:tractable} can be extended in this more general setting. We omit their proofs due to the similarity.

First, similar to Lemma \ref{prop:sns-sub-tu}, for any given $(\tilde{w}, \tilde{y}, \tilde{d}) \in \Xi$, the value of $V(\tilde{w}, \tilde{y}, \tilde{d})$ remains unchanged if constraints \eqref{sns-sub-con-integer-u} are replaced by non-negativity restrictions.
Thus, we rewrite $V(\tilde{w}, \tilde{y}, \tilde{d})$ as the following dual formulation:
\begin{subequations}
	\begin{align}
	V(\tilde{w}, \tilde{y}, \tilde{d}) \ = \ \max_{\alpha, \beta, \psi} \ & \ \sum_{j=1}^J \Big\{ (\tilde{d}_j - \tilde{w}_j) \alpha_j - x^{\text{\tiny U}}_j \psi_j \Big\} + \sum_{i=1}^I \tilde{y}_i \beta_i \label{sns-sub-dual-obj-u} \\
	\mbox{s.t.} \ & \ \beta_i + \alpha_j \leq 0, \ \ \forall i \in [I], \ \forall j \in P_i, \label{sns-sub-dual-con-max-u} \\
	& \ \alpha_j - \psi_j \leq c^{\mbox{\tiny x}}_j, \ \ \forall j \in [J], \label{sns-sub-dual-con-x} \\
	& \ 0 \leq \alpha_j \leq c^{\mbox{\tiny x}}_{J+1}, \ \ \forall j \in [J], \label{sns-sub-dual-con-bd-u} \\
	& \ \psi_j \geq 0, \ \forall j \in [J], \label{sns-sub-dual-con-bd2}
	\end{align}
\end{subequations}
where we define $c^{\text{\tiny x}}_{J+1} := c^{\text{\tiny u}}$ in constraints \eqref{sns-sub-dual-con-bd-u} for ease of exposition. 
Dual variables $\alpha_j$, $\beta_i$, and $\psi_j$ are associated with primal constraints \eqref{sns-sub-con-demand-u}, \eqref{sns-sub-con-pool-u}, and \eqref{sns-sub-con-xu}, respectively, and dual constraints \eqref{sns-sub-dual-con-max-u}, \eqref{sns-sub-dual-con-x}, and \eqref{sns-sub-dual-con-bd-u} are associated with primal variables $z_{ij}$, $x_j$, and $(u_j, e_j)$, respectively. 
Strong duality between formulations \eqref{sns-sub-obj-u}--\eqref{sns-sub-con-integer-u} and \eqref{sns-sub-dual-obj-u}--\eqref{sns-sub-dual-con-bd2} holds because \eqref{sns-sub-obj-u}--\eqref{sns-sub-con-integer-u} has a finite optimal value.
It follows that $\psi_j = \big[ \alpha_j - c^{\text{\tiny x}}_j \big]_+$ because $x^{\text{\tiny U}}_j \geq 0$. Then,
\begin{align}
V(\tilde{w}, \tilde{y}, \tilde{d}) \ = \ \max_{(\alpha, \beta, \psi) \in \Lambda'} \ & \ \Bigg\{\sum_{j=1}^J \Big\{ (\tilde{d}_j - \tilde{w}_j) \alpha_j - x^{\text{\tiny U}}_j \big[ \alpha_j - c^{\text{\tiny x}}_j \big]_+ \Big\} + \sum_{i=1}^I \tilde{y}_i \beta_i \Bigg\} \label{sns-sub-dual1-obj},
\end{align}
where $\Lambda' := \{(\alpha, \beta, \psi): \mbox{\eqref{sns-sub-dual-con-max-u}--\eqref{sns-sub-dual-con-bd2}}\}$ represents the dual feasible region for variables $(\alpha, \beta, \psi)$.

Second, we have the following extension of Proposition \ref{prop:ref}.\\
\textbf{Proposition \ref{prop:ref}$^\prime$}
Under Assumption \ref{assump:technical}, the (DRNS) model \eqref{drns} yields the same optimal value and the same set of optimal solutions as the following min-max optimization problem:
\begin{subequations}
	\begin{align}
	\min_{\substack{w, y\\ \gamma, \lambda, \rho}} \ & \ \max_{(\alpha, \beta, \psi) \in \Lambda'} F(\alpha, \beta) + \sum_{i=1}^I (c^{\mbox{\tiny y}}_i y_i + g_i(y_i) \lambda_i) + \sum_{j=1}^J \Bigg[c^{\mbox{\tiny w}}_j w_j + \sum_{q=1}^Q \mu_{jq} \rho_{jq} + f_j(w_j) \gamma_j \Bigg] \label{ref-note-8-u} \\
	\mbox{s.t.} \ & \ \mbox{\eqref{sns-con-w-bd}--\eqref{sns-con-integer}}, \label{ref-note-con-u}
	\end{align}
	where
	\begin{equation}
	F(\alpha, \beta) \ := \ \sum_{j=1}^J \Bigg\{ -x^{\text{\tiny U}}_j \big[ \alpha_j - c^{\text{\tiny x}}_j \big]_+  + \big[(- \alpha_j - \gamma_j) w_j \big]_+ + \sup_{\tilde{d}_j \in [d^{\tinyl}_j, d^{\tinyu}_j]_{\mathbb{Z}}} \Bigl\{\alpha_j \tilde{d}_j - \sum_{q=1}^Q\rho_{jq}\tilde{d}^q_j \Bigr\} \Bigg\} + \sum_{i=1}^I \big[ (\beta_i - \lambda_i) y_i \big]_+. \label{ref-f-def-u}
	\end{equation}
\end{subequations}

An extension of Lemma \ref{thm:extreme-points} follows.\\
\textbf{Lemma \ref{thm:extreme-points}$^\prime$}
Without loss of generality, we assume that $0 \leq c^{\text{\tiny x}}_1 \leq \ldots \leq c^{\text{\tiny x}}_{J+1} \equiv c^{\text{\tiny u}}$.	
For fixed $(w, y, \gamma, \lambda, \rho)$, there exists an optimal solution $(\bar{\alpha}, \bar{\beta}, \bar{\psi})$ to problem $\max_{(\alpha, \beta, \psi) \in \Lambda'} F(\alpha, \beta)$ such that 
(a) $\bar{\alpha}_j \in \{0, c^{\mbox{\tiny x}}_1, \ldots, c^{\mbox{\tiny x}}_{J+1} \}$ for all $j \in [J]$ and
(b) $\bar{\beta}_i = - \max\{\alpha_j: j \in P_i\}$ for all $i \in [I]$. \\

\indent Using Lemma \ref{thm:extreme-points}$^\prime$, we introduce binary variables to encode the $(\bar{\alpha}, \bar{\beta})$ values identified in the optimality conditions. Specifically, for all $j \in [J]$ and $k \in [J+1]$, binary variable $t_{jk} = 1$ if $\bar{\alpha}_j = c^{\mbox{\tiny x}}_k$ and $t_{jk} = 0$ otherwise.
For all $i \in [I]$ and $k \in [J+1]$, binary variable $s_{ik} = 1$ if $\bar{\beta}_i = -c^{\mbox{\tiny x}}_k$ 
and $s_{ik} = 0$ otherwise. In other words, $s_{ik} = 1$ if $k$ is the largest index in $[J+1]$ such that $\sum_{j \in P_i} t_{jk} > 0$.
Variables $(t,s)$ need to satisfy the following constraints to make the encoding well-defined:
\begin{subequations}
	\begin{align}
	& \sum_{k = 1}^{J+1} t_{jk} \leq 1, \ \forall j \in [J], \label{encode-0-u} \\
	& \sum_{k=1}^{J+1} s_{ik} \leq 1, \ \forall i \in [I], \label{encode-1-u}\\
	& s_{ik} \leq \sum_{j \in P_i} t_{jk}, \ \forall i \in [I], \forall k \in [J+1], \label{encode-2-u} \\
	& s_{ik} + t_{j \ell} \leq 1, \ \forall i \in [I], \forall k  \in [J], \forall j \in P_i, \forall \ell \in [k+1,J+1]_{\mathbb{Z}}, \label{encode-3-u} \\
	& t_{jk} \leq \sum_{\ell = 1}^{J+1} s_{i \ell}, \ \forall i \in [I], \forall j \in P_i, \forall k \in [J+1], \label{encode-4-u}\\
	& s_{ik} \in \mathbb{B}, \ \forall i \in [I], \forall k \in [J+1], \ \ t_{jk} \in \mathbb{B}, \ \forall j \in [J], \forall k \in [J+1].	\label{encode-5-u} 
	\end{align}	
\end{subequations}

The binary encoding extends Theorem \ref{thm:ip} as follows.\\
\textbf{Theorem \ref{thm:ip}$^\prime$} For fixed $(w, y, \gamma, \lambda, \rho)$, problem $\max_{(\alpha, \beta, \psi) \in \Lambda'} F(\alpha, \beta)$ yields the same optimal value as the following integer linear program:
\begin{subequations}
	\label{IntegerProgram-u}
	\begin{align}
	\max_{t,s,r,p} \ & \sum_{j=1}^J \bigg( c ^{\text{\tiny r}}_j r_j + \sum_{k=1}^{J+1} c^{\text{\tiny t}}_{jk} t_{jk} \bigg) + \sum_{i=1}^I \bigg( c^{\text{\tiny p}}_i p_i + \sum_{k=1}^{J+1} c^{\text{\tiny s}}_{ik} s_{ik} \bigg) \label{IntegerProgram-0-u} \\
	\mbox{s.t.} \
	& (t,s,r,p) \in \mathcal{H}' : = \Big\{  \eqref{encode-0}-\eqref{encode-5}, \label{IntegerProgram-1-u} \\ 
	& \hspace{33mm}  r_j + \sum_{k=1}^{J+1} t_{jk} = 1, \ \forall j \in [J],  \label{IntegerProgram-2-u} \\
	& \hspace{33mm}  p_i + \sum_{k=1}^{J+1} s_{ik} = 1, \ \forall i \in [I] \label{IntegerProgram-3-u}
	\Big\},
	\end{align}
	where
	$c^{\text{\tiny r}}_j := [ - \gamma_j w_j ]_+ + \sup_{\tilde{d}_j \in [d^{\tinyl}_j, d^{\tinyu}_j]_{\mathbb{Z}}} \bigg\{ - \sum_{q=1}^Q \rho_{jq} \tilde{d}^q_j \bigg\}$, 
	$c^{\text{\tiny t}}_{jk} := -x^{\text{\tiny U}}_j [c^{\text{\tiny x}}_k - c^{\text{\tiny x}}_j ]_+ + [ (-c^{\text{\tiny x}}_k - \gamma_j ) w_j ]_+ + \sup_{\tilde{d}_j \in [d^{\tinyl}_j, d^{\tinyu}_j]_{\mathbb{Z}}} \bigg\{ c^{\text{\tiny x}}_k \tilde{d}_j - \sum_{q=1}^Q \rho_{jq} \tilde{d}^q_j \bigg\}$, 
	$c^{\text{\tiny p}}_i := [ - \lambda_i y_i  ]_+$, and
	$c^{\text{\tiny s}}_{ik} := [ (-c^{\text{\tiny x}}_k - \lambda_i ) y_i  ]_+$. 
\end{subequations}
\\

This leads to the following MILP reformulation of (DRNS) in this more general setting, which extends Theorem \ref{thm:ref-linearized}. An immediate extension of Algorithm \ref{algo-sep} can be applied to solve this reformulation.\\
\textbf{Theorem \ref{thm:ref-linearized}$^\prime$}
Under Assumption \ref{assump:technical}, the (DRNS) model \eqref{drns} yields the same optimal value as the following mixed-integer program:
\begin{subequations}
	\begin{align}
	\min_{\substack{u, v, \varphi, \nu\\ \gamma, \lambda, \rho, \theta}} \ & \ \theta + \sum_{i=1}^I \Bigl(c^{\mbox{\tiny y}}_i y^{\tinyl}_i + g_i(y^{\tinyl}_i) \lambda_i + \sum_{\ell=1}^{y^{\mbox{\tinyu}}_i - y^{\mbox{\tinyl}}_i} \bigl(c^{\mbox{\tiny y}}_i v_{i\ell} + \delta_{i\ell} \nu_{i\ell}  \bigr) \Bigr) \nonumber \\
	& \ + \sum_{j=1}^J \Biggl[ \sum_{q=1}^Q \mu_{jq} \rho_{jq} + c^{\mbox{\tiny w}}_j w^{\tinyl}_j + f_j(w^{\tinyl}_j) \gamma_j + \sum_{\ell=1}^{w^{\mbox{\tinyu}}_j - w^{\mbox{\tinyl}}_j} \bigl( c^{\mbox{\tiny w}}_j u_{j \ell} + \Delta_{j \ell} \varphi_{j \ell} \bigr) \Biggr] \label{ref-linear-1-u} \\
	\mbox{s.t.} \ & \ \mbox{\eqref{ref-note-36}--\eqref{ref-note-33}}, \label{ref-linear-2-u} \\
	& \ \theta \geq \sum_{j=1}^J \bigg( c ^{\text{\tiny r}}_j r_j + \sum_{k =1}^{J+1} c^{\text{\tiny t}}_{jk} t_{jk} \bigg) + \sum_{i=1}^I \bigg( c^{\text{\tiny p}}_i p_i + \sum_{k = 1}^{J+1} c^{\text{\tiny s}}_{ik} s_{ik} \bigg), \ \ \forall (t, s, r, p) \in \mathcal{H}', \label{ref-linear-3-u} \\
	& \ u_{j \ell} \in \mathbb{B}, \ \ \forall j \in [J], \ \forall \ell \in [w^{\mbox{\tinyu}}_j - w^{\mbox{\tinyl}}_j], \ \ v_{i\ell} \in \mathbb{B}, \ \ \forall i \in [I], \ \forall \ell \in [y^{\mbox{\tinyu}}_i - y^{\mbox{\tinyl}}_i], \label{ref-linear-4-u}
	\end{align}
	where set $\mathcal{H}'$ is defined in \eqref{IntegerProgram-1-u}--\eqref{IntegerProgram-3-u} and coefficients $c^{\mbox{\tiny r}}_j$, $c^{\mbox{\tiny t}}_{jk}$, $c^{\mbox{\tiny p}}_i$, and $c^{\mbox{\tiny s}}_{ik}$ are represented through
	\begin{align}
	c^{\text{\tiny r}}_j &= \Bigl[- \gamma_j w^{\tinyl}_j - \sum_{\ell=1}^{w^{\mbox{\tinyu}}_j - w^{\mbox{\tinyl}}_j} \varphi_{j \ell} \Bigr]_+ + \sup_{\tilde{d}_j \in [d^{\tinyl}_j, d^{\tinyu}_j]_{\mathbb{Z}}} \bigg\{ - \sum_{q=1}^Q \rho_{jq} \tilde{d}^q_j \bigg\}, \label{ref-linear-5-u} \\ 
	c^{\text{\tiny t}}_{jk} &= -x^{\text{\tiny U}}_j \big[ c^{\text{\tiny x}}_k - c^{\text{\tiny x}}_j \big]_+ + \Bigl[(- c^{\mbox{\tiny x}}_k - \gamma_j) w^{\tinyl}_j - \sum_{\ell=1}^{w^{\mbox{\tinyu}}_j - w^{\mbox{\tinyl}}_j} (c^{\mbox{\tiny x}}_k u_{j \ell} + \varphi_{j \ell}) \Bigr]_+ + \sup_{\tilde{d}_j \in [d^{\tinyl}_j, d^{\tinyu}_j]_{\mathbb{Z}}} \bigg\{ c^{\text{\tiny x}}_k \tilde{d}_j - \sum_{q=1}^Q \rho_{jq} \tilde{d}^q_j \bigg\},\label{ref-linear-6-u} \\
	c^{\text{\tiny p}}_i &= \Bigl[ - \lambda_i y^{\tinyl}_i - \sum_{\ell=1}^{y^{\mbox{\tinyu}}_i - y^{\mbox{\tinyl}}_i} \nu_{i\ell} \Bigr]_+, \label{ref-linear-7-u}\\
	c^{\text{\tiny s}}_{ik} &= \Bigl[(- c^{\mbox{\tiny x}}_k - \lambda_i) y^{\tinyl}_i - \sum_{\ell=1}^{y^{\mbox{\tinyu}}_i - y^{\mbox{\tinyl}}_i} (c^{\mbox{\tiny x}}_k v_{i\ell} + \nu_{i\ell}) \Bigr]_+. \label{ref-linear-8-u}
	\end{align}
\end{subequations}
The valid inequalities \eqref{VI-general} are extended as follows.\\
\textbf{Lemma \ref{lem:valid}$^\prime$}
Under any nurse pool structure, inequalities \eqref{VI-general-u} are valid for all $(t, s, r, p) \in \mathcal{H}'$:
\begin{align}
\sum_{\ell = k }^{J+1} t_{j \ell} \leq \sum_{\ell = k}^{J+1} s_{i \ell}, \ \forall i \in [I], \forall j \in P_i, \forall k \in [J+1].
\label{VI-general-u}
\end{align}

Third, we extend the results in Section \ref{sec:tractable} under practical nurse pool structures. Under Structure~\ref{assump:one-pool}, we follow a similar proof to that of Proposition~\ref{prop:onepool-convex} to show that the valid inequalities \eqref{VI-general-u} suffice to describe an ideal formulation of $\mathcal{H}'$, i.e., $\overline{\mathcal{H}'} = \text{conv} (\mathcal{H}')$, where
\begin{subequations}
	\begin{align}
	\overline{\mathcal{H}'} := \bigg\{ (t,s,r,p) \geq 0 : \
	& r_j + \sum_{k =1}^{J+1} t_{jk} = 1, \ \forall j \in [J], \ \ p + \sum_{k =1}^{J+1} s_{k} = 1, \label{H-onepool-0-u} \\	
	& \sum_{k=1}^{J+1} s_{k} \leq 1, \ \sum_{k = 1}^{J+1} t_{jk} \leq 1, \ \forall j \in [J], \label{H-onepool-1-u} \\ 
	& s_{k} \leq \sum_{j = 1}^J t_{jk}, \ \forall k \in [J+1], \label{H-onepool-3-u}\\
	& \sum_{\ell = k }^{J+1} t_{j \ell} \leq \sum_{\ell = k}^{J+1} s_{\ell}, \ \forall j \in [J], \forall k \in [J+1] \bigg\}. \label{onepool-VI-u}
	\end{align}
\end{subequations}

An extension of Theorem \ref{thm:integrality-one} follows.\\
\noindent\textbf{Theorem \ref{thm:integrality-one}$^\prime$}
Under Structure \ref{assump:one-pool},
for fixed $(u, v, \gamma, \lambda, \rho)$, problem $\max_{(\alpha, \beta, \psi) \in \Lambda'} F(\alpha, \beta)$ yields the same optimal value as
\begin{align*}
\max \Bigg\{ c^{\text{\tiny p}}_1 + \sum_{j=1}^J c^{\text{\tiny r}}_j, \ \ \max_{\substack{j^* \in [J], \\ k^* \in [J+1] }} \bigg\{ \sum_{j=1}^J c^{\text{\tiny r}}_j + c^{\text{\tiny s}}_{1k^*}  +  (c^{\text{\tiny t}}_{j^* k^*} - c^{\text{\tiny r}}_{j^*}) + \sum_{j \in [J] \setminus \{j^*\} } \bigg( \max_{k \in [k^*] } \ \Big\{ (c^{\text{\tiny t}}_{jk} - c^{\text{\tiny r}}_j) \Big\} \bigg) \bigg\},
\Bigg\}
\end{align*}
where $c^{\mbox{\tiny r}}_j$, $c^{\mbox{\tiny t}}_{jk}$, $c^{\mbox{\tiny p}}_1$, and $c^{\mbox{\tiny s}}_{1k}$ are computed by \eqref{ref-linear-5-u}--\eqref{ref-linear-8-u}.\\

Theorem \ref{thm:integrality-one}$^\prime$ reduces the $(J+2)^J$ many constraints \eqref{ref-linear-3-u} in the reformulation of (DRNS) to $(J^2 + J + 1) = O(J^2)$ many. This yields the following MILP reformulation of (DRNS), which extends Proposition \ref{prop:one-pool-ref}.\\
\noindent\textbf{Proposition \ref{prop:one-pool-ref}$^\prime$}
Under Assumption \ref{assump:technical} and Structure \ref{assump:one-pool}, the (DRNS) model \eqref{drns} yields the same optimal objective value as the following MILP:
\begin{subequations}
	\begin{align}
	Z^{\star}_{1} := \min \ & \ \theta + c^{\mbox{\tiny y}}_1 y^{\tinyl}_1 + g_1(y^{\tinyl}_1) \lambda_1 + \sum_{\ell=1}^{y^{\mbox{\tinyu}}_1 - y^{\mbox{\tinyl}}_1} \bigl(c^{\mbox{\tiny y}}_1 v_{1\ell}  + \delta_{1\ell} \nu_{1\ell} \bigl) \nonumber \\
	& \ + \sum_{j=1}^J \Biggl[ \sum_{q=1}^Q \mu_{jq} \rho_{jq} + c^{\mbox{\tiny w}}_j w^{\tinyl}_j + f_j(w^{\tinyl}_j) \gamma_j + \sum_{\ell=1}^{w^{\mbox{\tinyu}}_j - w^{\mbox{\tinyl}}_j} \bigl( c^{\mbox{\tiny w}}_j u_{j\ell} + \Delta_{j\ell} \varphi_{j\ell} \bigr) \Biggr] \nonumber \\
	\mbox{s.t.} \ & \ \mbox{\eqref{ref-note-36}--\eqref{ref-note-33}}, \ \mbox{\eqref{ref-linear-4-u}}, \nonumber \\
	& \theta \geq \phi^{\mbox{\tiny e}}_1 + \sum_{j=1}^J (\zeta^{\mbox{\tiny e}}_j + \eta^{\mbox{\tiny e}}_j), \nonumber \\
	& \theta \geq \sum_{j=1}^J (\zeta^{\mbox{\tiny e}}_j + \eta^{\mbox{\tiny e}}_j) + \phi^{\mbox{\tiny x}}_{1k^*} + (- x^{\text{\tiny U}}_{j^*} \kappa^{\mbox{\tiny x}}_{j^* k^*} +\zeta^{\mbox{\tiny x}}_{j^* k^*} + \eta^{\mbox{\tiny x}}_{j^* k^*}) - (\zeta^{\mbox{\tiny e}}_{j^*} + \eta^{\mbox{\tiny e}}_{j^*}) + \sum_{j \in  [J] \setminus \{j^*\}} \chi_{jk^*}, \nonumber \\
	& \hspace{7mm}   \forall j^* \in [J], \forall k^* \in [J+1], \nonumber \\
	& \chi_{jk^*} \geq  (\zeta^{\mbox{\tiny x}}_{jk} + \eta^{\mbox{\tiny x}}_{jk}) - (\zeta^{\mbox{\tiny e}}_j + \eta^{\mbox{\tiny e}}_j), \ \forall j \in [J], \forall k^* \in [J+1], \forall k \in [ k^* ], \nonumber \\ 
	& \hspace{-0.15cm} \left.\begin{array}{l}
	\zeta^{\mbox{\tiny e}}_j \geq  - \gamma_j w^{\tinyl}_j - \sum_{k=1}^{w^{\mbox{\tinyu}}_j - w^{\mbox{\tinyl}}_j}  \varphi_{jk} , \ \ \ \zeta^{\mbox{\tiny e}}_j \geq 0, \\[0.1cm]
	\eta^{\mbox{\tiny e}}_j \geq  - \sum_{q=1}^Q\tilde{d}^q_j\rho_{jq}, \ \forall \tilde{d}_j \in [d^{\tinyl}_j, d^{\tinyu}_j]_{\mathbb{Z}},
	\end{array}\right\} \ \ \forall j \in [J], \nonumber \\ 
	& \hspace{-0.15cm} \left.\begin{array}{l}
	\kappa^{\mbox{\tiny x}}_{jk} \geq c^{\text{\tiny x}}_k - c^{\text{\tiny x}}_j, \ \ \ \kappa^{\mbox{\tiny x}}_{jk} \geq 0, \\ [0.1cm]
	\zeta^{\mbox{\tiny x}}_{jk} \geq (- c^{\mbox{\tiny x}}_k - \gamma_j) w^{\tinyl}_j - \sum_{k=1}^{w^{\mbox{\tinyu}}_j - w^{\mbox{\tinyl}}_j} ( c^{\mbox{\tiny x}}_k u_{jk} + \varphi_{jk}), \ \ \ \zeta^{\mbox{\tiny x}}_{jk} \geq 0, \\[0.1cm]
	\eta^{\mbox{\tiny x}}_{jk} \geq c^{\mbox{\tiny x}}_k \tilde{d}_j - \sum_{q=1}^Q\tilde{d}^q_j\rho_{jq}, \ \forall \tilde{d}_j \in [d^{\tinyl}_j, d^{\tinyu}_j]_{\mathbb{Z}}, \\[0.1cm]
	\end{array}\right\} \ \ \forall j \in [J], \forall k \in [J+1], \nonumber \\ 
	& \ \phi^{\mbox{\tiny e}}_1 \geq  - \lambda_1 y^{\tinyl}_1 - \sum_{\ell=1}^{y^{\mbox{\tinyu}}_1 - y^{\mbox{\tinyl}}_1} \nu_{1\ell}, \ \ \ \phi^{\mbox{\tiny e}}_1 \geq 0, \nonumber \\
	& \ \phi^{\mbox{\tiny x}}_{1k} \geq (- c^{\mbox{\tiny x}}_k - \lambda_1) y^{\tinyl}_1 - \sum_{\ell=1}^{y^{\mbox{\tinyu}}_1 - y^{\mbox{\tinyl}}_1} (c^{\mbox{\tiny x}}_k v_{1\ell} + \nu_{1\ell} ), \ \ \phi^{\mbox{\tiny x}}_{1k} \geq 0, \ \forall k \in [J], \nonumber
	\end{align}
\end{subequations}

Under Structure~\ref{assump:disjoint-pools}, we derive a similar MILP reformulation of (DRNS), which reduces the $(J+2)^J$ many constraints \eqref{ref-linear-3-u} to $O(J^2 + I)$ many. 

Under Structure~\ref{assump:chained-pools}, we reformulate (DRNS) using the longest path on an acyclic network, which consists of $O(J^3)$ nodes and $O(J^4)$ arcs. This yields a MILP reformulation that reduces the $(J+2)^J$ many constraints \eqref{ref-linear-3-u} to $O(J^4)$ many.
\end{appendices}

\baselineskip=18pt
\bibliographystyle{plain}
\bibliography{rjiang}

\end{document}